\documentclass[11pt]{article}
\usepackage{amsmath}\usepackage{amssymb}\usepackage{amscd}
\usepackage{theorem}
\usepackage{xspace}
\input xy
\xyoption{all}

\newtheorem{lemma}{Lemma}\newtheorem{proposition}{Proposition}

\newcommand{\beq}{\begin{equation}}\newcommand{\eeq}{\end{equation}}
\newcommand{\ba}{\begin{array}}\newcommand{\ea}{\end{array}}
\newcommand{\beqa}{\begin{eqnarray}}\newcommand{\eeqa}{\end{eqnarray}}

\begin{document}

\title{
{\hfill \small{CPT-P49-2006}}\\[.5cm]
{\large\bf R-MATRICES IN RIME}}
\author{Oleg Ogievetsky\footnote{On leave of absence from P.N. Lebedev Physical Institute, 
Theoretical Department, Leninsky prospekt 53, 119991 Moscow, Russia}\\[.6em] 
{\it Centre de Physique Th\'eorique\footnote{Unit\'e Mixte de Recherche (UMR 6207) du CNRS et des 
Universit\'es Aix--Marseille I, Aix--Marseille II et du Sud Toulon -- Var; laboratoire affili\'e 
\`a la FRUMAM (FR 2291)}, Luminy, 13288 Marseille, France } \\[1.8em] 
Todor  Popov \\[.6em] 
\it Institute for Nuclear Research and Nuclear Energy, \\ 
\it Bulgarian Academy of Sciences, Sofia, 
BG-1784, Bulgaria}
\date{}\maketitle

\begin{abstract}\noindent
We replace the ice Ansatz on matrix solutions of the Yang--Baxter equation by a weaker condition 
which we call {\it rime}. Rime solutions include the standard Drinfeld--Jimbo $R$-matrix. 
Solutions of the Yang--Baxter equation within the rime Ansatz which are maximally different from 
the standard one we call {\it strict rime}. A strict rime non-unitary solution is parameterized by 
a projective vector $\vec{\phi}$. We show that in the finite dimension this solution transforms to 
the Cremmer--Gervais $R$-matrix by a change of basis with a matrix containing symmetric functions 
in the components of $\vec{\phi}$.  A strict unitary solution (the rime Ansatz is well adapted for 
taking a unitary limit) in the finite dimension is shown to be equivalent to a quantization of a 
classical "boundary" $r$-matrix of Gerstenhaber and Giaquinto. We analyze the structure of the 
elementary rime blocks and find, as a by-product, that all non-standard $R$-matrices of 
$GL(1|1)$-type can be uniformly described  in a rime form. 
 
\vskip .1cm
\noindent We discuss then connections of the classical rime solutions with the B\'ezout operators. 
The B\'ezout operators satisfy the (non-)homogeneous associative classical Yang--Baxter equation
which is related to the Rota--Baxter operators. We calculate the Rota--Baxter operators 
corresponding to the B\'ezout operators.

\vskip .1cm
\noindent We classify the rime Poisson brackets: they form a 3-dimensional pencil. A normal form 
of each individual member of the pencil depends on the discriminant of a certain quadratic 
polynomial. We also classify orderable quadratic rime associative algebras 
 
\vskip .1cm
\noindent For the standard Drinfeld--Jimbo solution, there is a choice of the multiparameters, 
for which it can be non-trivially rimed. However, not every Belavin--Drinfeld triple admits a 
choice of the multiparameters for which it can be rimed. We give a minimal example. 
\end{abstract}
 
\newpage
\tableofcontents 

\section{From ice to rime}

A well known class of solutions $\hat{R}\in\,$End$\, (V\otimes V)$, $V$ is a vector space,
of the Yang--Baxter equation $Y\!\! B(\hat{R})=0$, where
\beq Y\!\! B(\hat{R}):=(\hat{R} \otimes 1 \!\! 1) (1 \!\! 1 \otimes \hat{R}) (\hat{R} \otimes 1 \!\! 1)-
(1 \!\! 1 \otimes \hat{R}) (\hat{R} \otimes 1 \!\! 1)(1 \!\! 1 \otimes \hat{R})\ ,\eeq
is characterized by the so called {\it ice condition} (see lectures \cite{QS} for details)
which says that $\hat{R}^{i j }_{k l}$ can be different from zero
only if the set of the upper and the set of the lower indices coincide,
\beq \hat{R}^{ij}_{kl}\neq 0\qquad\Rightarrow\qquad\{ i,j\}\equiv\{ k,l\} .\eeq
We introduce the {\it "rime"} Ansatz, relaxing the ice condition: the entry $\hat{R}^{i j }_{k l}$ 
can be different from zero if the set of the lower indices is a subset of the set of the upper 
indices,
\beq \hat{R}^{ij}_{kl}\neq 0\qquad\Rightarrow\qquad\{ k,l\}\subset\{ i,j\}\ . \eeq
Matrices for which it holds will be referred to as {\it ``rime''} matrices. Figuratively,
in the rime, in contrast to the ice, situation, putting an apple and a banana in a fridge, there 
is a non-zero amplitude to find next morning two apples instead (but never an apple and an 
orange). 
 
\vskip .2cm
The Yang--Baxter equation for a matrix $\hat{R}$ is equivalent to the equality of two different 
reorderings of $x^iy^jz^k$, using $x^iy^j=\hat{R}^{ij}_{kl}y^kx^l$, 
$x^iz^j=\hat{R}^{ij}_{kl}z^kx^l$ and $y^iz^j=\hat{R}^{ij}_{kl}z^ky^l$, to the form
$z^\bullet y^\bullet x^\bullet$. One of advantages of the rime Ansatz is that only indices $i,j$ 
and $k$ appear in the latter expression. Another advantage is that for fixed values of $i$ and 
$j$, the elements $x^\bullet$ and $y^\bullet$ with these values of indices form a subsystem.

\vskip .2cm
A rime $R$-matrix has the following structure
\beq\hat{R}_{kl}^{ij}=\alpha_{ij}\delta^i_l\delta^j_k +\beta_{ij}\delta^i_k\delta^j_l+\gamma_{ij} 
\delta^i_k\delta^i_l+\gamma'_{ij}\delta^j_k\delta^j_l\ \qquad\mbox{(no summation)}
\ .\label{rice}\eeq
To avoid redundancy, fix $\beta_{ii}=0$, $\gamma_{ii}=0=\gamma'_{ii}$. We denote by $\alpha_i$ the 
diagonal elements $\hat{R}_{ii}^{ii}$, $\alpha_i=\alpha_{i i }$. Throughout the text we shall 
assume that the matrix $\hat{R}$ is invertible which, in particular, implies that $\alpha_i\neq 0$ 
for all $i$. 

\vskip .2cm
The order of growth of the number of unknowns in the Yang--Baxter system for a rime matrix is 
$n^2$, where $n=\,$dim$\, V$.

\vskip .2cm
Arbitrary permutations and rescalings of coordinates preserve the rime condition.

\vskip .2cm
The ice and rime matrices are made of $4 \times 4$ elementary building blocks, respectively,
\beq\label{fri}\hat{R}^{ice}=\left(\ba{cccc}\alpha_1&0&0&0\\0&\beta_{12}&\alpha_{12}&0\\
0&\alpha_{21}&\beta_{21}&0\\0&0&0&\alpha_2\ea\right)\qquad {\mathrm{and}}\qquad 
\hat{R}^{rime}=\left(\ba{cccc}\alpha_1&0&0&0\\ \gamma_{12}&\beta_{12}&\alpha_{12}&\gamma'_{12}\\
\gamma'_{21}&\alpha_{21}&\beta_{21}&\gamma_{21}\\0&0&0&\alpha_2\ea\right)\ .\eeq

In appendix B we analyze the structure of the $4 \times 4$ rime blocks.

\vskip .2cm
We call a rime matrix {\it strict} if  $\alpha_{ij}\gamma_{ij} \neq 0$ $\forall\ i$ and $j$,
$i\neq j$. Note that strict rime matrices are necessarily not ice.

\begin{proposition}$\!\!\!${\bf .} Let $\hat{R}$ be a rime matrix (\ref{rice}). Then $\hat{R}$ 
is a solution of the  Yang--Baxter equation if it is of the form 
\beq\hat{R}_{k l }^{i j}=(1 - \beta_{j i} )\delta^i_l \delta^j_k +\beta_{i j} \delta^i_k \delta^j_l +\gamma_{i j} \delta^i_k \delta^i_l - \gamma_{  j i } \delta^j_k \delta^j_l\ ,
\label{mat}\eeq
where  $\beta_{ij}$ and $\gamma_{ij}$ satisfy the system
\begin{eqnarray}
\beta_{ij}\beta_{ji}&=&\gamma_{ji}\gamma_{ij}\ ,\label{bg}\\
\beta_{ij}+\beta_{ji}&=&\beta_{jk}+\beta_{kj}=:\beta\ ,\label{beta}\\
\beta_{ij}\beta_{jk}&=&(\beta_{jk}-\beta_{ji})\beta_{ik}=
(\beta_{ij}-\beta_{kj})\beta_{ik}\ ,\label{betaeq}\label{betab}\\
\gamma_{ij}\gamma_{jk}&=&(\beta_{ji}-\beta_{jk})\gamma_{ik}=
(\beta_{kj}-\beta_{ij})\gamma_{ik}\label{gammaeq}\ .
\end{eqnarray} \end{proposition}

\noindent
{\bf Proof.} The Yang--Baxter system of equations $Y\!\! B(\hat{R})^{i j k}_{a b c }=0$
for a rime matrix is given in the appendix A. 
The subset (\ref{yb1}) - (\ref{yb3}) together with its image under the involution 
(\ref{invo}) reads
\begin{eqnarray}
\label{eII}\alpha_{ij}\gamma'_{ij} (\gamma_{ij}+\gamma'_{ji})&=&0 
\quad=\alpha_{ij}\gamma_{ij}(\gamma_{ji}+\gamma'_{ij})\ ,\\	  
\label{eIII}\alpha_{ij}(\beta_{ij}\beta_{ji}+\gamma_{ij}\gamma'_{ij})&=&0\quad =
\alpha_{ij}(\beta_{ij}\beta_{ji}-\gamma_{ij}\gamma_{ji})\ ,\\
\label{eIV} \alpha_{ij}\gamma'_{ij}(\alpha_{ij}+\beta_{ji}-\alpha_i)&=&0 
\quad =\alpha_{ij}\gamma_{ij}(\alpha_{ij}+\beta_{ji}-\alpha_j)\ ,\\	
\label{eV}\alpha_{ij}\gamma'_{ij}(\alpha_{ji}+\beta_{ij}-\alpha_i)&=&0\quad =
\alpha_{ij}\gamma_{ij}(\alpha_{ji}+\beta_{ij}-\alpha_j)\ .\end{eqnarray}
These equations are implied by (and, in the strict rime situation, are equivalent to) 
the following system
\beqa\label{subst}
&\gamma'_{i j} = - \gamma_{j i} \ ,\qquad
\alpha_{i j} + \beta_{j i} =\alpha_i \ ,\qquad \alpha_{j i} +\beta_{i j} =\alpha_i\ ,&\\
& \label{subst2} \beta_{i j}  \beta_{ j i} = \gamma_{j i}\gamma_{i j}\ .&
\eeqa
One checks that other equations $Y\!\! B(\hat{R})^{i j k}_{a b c }=0$, for which
two indices among $\{i,j,k \}$ are different, follow from (\ref{subst}) and (\ref{subst2}).
The last two equations from (\ref{subst}) imply $\alpha_i =\alpha_j$ for all $i$ and $j$.
As an overall rescaling of a solution of the Yang--Baxter equation by a constant is again 
a solution of the Yang--Baxter equation, we can, without loss of generality, set it to one, 
\beq\alpha_i=1\ .\label{subst3}\eeq
Eqs. (\ref{subst}) and (\ref{subst3}) yield the form (\ref{mat}) of the matrix $\hat{R}$
and eq. (\ref{bg}).

Using (\ref{subst}), we rewrite the subset (\ref{ee1}) - (\ref{ee3}) together with its 
image under the involution (\ref{invo}) in the form
\begin{eqnarray}
&(\beta_{ij} +\beta_{ji}-\beta_{ik}-\beta_{ki})\gamma_{ij}\gamma_{ik}=0\ ,&\\
&\alpha_{ij}(\beta_{ij}\beta_{jk}+\beta_{ik}\beta_{ji}-\beta_{ik}\beta_{jk})=0 
=\alpha_{ji}(\beta_{ji}\beta_{kj}+\beta_{ki}\beta_{ij}-\beta_{ki}\beta_{kj})\ ,&\\
&\alpha_{ij}(\gamma_{ij}\gamma_{jk}+\gamma_{ik}(\beta_{jk}-\beta_{ji}))=0 
=\alpha_{ji}(\gamma_{ji}\gamma_{kj}+\gamma_{ki}(\beta_{kj}-\beta_{ij}))\ .&
\end{eqnarray}
These equations are implied by (and, in the strict rime situation, are equivalent to)
eqs. (\ref{beta}), (\ref{betaeq}) and (\ref{gammaeq}). 
One checks that other equations $Y\!\! B(\hat{R})^{i j k}_{a b c }=0$ with three 
different indices $\{i,j,k \}$ follow from the system (\ref{bg})--(\ref{gammaeq}).
The proof is finished. \hfill $\Box$

\begin{lemma}$\!\!\!${\bf .} \label{Heck}
The rime Yang--Baxter solution $R$ (\ref{mat}) is of Hecke type,
\beq\label{Hecke}\hat{R}^2= \beta \hat{R} + (1-\beta) 1\!\!1 \otimes 1\!\!1\ .\eeq
Moreover, when $\beta\neq2$, $R$ is of $GL$-type: it has two eigenvalues $1$ and $\beta -1$ 
with multiplicities $\frac{n(n+1)}{2}$ and $\frac{n(n-1)}{2}$, respectively. When $\beta=2$ 
the matrix $\hat{R}$ has a nontrivial Jordanian structure.
\end{lemma}

\noindent
{\bf Proof.} In view of the block structure of rime matrices it is enough to check the 
Hecke relation (\ref{Hecke}) 
for one elementary ($4 \times 4$) block which follows from (\ref{bg}) and (\ref{beta}). 
When $\beta \neq 2$ the multiplicities $m_1$ and $m_{\beta -1}$ are solutions 
of the system
\beq m_1 +m_{\beta -1} =n^2\ ,\qquad 
m_1+(\beta -1)m_{\beta -1}=n+\frac{n(n-1)}{2}\beta \ \  (\equiv Tr \hat{R})\ .\eeq
When $\beta=2$ the matrix $\hat{R}$ has only one eigenvalue $1$ but 
$\hat{R}\neq 1\!\!1 \otimes 1\!\!1$ due to (\ref{bg}) and (\ref{beta}). \hfill $\Box$	 

\vskip .2cm
Unitary solutions, $\hat{R}^2=1\!\!1\otimes 1\!\!1$, are singled out by the value of
the parameter $\beta=0$.
  
\begin{lemma}$\!\!\!${\bf .} A strict  rime Yang--Baxter solution $R$ (\ref{mat}) can be brought 
to a rime matrix
\beq\label{R}
\hat{R}_{k l }^{i j}=(1 - \beta_{j i} )\delta^i_l \delta^j_k + 
\beta_{i j} \delta^i_k \delta^j_l - 
\beta_{i j} \delta^i_k \delta^i_l + \beta_{ j i} \delta^j_k \delta^j_l\ ,
\eeq
that is, to a solution (\ref{mat}) with $\gamma_{ij}=-\beta_{ij}$, by a change of basis.
\end{lemma} 

\noindent
{\bf Proof.} The strict rime condition $\alpha_{ij}\gamma_{ij}\neq 0$ implies 
 $\beta_{ij}\beta_{ji}\neq 0$ in view of (\ref{bg}). Thus for a strict rime $R$-matrix all 
$\beta_{ij}$ and $\gamma_{ij}$ are nonvanishing. The ratio of eqs. (\ref{betaeq}) and (\ref{gammaeq}) is well-defined and it follows from eqs. (\ref{bg}) and (\ref{beta}) that
\beq
\frac{\gamma_{i j}  \gamma_{j k}}{\beta_{i j}  \beta_{j k}}= - 
\frac{( \beta_{ j i} -  \beta_{ j k})\gamma_{i k}}{( \beta_{ j i} -  \beta_{ j k}) \beta_{ i k}}
=-\frac{\gamma_{i k}}{ \beta_{ i k}}\qquad \qquad 
\frac{\gamma_{i j} \gamma_{j i}}{\beta_{i j}\beta_{j i }}=1\ ,\label{cocy}
\eeq
or
\beq\xi_{ij} \xi_{jk}=\xi_{ik} \qquad \qquad \xi_{ij} \xi_{ji}=1 \ ,\label{exi}\eeq 
where $\xi_{ij}=-\frac{\gamma_{i j}}{ \beta_{ i j}}$.
Eq. (\ref{exi}) is solved by $\xi_{ij}= \frac{d_i}{d_j}$ with $d_i\neq 0$, $i=1, \ldots,n$, hence
$\beta$'s and $\gamma$'s are related by
\beq\gamma_{ij}=-  \frac{d_i}{d_j} \beta_{ij}\ .\eeq

A change of basis with a matrix $D$, 
\beq\hat{R}\longmapsto (D \otimes D) \, \hat{R} \,\, (D^{-1} \otimes D^{-1}) \ ,\eeq 
where $D^i_j= d_j \delta^i_j$, transforms $R$ to the form (\ref{R}).\hfill $\Box$

\vskip .3cm
Under the strict rime condition, the Yang--Baxter system of equations (see appendix A) 
reduces to eqs. (\ref{beta}) and (\ref{betab}). However, the matrix (\ref{R}), where the 
parameters $\beta_{ij}$ are subject to eqs. (\ref{beta}) and (\ref{betab}), is a solution of the 
Yang--Baxter equation without a strict rime assumption. 

\vskip .3cm
\noindent{\bf Remark.} Right and left even quantum spaces are defined by, respectively,
\beq R^{ij}_{kl}x^kx^l=x^i x^j\ ,\ x_jx_iR^{ij}_{kl}=x_lx_k\ ;\eeq
right and left odd quantum spaces are defined by, respectively,
\beq R^{ij}_{kl}\xi^k\xi^l=(\beta -1)\xi^i\xi^j\ ,\ \xi_j\xi_iR^{ij}_{kl}=(\beta -1)\xi_l\xi_k 
\ .\eeq
Assume that $\beta \neq 2$.
The left even space is classical\footnote{Let $\hat{R}$ be a rime $R$-matrix (not necessarily 
strict). When $\beta \neq 2$, the following statement holds. If (i) the left even space is 
classical (which implies that $\gamma_{ij}'=-\gamma_{ji}$, $\alpha_{ij}+\beta_{ji}=1$ and
$\alpha_i=1$ in our normalization) and (ii) the $R$-matrix is Hecke (which implies that 
$\beta_{ij}+\beta_{ji}=\beta$) then the system of equations from the appendix A again reduces
to (\ref{bg}), (\ref{betaeq}) and (\ref{gammaeq}) as in the strict rime situation.}
as well as the right odd space
\beq [ x_i,x_j ] =0\ ,\qquad [ \xi^i,\xi^j]_{_+} =0\ ,\eeq
where $ [ \, , ] $ and $[\, ,]_{_+}$ stand for the commutator and the anti-commutator.
The relations for the right even space are
\beqa \label{qp}
 [ x^i,x^j ] +( \beta_{ij}  x^i+\beta_{ji} x^j)(x^i-x^j) =0\ ;\eeqa
the relations for the left odd space read
\beqa
(2-\beta )\xi_i^2+ \xi_i \rho + (1-\beta) \rho \xi_i & =& 0  \ ,\\
\ [ \xi_i,\xi_j]_{_+} -\beta_{ij}\xi_i\xi_j -\beta_{ji}\xi_j\xi_i&=&0\ ,
\qquad i\neq j\ ,\eeqa
where $\rho=\sum_j \xi_j$. 

\section{Rime Yang--Baxter solutions\vspace{.25cm}} 

In this section we  solve  eqs. (\ref{beta}) and (\ref{betab}) thus obtaining explicitly
rime Yang--Baxter solutions.  

\subsection{Non-unitary rime R-matrices}
 
\begin{proposition}$\!\!\!${\bf .} \label{prop:rnu}
The non-unitary  strict rime  Yang--Baxter  solutions  (\ref{R}) with a parameter 
$\beta =\beta_{ji}+\beta_{ij}\neq 0$ are parameterized by a point $\phi\in \mathbb P \mathbb C^n$ 
in a projective space, $\phi=(\phi_1: \phi_2: \ldots :\phi_n)$, such that $\phi_i\neq 0$ for 
all $i$ and $\phi_i\neq \phi_j$ for all $i$ and $j$, $i\neq j$. These solutions are given by
\beq \label{bphi}\beta_{ij}=\frac{\beta\phi_i}{\phi_i-\phi_j}\ .\eeq\end{proposition}
 
\noindent{\bf Proof.} Taking the ratio of the following pairs of equations from (\ref{betab}) 
\beq\beta_{ij}\beta_{jk}=(\beta_{jk}-\beta_{ji})\beta_{i k}\ ,
\qquad\beta_{kj}\beta_{ji}=(\beta_{ji}-\beta_{jk})\beta_{k i}\eeq
we find that quantities $\eta_{i j}=-\beta_{ij}/\beta_{ji}$ verify equations
\beq\eta_{i j} \eta_{ j k}= \eta_{ i k}\ , \qquad\eta_{i j} \eta_{ j i}=1\ ,\eeq
whose solutions are $\eta_{ij}=\phi_i/\phi_j$ for some constants $\phi_i\neq 0$, $i=1,\ldots,n$.

Substituting the relation $\beta_{ji}=-\displaystyle{\frac{\phi_j}{\phi_i}\beta_{ij}}$ into 
$\beta =\beta_{ij}+\beta_{ji}$, we obtain $\beta_{ij}-\displaystyle{\frac{\phi_j}{\phi_i}} \beta_{ij}=\beta$ which establishes (\ref{bphi}).\hfill $\Box$

\vskip .2cm\noindent {\bf Remark.} There is a different parameterization, $\beta_{ij}= 
-\displaystyle{\frac{\beta\phi_j}{\phi_i-\phi_j}}$, of strict rime solutions; it is related to the 
parameterization (\ref{bphi}) by $\phi_i\longmapsto (\phi_i)^{-1}$.

\vskip .2cm
A direct check shows that the condition $\phi_i\neq 0$ is not necessary:
the formula (\ref{bphi}) with $\phi_i\neq\phi_j$ for all $i$ and $j$, $i\neq j$, 
gives a rime solution of the Yang--Baxter equation. 
However when one of $\phi_i$ is 0, the matrix (\ref{R}) is no more strict.

\subsection{Unitary rime R-matrices}
 
For a unitary strict rime Yang--Baxter solution (\ref{R}), $\hat{R}^2=1\!\!1$, we have $\beta =0$, 
so $\beta_{ij}=-\beta_{ji}$.

\begin{proposition}$\!\!\!${\bf .}\label{prop: ru}
The unitary  strict rime  Yang--Baxter  solutions (\ref{R}) are parameterized by a vector 
$(\mu_1,\ldots ,\mu_n)$ such that $\mu_i \neq \mu_j$,
\beq\label{unitary}\beta_{ij} =\frac{1}{\mu_i - \mu_j}\ .\eeq
\end{proposition}

\noindent{\bf Proof.} Since $\beta_{ij}=-\beta_{ji}$ we can rewrite  
$\beta_{ij}\beta_{j k} =(\beta_{jk}-\beta_{ji})\beta_{ik}$
as $\beta_{ij}\beta_{jk}=(\beta_{ij}+\beta_{jk})\beta_{ik}$ or
\beq\frac{1}{\beta_{ik}}=\frac{1}{\beta_{ij}}+\frac{1}{ \beta_{jk}}.\label{iter}\eeq
These equations are solved by
\beq\frac{1}{\beta_{ij}} =\mu_i - \mu_j\ ,\eeq
which is equivalent to (\ref{unitary}).\hfill$\Box$ 

\vskip .3cm 
\noindent {\bf Remark.} 
The unitary $R$-matrices of Proposition \ref{prop: ru} can be obtained as a limit 
$\beta \rightarrow 0$ of the non-unitary $R$-matrices of Proposition \ref{prop:rnu}.
Indeed, for the following expansion of the parameters $\phi_i$ in
the ``small'' parameter $\beta$, 
\beq\phi_i = 1 + \beta \mu_i  + o(\beta)\ , \eeq
the expression (\ref{bphi}) has a limit (\ref{unitary}), 
\beq\beta_{ij} =\frac{\beta (1+\beta\mu_i+o(\beta))}{\beta\mu_i-\beta\mu_j+o(\beta)}\qquad
\stackrel{\beta\rightarrow 0}{\longrightarrow}\qquad 
\beta_{ij}=\frac{1}{\mu_i-\mu_j}\ .\label{liu}\eeq

\subsection{Properties}

\paragraph{1.}
Denote the $R$-matrix (\ref{R}) with $\beta_{ij}$ as in (\ref{bphi}) by $\hat{R}(\vec{\phi}\, )$. 
Let $\hat{R}_{21}=P\hat{R}_{12}P$, where $P$ is the permutation operator. 
Then the following holds:
\beq \hat{R}_{21}(\vec{\phi}\, )=F^{-1}\otimes F^{-1}\ 
R_{12}(\vec{\phi}{}^{\scriptscriptstyle{-1}})\ F\otimes F\ ,\eeq
where $F={\mathrm{diag}}(\phi_1,\phi_2,\dots ,\phi_n)$ and $\vec{\phi}{}^{\scriptscriptstyle{-1}}$ 
is a vector with components $\phi_i^{-1}$.

\vskip .2cm
Denote the $R$-matrix (\ref{R}) with $\beta_{ij}$ as in (\ref{unitary}) by $\hat{R}(\vec{\mu}\, )$.
Then the following holds:
\beq \hat{R}_{21}(\vec{\mu} )=\hat{R}_{12}(-\vec{\mu} )\ .\eeq

\paragraph{2.} The $R$-matrix (\ref{R}) is skew invertible in the sense that there exists an
operator $\hat{\Psi}_R$, which satisfies (see, {\em e.g.} \cite{QS})
\beq {\mathrm{Tr}}_2 ( \hat{R}_{12} (\hat{\Psi}_R)_{23}) = P_{13}\ .\eeq
The matrices of the left and right quantum traces (that is, the left and right traces of the skew inverse $\hat{\Psi}_R$), $(Q_R)_1={\mathrm{Tr}}_2 ((\hat{\Psi}_R)_{12})$
and $(\tilde{Q}_R)_2={\mathrm{Tr}}_1((\hat{\Psi}_R)_{12})$, are given by the formulas
\beqa\label{qtq1}
(Q_R)^k_j=-\beta_{jk}{\displaystyle\prod_{l:\ l\neq k}(1-\beta_{jl})}\ ,\ \ 
k \neq j\ ,&{\mathrm{and}}& (Q_R)^j_j={\displaystyle\prod_l(1-\beta_{jl})}\ ;\\[1em]
(\tilde{Q}_R)^k_j=\ \ \beta_{jk}{\displaystyle\prod_{l:\ l\neq k}(1-\beta_{lj})}\ ,\ \ k\neq j\ ,
&{\mathrm{and}}&(\tilde{Q}_R)^j_j={\displaystyle\prod_l(1-\beta_{lj})}\ .\label{qtq2}\eeqa
The matrices $Q_R$ and $\tilde{Q}_R$ satisfy $Q_R\tilde{Q}_R=(1-\beta )^{n-1}1\!\! 1$.

\vskip .2cm
For (\ref{bphi}), one has $Spec\ Q_R=Spec\ \tilde{Q}_R=\{ (1-\beta)^a ,\ a=0,\dots ,n-1\}$. The
eigenvector $w_a(\vec{\phi})$ of the matrix $Q_R$ with the eigenvalue $(1-\beta)^{n-1-a}$ 
coincides with the eigenvector of the matrix $\tilde{Q}_R$ with the eigenvalue $(1-\beta)^{a}$. 
One has $(w_a(\vec{\phi}))^j= e_{a}^{\hat{j}}(\vec{\phi})$, where $e_i^{\hat{j}}(\vec{\phi})$ is 
the $i$-th elementary symmetric function of $(\phi_1,\phi_2,\dots ,\phi_n )$ with $\phi_j$ omitted.

\vskip .2cm
For (\ref{unitary}), the Jordanian form of the matrix $Q_R$, as well as of $\tilde{Q}_R$, is 
non trivial: it is a single block. In the basis $\{ w_i(\vec{\mu})\}$, $i=0,1,2,\dots ,n-1$, where 
$(w_i(\vec{\mu}))^j=e_{i}^{\hat{j}}(\vec{\mu})$, one has 
\beq Q_R\,  w_i(\vec{\mu}) 
=\sum_{s=0}^{i}\left(\begin{array}{c}n-1-s\\ i -s \end{array}\right) w_s(\vec{\mu})\ .\eeq

\paragraph{3.} For an $R$-matrix $\hat{R}$, the group of invertible matrices $Y$ satisfying
\beq \hat{R}_{12}Y_1Y_2=Y_1Y_2\hat{R}_{12}\label{ingr}\eeq 
form the invariance group $G_R$ of $\hat{R}$. The matrices $Q_R$ and 
$\tilde{Q}_R$ belong to the invariance group as well as the matrices proportional to the 
identity matrix. One can write down formulas for the group $G_R$ for 
a rime $R$-matrix (\ref{R}) uniformly in terms of $\beta_{ij}$ as in (\ref{qtq1}) and (\ref{qtq2}) 
but the properties are different in the non-unitary and unitary cases and we describe them 
separately.
 
\vskip .2cm
{\bf 3a.} The invariance group $G_{R(\vec{\phi}\, )}$ for the $R$-matrix $\hat{R}(\vec{\phi}\, )$ 
is 2-parametric. It consists of matrices $Y(u,v)$, $u,v\neq 0$, where
\beq Y(u,v)^j_j=\prod_{l:l\neq j}\frac{u\phi_j-v\phi_l}{\phi_j-\phi_l}\ \ \ {\mathrm{and}}
\ \ \ Y(u,v)^i_j=\frac{(u-v)\phi_j}{\phi_j-\phi_i}
\prod_{l:l\neq i,j}\frac{u\phi_j-v\phi_l}{\phi_j-\phi_l}\ ,\ i\neq j\ .\label{inr1}\eeq

One has
\beq Q_{R(\vec{\phi}\, )}=Y(1-\beta ,1)\ \ ,\ \ \tilde{Q}_{R(\vec{\phi}\, )}=Y(1,1-\beta )\ .\eeq
The composition law is the component-wise multiplication of the parameters $\{ u,v\}$,
\beq Y(u_1,v_1)Y(u_2,v_2)=Y(u_1u_2,v_1v_2)\ .\eeq
The point $u=v=1$ corresponds to the identity matrix, $Y(1,1)=1\!\!1$; the determinant of $Y(u,v)$ 
is $(uv)^{n(n-1)/2}$; $u=v$ corresponds to global rescalings; the connected component of unity 
of the subgroup $SG_{R(\vec{\phi}\, )}$ consisting of matrices with determinant 1 is $uv=1$; the 
generator $\eta$ of the connected component of unity of the subgroup $SG_{R(\vec{\phi}\, )}$ 
is traceless and reads
\beq \eta^i_j=\frac{\phi_j}{\phi_j-\phi_i}\ ,\ i\neq j\ ,\ \ {\mathrm{and}}\ \ 
\eta^j_j=-\frac{n-1}{2}+\sum_{l:l\neq j}\frac{1}{\phi_j-\phi_l}\ .\label{inr2}\eeq

\vskip .1cm
{\bf 3b.} For the $R$-matrix $\hat{R}(\vec{\mu}\, )$, the group $SG_{R(\vec{\mu}\, )}$, 
consisting of matrices with determinant 1 is 1-parametric as well. It is formed by matrices 
$Y^{(0)}(a)$, where
\beq Y^{(0)}(a)^j_j=\prod_{l:l\neq j}(1+\frac{a}{\mu_j-\mu_l})\ \ {\mathrm{and}}\ \ Y^{(0)}(a)^i_j
=\frac{a}{\mu_j-\mu_i}\prod_{l:l\neq i,j}(1+\frac{a}{\mu_j-\mu_l})\ ,\ i\neq j\ .\label{inr3}\eeq
The expression (\ref{inr3}) can be obtained by taking a limit of (\ref{inr1}), similarly to 
(\ref{liu}) and letting additionally $u=1+a\beta /2 +o(\beta )$ and $v=1-a\beta /2 +o(\beta )$.

\vskip .2cm
One has
\beq Q_{R(\vec{\mu}\, )}=Y^{(0)}(-1)\ \ ,\ \ \tilde{Q}_{R(\vec{\mu}\, )}=Y^{(0)}(1)\ .\eeq

The composition law is $Y^{(0)}(a_1)Y^{(0)}(a_2)=Y^{(0)}(a_1+a_2)$.

\vskip .2cm
The point $a=0$ in (\ref{inr3}) corresponds to the identity matrix, $Y^{(0)}(0)=1\!\!1$; the 
generator $\eta^{(0)}$ of the invariance group $SG_{R(\vec{\mu}\, )}$ is
\beq (\eta^{(0)})^i_j=\frac{1}{\mu_j-\mu_i}\ ,\ i\neq j\ ,\ \ {\mathrm{and}}\ \ 
(\eta^{(0)})^j_j=\sum_{l:l\neq j}\frac{1}{\mu_j-\mu_l}\ .\label{inr4}\eeq

\section{Rime and Cremmer--Gervais R-matrices}\label{sectionCG}

The Cremmer--Gervais $R$-matrix arises in the exchange relations of the chiral vertex operators
in the non-linearly $W$-extended Virasoro algebra \cite{CG}. The  Cremmer--Gervais solution  
\cite{CG}  of the Yang--Baxter equation in its general two-parametric form reads (see, {\em e.g.} 
\cite{Ho}; we use a rescaled matrix with eigenvalues $1$ and $ -q^{-2}$) 
\beq\label{CG}
(\hat{R}_{CG,p})^{ij}_{kl} = q^{-2\theta_{ij}} p^{i-j} \delta^{i}_{l} \delta^{j}_{k}
+ (1 - q^{-2}) \sum_{s:\, i \leq s < j} p^{i-s} \delta^{s}_{k} \delta^{i + j -s }_{l} 
- (1 - q^{-2}) \sum_{s:\, j < s < i} p^{i-s}   \delta^{s}_{k} \delta^{i + j -s }_l\ ,\eeq
where  $\theta_{ij}$ is the step function ($\theta_{ij} =1$ when $i >j$ and $\theta_{ij} =0$ 
when $i \leq j$). 

\vskip .2cm
The parameter value $p=q^{2/n}$ specifies the $SL(n)$ Cremmer--Gervais $R$-matrix 
(its diagonal twist being the $GL(n)$ solution (\ref{CG})). The Cremmer--Gervais solution is a 
non-diagonal twist of the standard Drinfeld--Jimbo solution \cite{IO,ESS}. 

\vskip .2cm
Let $\hat{R}_{CG}:=\hat{R}_{CG,1}$, that is, the solution (\ref{CG}) with $p=1$. The matrix 
\beq D(p)^i_j=\delta^i_j p^{i-1}\label{invdcg}\eeq
with arbitrary $p$ satisfies $(\hat{R}_{CG})_{12}D(p)_1D(p)_2=D(p)_1D(p)_2 (\hat{R}_{CG})_{12}$. 
It was observed in \cite{EO} that if $\hat{R}_{12}D_1D_2=D_1D_2\hat{R}_{12}$ for some $R$-matrix 
$\hat{R}$ and operator $D$ then $D_1\hat{R}_{12}D^{-1}_1$ is again an $R$-matrix (this operation 
was also used in \cite{GIO} to partially change the statistics of ghosts in the super-symmetric 
situation). The two-parametric matrix $\hat{R}_{CG,p}$ (\ref{CG}) can be obtained from the 
Cremmer--Gervais matrix $\hat{R}_{CG}$ by this operation as well,
\beq (\hat{R}_{CG,p})_{12}= D(p)_1(\hat{R}_{CG})_{12}D(p)^{-1}_1\ .\label{chst}\eeq
 
Let $\hat{R}$ be the non-unitary rime matrix from Proposition \ref{prop:rnu} with 
$\phi_i\neq\phi_j$.

\begin{proposition}$\!\!\!${\bf .} \label{transem} The matrix $\hat{R}$ transforms into the 
Cremmer--Gervais solution $\hat{R}_{CG}$ 
\beq\label{change}\hat{R}=(X\otimes X)\,\hat{R}_{CG}\,\, (X^{-1}\otimes X^{-1})\eeq
by a change of basis with the invertible matrix 
\beq X^k_j = e_{j-1} \, (\phi_{1},\ldots, \hat{\phi}_{k}, \ldots,\phi_{n} )=:e_{j-1}^{\hat{k}}
\label{matX}\eeq 
whose inverse is
\beq
(X^{-1})^j_i=\frac{ (-1)^{j-1} \phi_i^{n-j}}
{\displaystyle{\prod_{k:k\neq i}} (\phi_{i}- \phi_{k})}\ . \label{transe}\eeq 
Here the hat over ${\phi}_{j}$ means that this entry is omitted in the expression and $e_i$ are 
the elementary symmetric polynomials 
$e_i(x_1,\ldots ,x_N)=\mathop{\sum}\limits_{s_1< \ldots <s_{i}}x_{s_1}x_{s_2}\ldots x_{s_{i}}.$
The projective parameters $(\phi_1: \phi_2: \ldots :\phi_n)$ are removed by the transformation 
$X$ and the only essential parameter $\beta$ in $\hat{R}$ is related to the parameter $q$ in 
$\hat{R}_{CG}$ by \beq q^{-2}=1- \beta\ .\eeq\end{proposition}
 
\noindent {\bf Proof.} Due to the Lagrange interpolation formula, the matrix, inverse to the 
Vandermonde matrix 
\beq ||V_k^j ||_{j,k=1}^n = \phi^{n-j}_k\quad\quad\quad
{\mathrm{is}}\quad\quad\quad (V^{-1})^k_j=
\frac{(-1)^{j-1}e_{j-1}^{\hat{k}}}{\displaystyle{\prod_{l:l\neq k}}(\phi_k-\phi_l)}\ .\eeq
The matrix $X$ (\ref{matX}) has the form $X=D\, V^{-1}\, d$, 
where $D^m_k = \delta^m_k {\prod_{l:l\neq k} (\phi_k -\phi_l)} $ and $d^i_j = (-1)^{j-1} 
\delta^i_j$ are diagonal $n\times n$ matrices. Thus, its inverse is
$X^{-1} = d^{-1}\, V \, D^{-1}$, which establishes (\ref{transe}).

We now prove the matrix identity (\ref{change}) in the form 
\beq\label{cha0}\hat{R}(X\otimes X)=(X\otimes X)\,\hat{R}_{CG}\ . \eeq
The substitution of the explicit form of the rime matrix $\hat{R}$ 
(\ref{R}) with $\beta_{ij}= \beta \phi_i /(\phi_{i}-\phi_j)$ and $\hat{R}_{CG}$ (\ref{CG})
reduces (\ref{cha0}) to a set of relations between the symmetrical polynomials
$e_{k-1}^{\hat{a}}$
\beq\label{cha}\sum_{a,b}\hat{R}^{ij}_{ab} e_{k-1}^{\hat{a}}e_{l-1}^{\hat{b}}=\sum_{a,b}
e_{a-1}^{\hat{i}}e_{b-1}^{\hat{j}}   (\hat{R}_{CG})^{ab}_{k l}\ .\eeq
There are two subcases: i) $i=j$ and ii) $i\neq j$.

\vskip .3cm
i) The left hand side of eq. (\ref{cha}) with $i=j$ is just $e_{k-1}^{\hat{i}}e_{l-1}^{\hat{i}}$ 
due to the rime condition. Eq. (\ref{cha}) is satisfied because of the symmetry relation 
$(\hat{R}_{CG})^{ab}_{kl}=\delta^a_k \delta^b_l + \delta^a_l \delta^b_k
-(\hat{R}_{CG})^{ba}_{kl}$. 

\vskip .3cm
ii) For $i \neq j$ eq. (\ref{cha}), where $q^{-2}=1- \beta$, reduces, after some algebraic manipulations, to
\beq\label{sectype}\begin{array}{l}
\displaystyle{\frac{1}{\phi_i - \phi_j}}(\phi_i e_{k-1}^{\hat{i}}-\phi_j e_{k-1}^{\hat{j}})
(e_{l-1}^{\hat{i}}- e_{l-1}^{\hat{j}} )\\[.5em]
\hspace{2cm} =\ \ \displaystyle{\sum_{s:\, s\geq\max (1,k-l+2)}^{}} \ \  
(e_{l+s-2}^{\hat{i}}e_{k-s}^{\hat{j}} -
e_{l+s-2}^{\hat{j}}e_{k-s}^{\hat{i}})\ , \quad  1\leq i,j,k,l \leq n\ .\end{array}\eeq
In fact, the sum in the right hand side goes till $s=\min (k,n+1-l)$ since $e_{r}^{\hat{j}}=0$
when $r\geq n-1$; moreover we can start the summation from $s=1$ 
because when $1<k-l+2$ the sum for $1\leq s \leq k-l+1$ is anti-symmetric under 
$s\longleftrightarrow k-l+2-s$ and thus vanishes.
 
To prove (\ref{sectype}) we write $e_r= e_r^{\hat{i}} +\phi_i e_{r-1}^{\hat{i}}$; therefore
$e_r^{\hat{i}}= e_r^{\hat{i}\hat{j}} +\phi_j e_{r-1}^{\hat{i}\hat{j}}$ and
$e_r^{\hat{j}}= e_r^{\hat{i}\hat{j}} +\phi_i e_{r-1}^{\hat{i}\hat{j}}$ and
eq. (\ref{sectype}) becomes
\beq -(\phi_i -\phi_j)e_{k-1}^{\hat{i}\hat{j}}e_{l-2}^{\hat{i}\hat{j}}=
(\phi_i - \phi_j) \sum_{s\geq 1}(e_{l+s-2}^{\hat{i}\hat{j}}e_{k-s-1}^{\hat{i}\hat{j}}
-e_{l+s-3}^{\hat{i}\hat{j}}e_{k-s}^{\hat{i}\hat{j}})\ .\eeq
The sum in the right hand side telescopes to the value of 
$(-e_{l+s-3}^{\hat{i}\hat{j}}e_{k-s}^{\hat{i}\hat{j}})$ at $s=1$, that is, to 
$(-e_{k-1}^{\hat{i}\hat{j}}e_{l-2}^{\hat{i}\hat{j}})$. The proof is complete. \hfill$\Box$ 

\vskip .2cm
It should be noted that the matrix $X=X(\vec{\phi}\,)$ does not depend on $q$. The change of the 
basis with the matrix $X(\vec{\phi'}\,)X(\vec{\phi}\,)^{-1}$ transforms the $R$-matrix 
$\hat{R}(\vec{\phi}\,)$ to $\hat{R}(\vec{\phi'}\,)$. We have
\beq (X(\vec{\phi'}\,)X(\vec{\phi}\,)^{-1})^i_j=\frac{1}{\phi_j-\phi_i'}\ 
\frac{{\displaystyle{\prod_k(\phi_j-\phi_k')}}}{{\displaystyle{\prod_{l:l\neq j}(\phi_j-\phi_l)}}}
\eeq

The structure of the matrices $X$ and $X^{-1}$ shows that when the dimension is infinite,
the $R$-matrices $\hat{R}_{CG,1}$ and $\hat{R}(\vec{\phi}\,)$ (as well as the $R$-matrices
$\hat{R}(\vec{\phi}\,)$ and $\hat{R}(\vec{\phi'}\,)$ for different $\phi$ and $\phi'$) are in 
general not equivalent. 

\vskip .3cm
The  right even quantum plane for the Cremmer--Gervais matrix $\hat{R}_{CG,1}$ 
is defined by the following equations 
\beq\label{qpcg}
y^i y^j = q^2 y^j y^i + (q^2 -1) (y^{i+1} y^{j-1}+ \ldots + y^{j-1}y^{i+1}), \qquad  i < j \ .
\eeq
If $i+1<j-1$, one uses the formula (\ref{qpcg}) recursively to get the ordering relations.

The change of basis with the matrix $X$,
\beq x^i = \sum_{j =1}^n e^{\hat{i}}_{j-1} y^j\ ,\label{xty}\eeq
transformes the quantum plane (\ref{qpcg}) into the rime quantum plane (\ref{qp})
exhibiting coordinate two-dimensional subplanes. The change of basis (\ref{xty}) can be written
in terms of a "generating function": let 
\beq G:=\sum_j e_j(\phi_1,\dots ,\phi_n )\, y^j\ .\label{gxty1}\eeq 
Then
\beq x^i=\frac{\partial G}{\partial\phi_i}\ .\label{gxty2}\eeq

\vskip .3cm\noindent
{\bf Remark.} The standard Drinfeld--Jimbo $R$-matrix admits, for a certain choice of
multi-parameters, a different rime form. The relations $u^iv^j=(\hat{R_c})^{ij}_{kl}v^ku^l$
for this choice are
\beq\label{stcl}\begin{array}{l}
u^iv^i=v^iu^i\ \ ,\\[.3em] u^iv^j=v^ju^i+(1-q^{-2})v^iu^j\ \ ,\ \ i<j\ \ ,\\[.3em]
u^iv^j=q^{-2}\, v^ju^i \ \ ,\ \ i>j\ \ . \end{array}\eeq
The left even space for this $R$-matrix is classical.

The change of variables with the matrix $\tilde{X}^i_j=1-\theta_{ji}$,
\beq U^i:=u^1+u^2+\dots +u^i\ ,\ V^i:=v^1+v^2+\dots +v^i\ ,\eeq
transforms the relations (\ref{stcl}) into
\beq\label{rstcl}\begin{array}{l}
U^iV^i=V^iU^i\ \ ,\\[.3em] U^iV^j=V^jU^i+(1-q^{-2})V^iU^j-(1-q^{-2})V^iU^i\ \ ,\ \ i<j\ \ ,\\[.3em]
U^iV^j=q^{-2}\, V^jU^i+ (1-q^{-2})V^jU^j\ \ ,\ \ i>j\ \ . \end{array}\eeq

The matrix $X$, defined by eq. (\ref{matX}), degenerates if $\phi_i=\phi_j$ for some $i$ and $j$. 
Interestingly, the $R$-matrix $X\otimes X\hat{R_c}\ X^{-1}\otimes X^{-1}$ admits
limits $\lim_{\phi_{\sigma (2)}\rightarrow 0}\lim_{\phi_{\sigma (3)}\rightarrow 0}\dots \lim_{\phi_{\sigma (n)}\rightarrow 0}$ for an arbitrary permutation $\sigma\in S_n$ and the result 
is always rime. In particular,
\beq \tilde{X}\otimes \tilde{X}\hat{R_c}\ \tilde{X}^{-1}\otimes \tilde{X}^{-1}=
\lim_{\phi_2\rightarrow 0}\lim_{\phi_3\rightarrow 0}\dots \lim_{\phi_n\rightarrow 0}
X\otimes X\hat{R_c}\ X^{-1}\otimes X^{-1}\ .\eeq
 
\section{Classical rime r-matrices\vspace{.25cm}}

The classical limit of an $R$-matrix is a {\it classical $r$-matrix}, a solution of the
{\it classical} Yang--Baxter (cYB) equation 
\beq [ r_{12},r_{13} ] + [ r_{12},r_{23} ] + [ r_{13},r_{23} ] =0\ .\eeq

\vskip .2cm
We are going to show that the classical limits of the rime $R$-matrices from Section 2 are
equivalent to the Cremmer-Gervais $r$-matrices in the non-skew-symmetric case and to the 
''boundary'' $r$-matrix of Gerstenhaber and Giaquinto \cite{GG} (see also \cite{BD2}; this 
$r$-matrix is attributed to A. G. Elashvili there) in the
skew-symmetric case. Similar equivalences appeared in the study of the gauge transformations of 
the dynamical $r$-matrices in the Calogero-Moser model \cite{FP1,FP2} \footnote{We thank 
L\'aszl\'o Feh\'er for drawing our attention to the references \cite{FP1,FP2}.}. 

\vskip .2cm
In the sequel we use the following conventions. An $R$-matrix acts in a space $V\otimes V$.
A basis of $V$ is $\{ e_i\}$ (labeled by a lower index); an operator $A$ in $V$ has matrix 
coefficients $A_i^j$, $A(e_i)=A_i^je_j$, so for a vector $\vec{v}=v^ie_i$ one has 
$(A\vec{v})^i=A^i_j\vec{v}^j$; the matrix units are $e^i_j$, $e^i_j(e_k)=\delta^i_ke_j$, so the 
multiplication rule is $e^i_je^k_l=\delta^i_l e^k_j$; $e_{\alpha_i}$  are the $\mathfrak{sl}(n)$ 
simple positive root elements, $e_{\alpha_i}=e^{i+1}_i$; $P$ is the permutation operator, 
$P(e_i\otimes e_j)=e_j\otimes e_i$, so $P(e^i_j\otimes e^k_l)= e^i_l\otimes e^k_j$ and 
$(PB)^{kl}_{ij}=B^{lk}_{ij}$ for an operator $B$ in $V\otimes V$ having matrix coefficients 
$B^{kl}_{ij}$, $B(e_i\otimes e_j)=B^{kl}_{ij}e_k\otimes e_l$.

\subsection{Non-skew-symmetric case}\label{nskc}

\begin{proposition}$\!\!\!${\bf .} The  non-unitary rime $R$-matrix (Proposition \ref{prop:rnu}) 
is a quantization of the non-skew-sym\-metric $r$-matrix
\beq\label{rcb} r=\sum_{i,j:i\neq j}\ \frac{\phi_i}{\phi_i -\phi_j}(e^{i}_{j}\otimes e^{j}_{i}
-e^{i}_{i}\otimes e^{j}_{j}+e^{i}_{i}\wedge e^{i}_{j})\ ,\eeq
where $x\wedge y:=x\otimes y-y\otimes x$. 
The change of basis with the matrix  $X^j_k = e_{k-1} \, (\phi_{1},\ldots, \hat{\phi}_{j}, 
\ldots,\phi_{n} )$ transforms $r$ into the parameter-free cYB solution $r_{CG}$
\beq \label{clcr}r_{CG}=\sum_{i,j:i< j}\ \sum_{s=1}^{j-i}(e_j^{i+s-1}\otimes 
e_i^{j-s+1}-e_i^{i+s-1}\otimes e_j^{j-s+1})\ .\eeq
\end{proposition}

\noindent {\bf Proof.}
The coefficients $\beta_{ij}$ (\ref{bphi}) are linear in the deformation parameter $\beta$
($\beta=0$ is the classical point).
Hence  
\beq R =1\!\!1\otimes 1\!\!1 + \beta r\ ,\eeq
where $R=P\hat{R}$ and $r$ is given by (\ref{rcb}).

The matrix $R_{CG}-1\!\! 1\otimes 1\!\!1$, where $R_{CG}= P \hat{R}_{CG}$, is linear with respect 
to the parameter $\beta =1-q^{-2}$ as well,
\beq R_{CG} = 1\!\! 1\otimes 1\!\!1 + \beta \, r_{CG}\eeq
thus the formula (\ref{change}) implies  $r  = (X \otimes X) \,r_{CG} \,\, (X^{-1} \otimes 
X^{-1})$.\hfill $\Box$

\vskip .2cm
We mentioned two ways of obtaining the numerical two-parametric $R$-matrix $(\hat{R}_{CG,p})$ from 
the $R$-matrix $(\hat{R}_{CG,1})$: by a diagonal twist and by the operation (\ref{chst}). There is 
one more way which consists of changing the representation. We shall illustrate it on the example 
of the classical $GL$ $r$-matrix (\ref{clcr}). A change of representation of the Lie algebra $GL$,
\beq e^i_j\mapsto e^i_j+c\,\delta^i_j 1\!\!1\ ,\label{chre}\eeq
where $c$ is a constant, produces the following effect on the $r$-matrix (\ref{clcr}):
\beq r_{CG}\mapsto r_{CG}+c\left(\eta\otimes 1\!\!1-1\!\!1\otimes\eta-(n-1)1\!\!1\otimes 1\!\!1
\phantom{\frac{}{}}\right)\ ,
\label{chrecg}\eeq
where $n=\,$dim$\,V$ and
\beq \eta =-\frac{n(n+1)}{2} 1\!\! 1
+\sum j \, e^j_j\ ,\ \ {\mathrm{tr}}\,\eta =0\ .\label{geinv}\eeq

The classical version of the operation (\ref{chst}) is as follows. Let $\eta$ be an arbitrary 
generator of the invariance group of an $r$-matrix $r$,
\beq [ r,\eta_1+\eta_2 ] =0\ .\label{chstc1}\eeq
Then the operator 
\beq r_{(c)}=r+c(\eta_1-\eta_2)\ ,\label{chstc2}\eeq 
where $c$ is a constant, is again a classical $r$-matrix (a solution of the cYBe). 

\vskip .2cm
The operator $\eta$ in (\ref{geinv}) is, up to a scale, the unique traceless generator
of the invariance group (see (\ref{invdcg})) of the $r$-matrix (\ref{clcr}). Thus, the 
representation change and the operation (\ref{chstc2}) give the same family of $r$-matrices (up to 
an addition of a multiple of the identity operator, which does not violate the cYBe). 

\subsection{BD triples.} 

Each block in the strict rime classical $r$-matrix (\ref{rcb}) looks even more "rimed",
\beq \label{crm} \left(\ba{rrrr}0 \, & 0\,& 0\,&0\, \\
\beta'_{12}&-\beta'_{12}& \beta'_{21}& - \beta'_{21} \\
-\beta'_{12}&\beta'_{12}& -\beta'_{21}&\beta'_{21} \\ 0\, &0 \, &0\, & 0\, \ea \right)\ ,\eeq
where $\beta'_{i j}= {\beta_{i j}}/{ \beta}={\phi_i}/{(\phi_i -\phi_j)}$.
The multiplication from the left by $P$ acts on each block as a permutation of the
second and third lines, so the rime $r$-matrix (\ref{crm}) enjoys  the symmetry $Pr= -r$. 
We shall now discuss this symmetry property in the context of Belavin--Drinfeld triples.
  
\vskip .2cm
In \cite{BD} Belavin and Drinfeld gave, for a simple Lie algebra 
$\mathfrak {g}$, a description of non-unitary (non-skew-symmetric) cYB solutions $r\in 
\mathfrak {g}\otimes \mathfrak {g}$, satisfying $r_{12}+r_{21}= t$, where
$t\in\mathfrak{g}\otimes\mathfrak{g}$ is the $\mathfrak{g}$-invariant element. 
The non-unitary solutions are put into correspondence with combinatorial objects called 
{\it Belavin--Drinfeld triples} (BD-triples for short). The Belavin--Drinfeld  triple  
$(\Pi_1,\Pi_2,\tau)$ for a simple Lie algebra $\mathfrak {g}$ consists of the following data: 
$\Pi_1,\Pi_2$ are subsets of the set of simple positive roots $\Pi$ of the algebra $\mathfrak {g}$ 
and $\tau$ is an invertible mapping: $\Pi_1 \rightarrow\Pi_2$ such that
$\langle \tau(\rho),\tau(\rho') \rangle = \langle \rho,\rho' \rangle$ for any 
$\rho,\rho' \in \Pi_1$ and $\tau^k(\rho) \neq \rho$ for any $\rho\in \Pi_1$ and any natural $k$
for which $\tau^k(\rho)$ is defined.

The $r$-matrix for a triple $(\Pi_1,\Pi_2,\tau)$ has the form
\beq\label{BRr} r=r_0+\sum_{\alpha\in\Delta_+}e_{-\alpha}\otimes e_{\alpha}
+\sum_{\alpha ,\beta\in \Delta_+ :\alpha <\beta}e_{-\alpha}\wedge e_{\beta}\ ,\eeq
where $\, <\,$ is a partial order on the set of positive roots $\Delta_+$ defined by the rule:
$\alpha <\beta$ for $\alpha,\beta \in \Delta_+$ if there exists a natural $k$ such that  
$\tau^k(\alpha )=\beta$. The part $r_0$ belongs to ${\mathfrak{h}}\otimes {\mathfrak{h}}$,
where ${\mathfrak{h}}$ is the Cartan subalgebra of $\mathfrak {g}$; $r_0$ contains 
continuous "multiparameters", which satisfy
\beq (\tau (\alpha )\otimes {\mathrm{id}}+{\mathrm{id}}\otimes \alpha )(r_0)=0
\ \ \ \ {\mathrm{for}}\ \ {\mathrm{all}}\ \ \ \alpha\in\Pi_1\ .\label{cocrt}\eeq

\vskip .2cm
We are dealing with matrix solutions $r$ of the cYB equation, $r\in\mathfrak{gl}(n)\otimes 
\mathfrak{gl}(n)$, so $r_{12} + r_{21}$ can be a linear combination of $P$ and 
$1\!\! 1\otimes 1\!\!1$.

Let $\Pi=\{ \alpha_1, \ldots, \alpha_{n-1} \}$ be the set of the positive simple roots
for the Lie algebra $\mathfrak {sl}(n)$. 

There are two Cremmer--Gervais BD triples, $\mathfrak{T}_+$ and 
$\mathfrak{T}_-$. For the Cremmer--Gervais triple BD-triple $\mathfrak{T}_+$,
$\Pi_1=\{\alpha_1,\alpha_2,\ldots ,\alpha_{n-2}\}$, $\Pi_2=\{\alpha_2,\alpha_3,\ldots ,
\alpha_{n-1} \}$ and $\tau(\alpha_i)=\alpha_{i+1}$.
 The data $(\Pi_1,\Pi_2,\tau)$ is encoded in the graph 
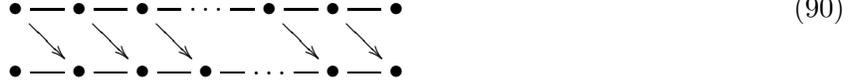
\beq\label{lrdiag}\xymatrix@!0{ 
\bullet \ar@{-}[r]\ar@{->}[dr] &\bullet \ar@{-}[r] \ar@{->}[dr] &\bullet 
\ar@{-}[r] \ar@{->}[dr]&  \ldots \ar@{-}[r] &\bullet \ar@{-}[r] \ar@{->}[dr] 
 &\bullet  \ar@{-}[r]\ar@{->}[dr] &\bullet\\
\bullet  \ar@{-}[r] &\bullet \ar@{-}[r]  &\bullet\ar@{-}[r] &  \bullet \ar@{-}[r] & \ldots   \ar@{-}[r] &\bullet  \ar@{-}[r] &\bullet } \eeq

\vskip .1cm
The triple $\mathfrak{T}_-$ can be obtained from the triple $\mathfrak{T}_+$ either
by setting $\Pi_1'=\Pi_2$, $\Pi_2'=\Pi_1$ and $\tau'=\tau^{-1}$ or by applying the outer
automorphism of the underlying  $A_{n-1}$ Dynkin diagram; the graph corresponding to
the triple $\mathfrak{T}_-$ is

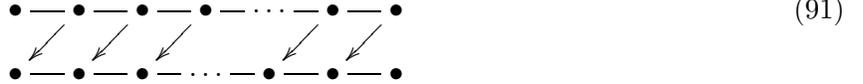
\beq\label{rldiag}\xymatrix@!0{ 
\bullet\ar@{-}[r] &\bullet \ar@{-}[r] \ar@{->}[dl] &\bullet\ar@{-}[r] \ar@{->}[dl] 
&\bullet\ar@{-}[r] \ar@{->}[dl]& \ldots   \ar@{-}[r] &\bullet  \ar@{-}[r]\ar@{->}[dl]
&\bullet\ar@{->}[dl]\\
\bullet\ar@{-}[r] &\bullet \ar@{-}[r] &\bullet \ar@{-}[r]&  \ldots \ar@{-}[r] 
&\bullet\ar@{-}[r]  &\bullet  \ar@{-}[r] &\bullet}\eeq

\vskip .2cm\noindent
The $r$-matrix (\ref{clcr}) corresponds to the triple (\ref{lrdiag}) for a certain choice of the 
multiparameters. Here is the $r$-matrix $r'$ corresponding to the triple 
(\ref{rldiag})
\beq \label{clcrr} r_{CG}'=\sum_{i,j:i< j}\ \sum_{s=1}^{j-i}(e^i_{j-s+1}\otimes e^j_{i+s-1}- 
e^j_{j-s+1}\otimes e^i_{i+s-1})\ \eeq
for a certain choice of the multiparameters, for which it satisfies $r'P=-r'$.

For the $r$-matrices (\ref{clcr}) and (\ref{clcrr}), one has $r_{12}+r_{21}=
P-1\!\!1\otimes 1\!\!1$. The Cartan part of the $r$-matrices (\ref{clcr}) and (\ref{clcrr}) are
\beq r_0=-\sum_{i,j:i<j}\ e^i_i\otimes e^j_j\ \ \ ,\ \ \ r_0'=-\sum_{i,j:i<j}\ e^j_j\otimes e^i_i
\ .\label{capar}\eeq

The following lemma shows that a classical $r$-matrix $r$ for a triple $\mathfrak{T}$ can have a 
symmetry with respect to the multiplication by $P$ from one side if and only if all segments 
(connected components) of $\Pi_1$ are mapped by $\tau$ according to either (\ref{lrdiag}) or (\ref{rldiag}). 

\begin{lemma}$\!\!\!${\bf .} A non-skew-symmetric classical $r$-matrix with a Belavin--Drinfeld 
data $(\Pi_1,\Pi_2,\tau)$ can satisfy $Pr=-r$ (respectively, $rP=-r$) for a certain choice of 
the multiparameters if and only if $\tau(\alpha_i)=\alpha_{i+1}$ (respectively, 
$\tau(\alpha_i)=\alpha_{i-1}$) for all $i\in \Pi_1$.
\end{lemma}

\noindent {\bf Proof.} Assume that $\tau(\alpha_{m})=\alpha_{m+k}$ for some natural $k$, 
$k \geq 1$. Then $r$ contains the term $e^{m+k}_{m+k+1}\wedge e^{m+1}_{m}$ with the coefficient
1. Such $r$-matrix cannot satisfy $rP=-r$ for if $rP=-r$ then $r$ contains the term 
$e^{m+1}_{m+k+1}\wedge e^{m+k}_{m}$ with the coefficient $(-1)$ but the coefficient in 
$e_{-\alpha}\wedge e_\beta$ is 1 in the formula (\ref{BRr}).

If $Pr=-r$ then $r$ should contain also the term $e^{m+1}_{m+k+1}\wedge e^{m+k}_{m}$. It then 
follows that

\vskip .2cm
\noindent (i) the Lie subalgebra generated by $\Pi_1$ contains $e^{m+k}_{m}$ therefore
the interval $[ \alpha_m,\alpha_{m+1},\dots ,\alpha_{m+k-1} ]$ is contained in $\Pi_1$;

\vskip .1cm
\noindent (ii) the Lie subalgebra generated by $\Pi_2$ contains $e^{m+1}_{m+k+1}$ therefore
the interval $[ \alpha_{m+1},\alpha_{m+2},\dots ,\alpha_{m+k} ]$ is contained in $\Pi_2$;

\vskip .1cm
\noindent (iii) the image of the interval $[ \alpha_m,\alpha_{m+1},\dots ,\alpha_{m+k-1} ]$
under $\tau$ is the interval $[ \alpha_{m+1},\alpha_{m+2},\dots ,\alpha_{m+k} ]$.

\vskip .2cm
This implies that the interval $[ \alpha_{m+1},\alpha_{m+2},\dots ,\alpha_{m+k-1} ]$
is $\tau$-invariant (since $\tau(\alpha_{m})=\alpha_{m+k}$) which contradicts to the nilpotency of 
$\tau$ unless this interval is empty, that is, $k=1$.

\vskip .2cm
Similarly, $rP=-r$ is possible only if $\tau(\alpha_i)=\alpha_{i-1}$ for all $i\in \Pi_1$.

\vskip .2cm
It  is left to show that when  $\tau(\alpha_i)=\alpha_{i+1}$ (respectively, 
$\tau(\alpha_i)=\alpha_{i-1}$) for all $i\in \Pi_1$ the multiparameters can indeed be adjusted 
to fulfill $Pr=-r$ (respectively, $rP=-r$). We leave it as an exercise for the reader
to check that with the assignment (\ref{capar}) for $r$ (respectively, for $r'$) the 
compatibility condition (\ref{cocrt}) is verified. The proof is finished. \hfill $\Box$

\vskip .3cm\noindent {\bf Remark.} Two extreme BD triples can be rimed, the empty 
(Drinfeld--Jimbo) one and the ``maximal'' Cremmer--Gervais one. However, not every triple
can be rimed: already the triple
\beq\label{ftr}\xymatrix@!0{\bullet \ar@{-}[r]\ar@{->}[drr] &\bullet \ar@{-}[r]& \bullet\\
\bullet \ar@{-}[r] & \bullet \ar@{-}[r] &\bullet}\eeq
provides a counterexample. We outline a computer-aided proof in appendix C. 

\subsection{Skew-symmetric case}\label{sksc}

A skew-symmetric classical $r$-matrix $r\in\mathfrak {g}\wedge\mathfrak{g}$ is canonically
associated with a quasi-Frobenius Lie subalgebra $(\mathfrak{f},\omega)$ of $\mathfrak{g}$ 
(see, {\em e.g.}, \cite{St}). A Lie algebra  $\mathfrak{f}$  which admits a non-degenerate 
2-cocycle $\omega$ is called {\it quasi-Frobenius}; it is {\it Frobenius} if $\omega$ is a coboundary, i.e., 
$\omega(X,Y)= \lambda([X,Y])$ for some $\lambda\in\mathfrak{f}^{\ast}$.

We describe now the skew-symmetric $r$-matrix arising in the classical limit of the unitary rime 
$R$-matrix from Proposition \ref{prop: ru}. 

\begin{proposition}$\!\!\!${\bf .} The  unitary rime $R$-matrix (Proposition \ref{prop: ru}) is a 
quantization of the skew-symmetric $r$-matrix
\beq\label{rcc} r=\sum_{i,j:i<j}\ \frac{1}{\mu_i-\mu_j}(e^i_j-e^j_j)\wedge (e^j_i-e^i_i) \ 
\in\, \mathfrak {gl}(n) \wedge \mathfrak {gl}(n)\ .\eeq
This skew-symmetric classical $r$-matrix corresponds to a Frobenius Lie algebra $(\mathfrak 
{g}_0(n),\delta\lambda_n)$ spanned by the generators $Z^i_j:=e^i_j- e^j_j$, $i\neq j$, 
with the Frobenius structure determined by the coboundary of the $1$-cochain $\lambda_n=-{\displaystyle{\sum_{i,j:i\neq j}}}\,\mu_i z^i_j$, where $\{ z^i_j\}$, $i\neq j$, is the basis in $\mathfrak {g}^{\ast}_0(n)$, dual to the basis $\{ Z^i_j\}$ in $\mathfrak{g}_0(n)$, 
$z^i_j(Z^k_l)=\delta^i_l\delta_j^k$.\end{proposition}
  
\noindent {\bf Proof.} An artificial introduction of a small parameter $c$ by a rescaling
$\mu_i\mapsto c^{-1}\mu_i$ in the formula for the $R$-matrix $\hat{R}$ in Proposition 
\ref{prop: ru} gives
\beq R =1\!\!1\otimes 1\!\!1 + c\, r\ ,\eeq
where $r$ is given by (\ref{rcc}).
 
\vskip .2cm
The $n(n-1)$ matrices $Z^i_j:=e^{i}_{j}- e^{j}_{j}$, $i\neq j$, form an associative subalgebra
of the matrix algebra,
\beq Z^j_iZ^k_l=(\delta^j_l-\delta^i_l)(Z^k_i-Z^l_i)\ \eeq 
(we set $Z^i_i=0$ for all $i$);
with respect to the commutators these matrices form a Lie subalgebra $\mathfrak {g}_0(n)$
of the Lie algebra $\mathfrak{gl}(n)$, $\mathfrak{g}_0(n)\subset\mathfrak{gl}(n)$:
\beq\label{zz} [ Z^i_j,Z^j_i ] =Z^j_i-Z^i_j\ ,\qquad [ Z_i^j, Z_i^k ] =Z_i^j-Z_i^k\ ,\qquad  
 [ Z^i_j,Z^j_k ] =Z^j_k-Z^i_k\ ,\quad i\neq j\neq k\neq i\ ,\eeq
all other brackets vanish. The skew-symmetric solution (\ref{rcc}) of the cYB equation,  
\beq\label{rcc2} r=\sum_{i,j:i<j}\,\frac{Z^i_j\wedge Z^j_i}{\mu_i-\mu_j}\ ,\eeq
is non-degenerate on the {\it carrier} subalgebra $\mathfrak {g}_0(n)$.
The  carrier subalgebra $\mathfrak {g}_0(n)$ is necessarily quasi-Frobenius, having a 2-cocycle 
$\omega$ given by the inverse of the $r$-matrix, that is,
\beq\omega(Z_A,Z_B)=r_{AB}\ ,\quad\mbox{where}\quad r^{AB}r_{BC}=\delta^A_C\ \ ,\ \  
r=\sum_{A,B}r^{AB}Z_A\wedge Z_B\ .\eeq
We have
\beq\omega(Z^i_j,Z^k_l)=-(\mu_i-\mu_j)\delta^l_i\delta^j_k\ .\eeq
It is easy to check that the  2-cycle $\omega$ is a coboundary,
\beq\omega(Z^i_j,Z^k_l)=\lambda_n([Z^i_j,Z^k_l])\ \ ,\quad\lambda_n=-\sum_{i,j:i\neq j}\,\mu_i
z^i_j \in \mathfrak {g}_0^{\ast}(n)\ ,\eeq
thus the subalgebra $\mathfrak {g}_0(n)$ is Frobenius. \hfill $\Box$

\vskip .2cm
The ''Frobenius'' $r$-matrix (\ref{rcc}) (and its quantization) was considered in the 
work \cite{AF}. 

\begin{proposition}$\!\!\!${\bf .}\label{btor}
The skew-symmetric rime classical $r$-matrix (\ref{rcc}), $r=\sum_{i<j}(\mu_i-\mu_j)^{-1} 
Z^i_j\wedge Z^j_i$, where $\mu=(\mu_1,\mu_2,\ldots ,\mu_n)$ is an arbitrary vector such that 
$\mu_i \neq \mu_j$, belongs to the orbit of the parameter-free classical $r$-matrix 
\beq b=\sum_{i,j:i<j}\ \sum_{k=1}^{j-i}e_i^{i+k}\wedge e_j^{j-k+1}\ .\label{bee}\eeq
More precisely,
\beq r=Ad_{X_{\mu}}\otimes Ad_{X_{\mu}}(b)\ ,\label{bee'}\eeq
where the element $X_{\mu} \in GL(n)$ is defined by $(X_{\mu})^j_k = e_{k-1}\, 
(\mu_{1},\ldots, \hat{\mu}_{j}, \ldots,\mu_{n} )$.\footnote{This matrix 
is the same $X$ as in Proposition \ref{transem} but depending on variables $\mu_i$.}
\end{proposition}

\noindent {\bf Proof.} The equality $r=Ad_{X_{\mu}}\otimes Ad_{X_{\mu}}(b)$ is equivalent to
a set of relations for the elementary symmetric functions $e_i$,
\beq\label{eqb} (X_{\mu}\otimes X_{\mu} )\, b=r\, (X_{\mu}\otimes X_{\mu})
\qquad\Leftrightarrow\qquad {\sum}_{r,s}e_{r-1}^{\hat{i}}e_{s-1}^{\hat{j}}\,\, b^{rs}_{kl}=
{\sum}_{a,b} {r}^{ij}_{ab} e_{k-1}^{\hat{a}}e_{l-1}^{\hat{b}}\ ,\eeq
where 
$$b^{ij}_{ab}= \sum_{k=1}^{j-i}\delta_{b}^{j-k+1}\delta^{i+k}_{a} 
-\sum_{k=1}^{i-j} \delta^{i-k+1}_{a}\delta^{j+k}_{b}\ \ {\mathrm{and}}\ \  
r^{ij}_{ab}=\left\{\begin{array}{ccc}\!
(\delta^{i}_{a}\delta^{i}_{b}+\delta^{j}_{a}\delta^{j}_{b}
-\delta^{i}_{a}\delta^{j}_{b}-\delta^{j}_{a}\delta^{i}_{b})/(\mu_i -\mu_j)&,&i\neq j\ ,\\[.5em]
\! 0&,&i=j\ .\end{array}\right. $$
Both operators $b^{ij}_{ab}$ and  $r^{ij}_{ab}$ are symmetric in the lower indices and 
anti-symmetric in the upper indices, that is, 
\beq Pb=-b\ ,\ bP=b  \ \ \ {\mathrm{and}}\ \ \ Pr=-r\ ,\ rP=r\ .\eeq
Eqs. (\ref{eqb}) have the following form
\beq -\sum_{s\geq 1}(e^{\hat{i}}_{b+s-2}e^{\hat{j}}_{a-s-1}-e^{\hat{j}}_{b+s-2} 
e^{\hat{i}}_{a-s-1})=\frac{1}{\mu_i-\mu_j}(e^{\hat{i}}_{a-1}-e^{\hat{j}}_{a-1})(e^{\hat{i}}_{b-1}- 
e^{\hat{j}}_{b-1})\ .\label{keke}\eeq
Due to (\ref{sectype}), the left hand side of (\ref{keke}) equals
\beq -\frac{1}{\mu_i-\mu_j}(\mu_i e^{\hat{i}}_{a-2}-\mu_j e^{\hat{j}}_{a-2})
(e^{\hat{i}}_{b-1}-e^{\hat{j}}_{b-1})\ .\label{keke2}\eeq
The right hand side of (\ref{keke}) equals (\ref{keke2}) as well because
$e^{\hat{i}}_{a-1}=e_{a-1}-\mu_i e^{\hat{i}}_{a-2}$.\hfill $\Box$
 
\vskip .2cm 
As in the non-skew-symmetric case, in the infinite dimension the operators $b$ and $r$ 
are in general not equivalent.

\vskip .3cm
{\bf The $\mathfrak{sl}(n)$ cYB solution.} Let $I=\sum_{i=1}^ne^i_i$ be the central element of 
$\mathfrak{gl}(n)$. The generators  $\tilde{Z}^{i}_{j}=Z^{i}_{j}+\frac{1}{n}I\in\mathfrak{sl}(n)$ satisfy the same relations (\ref{zz}) as ${Z}^{i}_{j}$ thus 
they form a subalgebra $\tilde{\mathfrak{g}}_0(n)$ of the Lie algebra $\mathfrak{sl}(n)$
which is isomorphic to $\mathfrak {g}_0(n)$, $\tilde{\mathfrak{g}}_0(n)\simeq\mathfrak{g}_0(n)$. This isomorphism gives rise to another solution $\tilde{r}\in \mathfrak{sl}(n)\wedge 
\mathfrak{sl}(n)$ of the cYB equation,
\beq\label{rcc2'}\tilde{r}=\sum_{i,j:i<j}\,\frac{\tilde{Z}^i_j\wedge\tilde{Z}^j_i}{\mu_i-\mu_j}   
\ \in\,\mathfrak{sl}(n)\wedge\mathfrak{sl}(n)\ .\eeq
We have the following lemma about the carrier Lie algebra of $\tilde{r}$ (the Lie subalgebra of 
$\mathfrak{sl}(n)$ spanned by the generators $\tilde{Z}^{i}_{j}$).
 
\begin{lemma}$\!\!\!${\bf .} The  subalgebra $\tilde{\mathfrak{g}}_0(n)\subset\mathfrak{sl}(n)$ of 
dimension $\dim\tilde{\mathfrak{g}}_0(n)=n(n-1)$ is isomorphic to the maximal parabolic subalgebra 
$\mathfrak{p}$ of $\mathfrak{sl}(n)$ obtained by deleting the first negative root.\end{lemma}
 
\noindent {\bf Proof.} The vector $v=\sum_{i=1}^n e_i$ is an eigenvector for all elements
$\tilde{Z}^i_j$,
\beq\tilde{Z}^i_j(v)=\frac{1}{n}v\ \ \ {\mathrm{for}}\ {\mathrm{all}}\ i\ {\mathrm{and}}\ j\ ,
\ i\neq j\ .\eeq
In a basis in which the first vector is $v$, the linear span of the generators $\tilde{Z}^i_j$ is
\beq\label{parabol}\left(\ba{cccc}\ast&\ast&\dots&\ast\\0&\ast&&\ast\\
\vdots&\vdots&&\vdots\\0 &\ast&\dots&\ast\ea\right)\ ,\eeq
with the traceless condition. The comparison of dimensions finishes the proof.\hfill $\Box$ 

\vskip .2cm
Gerstenhaber and Giaquinto \cite{GG} found a classical $r$-matrix  $b_{CG}$ which they called ``boundary'' because it lies in the closure of the solution space of the YB equation.
The cYB solution $b_{CG}$ corresponds to a Frobenius subalgebra $(\mathfrak{p},\Omega)$, where $\mathfrak{p}$ is the parabolic subalgebra of $\mathfrak{sl}(n)$ as above and the 2-cocycle $\Omega$ is a coboundary,
\beq \Omega=\delta \lambda_{b_{CG}} \qquad ,\qquad \lambda_{b_{CG}}=\sum_{i=1}^n (e_{i+1}^i)^{\ast}
\in \mathfrak{p}^{\ast}\ .\eeq
The $r$-matrix $b_{CG}$ is a twist of $b$ (see \cite{EH}). 

\vskip .2cm
Since the carriers of $\tilde{r}$ and $b_{CG}$  are isomorphic, the $r$-matrices are equivalent.
We shall now prove that the same matrix $X_{\mu}$ transforms $b_{CG}$ into $\tilde{r}$.

\begin{proposition}$\!\!\!${\bf .} The boundary classical $r$-matrix
$b_{CG}\in\mathfrak{sl}(n)\wedge\mathfrak{sl}(n)$,  
\beqa\label{bcg}b_{CG}&=&\sum_{i,j}(1-\frac{j}{n})\, e^i_i\wedge e_j^{j+1}+
\sum_{i,j:i<j}\ \sum_{k=1}^{j-i}e_i^{i+k}\wedge e_j^{j-k+1}\ ,\eeqa
transforms into the cYB solution $\tilde{r}\in\mathfrak{sl}(n)\wedge\mathfrak{sl}(n)$, 
\beq\label{rcc2''} \tilde{r}=\sum_{i,j:i<j}\ \frac{\tilde{Z}^i_j\wedge\tilde{Z}^j_i}{\mu_i-\mu_j}  
\ ,\ \ \mbox{where}\quad\tilde{Z}^{i}_{j}=e^i_j-e^j_j+\frac{1}{n}\sum_{i=1}^n e^i_i\ ,\eeq
by a change of basis with the matrix $X_{\mu}\in GL(n)$, 
\beq\tilde{r}=Ad_{X_{\mu}}\otimes Ad_{X_{\mu}}(b_{CG})\ .\eeq
\end{proposition}

\noindent {\bf Proof.} Due to Proposition \ref{btor} we have $r=Ad_{X_{\mu}}\otimes Ad_{X_{\mu}} 
(b)$. The cYB solution $b_{CG}$ is the sum of $b$ and other terms,
$b_{CG}=b+\sum_{i,j}(1-\frac{j}{n})\, e^i_i\wedge e^{j+1}_{j}$.
Therefore it is enough to show that $\tilde{r}-r = Ad_{X_{\mu}}\otimes Ad_{X_{\mu}}(b_{CG}-b)$.
One has
\beq\tilde{r}-r=\frac{1}{n}\ I\wedge\sum_{i,j:i\neq j}\,\frac{Z^j_i}{\mu_i-\mu_j}\ \ ,\ \ 
b_{CG}-b=I\wedge\sum_{j}(1-\frac{j}{n})\, e_j^{j+1}\ .\eeq
Thus we have to show that 
\beq X_{\mu}\ \sum_{j}(1-\frac{j}{n})\, e_j^{j+1}=\frac{1}{n}
\sum_{i,j:i\neq j}\,\frac{Z^j_i}{\mu_i-\mu_j}\ X_{\mu}\ ,\eeq
which amounts to the following identities for the elementary symmetric functions:
\beq (1-\frac{b-1}{n})e^{\hat{i}}_{b-2}=\frac{1}{n}\sum_{j:j\neq i} 
\frac{e^{\hat{j}}_{b-1}-e^{\hat{i}}_{b-1}}{\mu_i-\mu_j}\ .\eeq
Replacing, in the right hand side, $e^{\hat{j}}_{b-1}$ by $e^{\hat{i}\hat{j}}_{b-1}+\mu_i e^{\hat{i}\hat{j}}_{b-2}$, $e^{\hat{i}}_{b-1}$ by $e^{\hat{i}\hat{j}}_{b-1}+\mu_j e^{\hat{i}\hat{j}}_{b-2}$ and noticing that $\sum_i e^{\hat{i}}_c=(n-c)e_c$, $c=1,2,\dots ,n$ 
(for the elementary symmetric functions in $n$ variables) finishes the proof.\hfill $\Box$

\vskip .2cm
The passage to the $\mathfrak{sl}(n)$ solution is another instance of the representation change. The general representation change (\ref{chre}) produces the following effect on the 
numerical $r$-matrix (\ref{bee}):
\beq b\mapsto b-c\eta^{(0)}\wedge 1\!\!1\ ,\label{unopa}\eeq
where $\eta^{(0)}$ is the generator of the invariance group of the $r$-matrix (\ref{bee}),
\beq \eta^{(0)}=\sum (n-j)e^{j+1}_j\ .\label{guninv}\eeq
The representation change and the operation (\ref{chstc2}) produce the same 1-parametric family 
(\ref{unopa}) of skew-symmetric $r$-matrices. The choice $c=-1/n$ corresponds to the $r$-matrix
$b_{CG}$.

\section {B\'ezout operators\vspace{.25cm}}

The B\'ezout operator \cite{B} is the following endomorphism $\mathfrak{b}^{(0)}$ of the space 
$\mathfrak{P}$ of polynomials of two variables $x$ and $y$:
\beq \mathfrak{b}^{(0)}f(x,y)=\frac{f(x,y)-f(y,x)}{x-y}\ \ {\mathrm{or}}\ \ 
\mathfrak{b}^{(0)}=\frac{1}{x-y}(I-P)\ ,\label{bez1}\eeq
where $I$ is the identity operator and $P$ is a permutation, $Pf(x,y)=f(y,x)$. For any natural
$n$, the subspace $\mathfrak{P}_n$ of polynomials of degree less than $n$ in $x$ and less than 
$n$ in $y$ is invariant with respect to the operator $\mathfrak{b}^{(0)}$. The matrix of the
restriction of $\mathfrak{b}^{(0)}$ onto $\mathfrak{P}_n$, written in the basis $\{ x^ay^b\}$ of 
powers (in the decreasing order) coincides with the operator (\ref{bee}). 

\vskip .2cm
The non-skew-symmetric matrix (\ref{clcr}) is the matrix of the operator
\beq \mathfrak{b}=\frac{x}{x-y}(I-P)\ \label{bez2}\eeq
in this basis. The rime bases are formed by the non-normalized Lagrange polynomials 
$\{ l_i(x)l_j(y)\}$, $l_i(t)=\displaystyle{\prod_{s:s\neq i}}(t-\phi_s)$, at points $\{ \phi_i\}$, 
$i=1,2,\dots ,n$.

\vskip .2cm
We shall call the operators $\mathfrak{b}^{(0)}$ and $\mathfrak{b}$ {\it B\'ezout $r$-matrices}.
The B\'ezout $r$-matrices were rediscovered in several different contexts related to the 
Yang--Baxter equation (except the fact that they are the Cremmer--Gervais $r$-matrices, they 
appear, for instance, in \cite{D} and \cite{L}).

\vskip .2cm
The standard $r$-matrix $r^{(s)}$, for the choice of the multi-parameters for which it can be
non-trivially rimed (see the remark at the end of section \ref{sectionCG}), has the following
form in terms of polynomials
\beq r^{(s)}\, :\, x^iy^j\mapsto \theta (i-j)x^iy^j-\theta (j-i)x^jy^i\ .\label{bez3}\eeq
The subspaces $\mathfrak{P}_n$ are invariant with respect to $r^{(s)}$.

\vskip .2cm
The properties of the B\'ezout $r$-matrices $\mathfrak{b}^{(0)}$ and $\mathfrak{b}$ (and of the 
operator $r^{(s)}$) become more transparent when they are viewed as operators on polynomials. In 
particular, 
\beqa &(\mathfrak{b}^{(0)})^2=0\ ,\ \mathfrak{b}^{(0)}P=-\mathfrak{b}^{(0)}\ ,\ 
P\mathfrak{b}^{(0)}=\mathfrak{b}^{(0)}\ ,\ 
\mathfrak{b}^{(0)}+\mathfrak{b}^{(0)}_{21}=0\ ,&\label{bez4}\\[.3em]
&\mathfrak{b}^2=\mathfrak{b}\ ,\ \mathfrak{b}P=-\mathfrak{b}\ ,\ 
\mathfrak{b}+\mathfrak{b}_{21}=I-P\ ,&\label{bez5}\\[.3em]
&(r^{(s)})^2=r^{(s)}\ ,\ r^{(s)}P=-r^{(s)}\ ,\ r^{(s)}+r^{(s)}_{21}=I-P\ .&\label{bez6}\eeqa

The description of the invariance groups of the operators $\mathfrak{b}^{(0)}$ and 
$\mathfrak{b}$ is especially transparent when these operators are viewed as operators on the
spaces of polynomials. Let $\partial_x$ and $\partial_y$ be the derivatives in $x$ and $y$.
We have $(\partial_x+\partial_y)\Bigl(\displaystyle{\frac{1}{x-y}}\Bigr) =0$ which implies that 
$\partial_x$ is the generator of the invariance group of $\mathfrak{b}^{(0)}$; the group is formed 
by translations. Similarly, $(x\partial_x+y\partial_y)\Bigr(\displaystyle{\frac{x}{x-y}}\Bigr) =0$ 
which implies that $x\partial_x$ is the generator of the invariance group of $\mathfrak{b}$; the 
group is formed by dilatations. The operation (\ref{chstc2}) implies that the operators 
\beq \mathfrak{b}^{(0)}+c(\partial_x-\partial_y)\ ,\ \mathfrak{b}+c(x\partial_x-y\partial_y)
\label{bez31}\eeq
are solutions of the cYBe (the quantum version is easy as well) for an arbitrary constant $c$.

\subsection{Non-homogeneous associative classical Yang--Baxter equation}

The operators $\mathfrak{b}^{(0)}$, $\mathfrak{b}$ and $r^{(s)}$ satisfy an equation stronger
than the cYBe. For an endomorphism $r$ of $V\otimes V$, define
\beq r\circ r:=r_{12}r_{13}+r_{13}r_{23}-r_{23}r_{12}\ ,\ 
r\circ' r:=r_{13}r_{12}+r_{23}r_{13}-r_{12}r_{23}\ .\label{bez7}\eeq
The equation $r\circ r=0$ (as well as $r\circ' r=0$) is called {\it associative 
classical Yang--Baxter equation} (acYBe)
\cite{A,M}. 

\vskip .2cm
We introduce a {\it non-homogeneous associative classical Yang--Baxter equation} (nhacYBe):
\beq r\circ r =c r_{13}\ ,\label{bez8}\eeq
where $c$ is a constant.

\vskip .2cm
Let ${\cal{F}}$ be the space of polynomials in one variable. For the space 
${\cal{F}}\otimes {\cal{F}}$ of polynomials in two variables, we denote by $x$ (respectively, $y$) 
the generator of the first (respectively, second) copy of ${\cal{F}}$. For 
${\cal{F}}\otimes {\cal{F}}\otimes {\cal{F}}$, the generators are denoted by $x$, $y$ and $z$.

\begin{lemma}$\!\!\!${\bf .} \label{lemfornhacYBe}
1. Let $M$ be an operator on the space ${\cal{F}}\otimes {\cal{F}}$. Assume that 
\beqa &M(xf)=f+yM(f)\ ,&\label{b0b}\\
&M(yf)=-f+xM(f)&\label{b0b2}\eeqa
for an arbitrary $f\in {\cal{F}}\otimes {\cal{F}}$. Then\footnote{Eq. 
$M\!\circ\! M(xF)=z\, M\!\circ\! M(F)$ follows from (\ref{b0b}) alone.}
\beq M\!\circ\! M(xF)=z\, M\!\circ\! M(F)\ ,\ M\!\circ\! M(yF)=x\, M\!\circ\! M(F)\ ,\ 
M\!\circ\! M(zF)=y\, M\!\circ\! M(F)\ \label{b0b3}\eeq
for an arbitrary $F\in {\cal{F}}\otimes {\cal{F}}\otimes {\cal{F}}$.

\vskip .1cm\noindent  2. The operator $M=\mathfrak{b}^{(0)}$ verifies (\ref{b0b}) and (\ref{b0b2}).

\vskip .1cm\noindent  3. 
Moreover, the unique solution of eqs. (\ref{b0b}) and (\ref{b0b2}) (for the operator $M$
on the space ${\cal{F}}\otimes {\cal{F}}$)
together with the "initial" condition $M(1)=0$ is $M=\mathfrak{b}^{(0)}$.
\end{lemma}

\noindent {\bf Proof.} A direct calculation.\hfill $\Box$

\vskip .2cm
\begin{proposition}$\!\!\!${\bf .} \label{lemnhacYBe}
1. The B\'ezout operator $\mathfrak{b}^{(0)}$ satisfies the acYBe.

\vskip .1cm\noindent 2. The B\'ezout operator $\mathfrak{b}$ and the operator $r^{(s)}$ satisfy 
the nhacYBe with $c=1$.\end{proposition}

\noindent {\bf Proof.} A direct calculation for $\mathfrak{b}^{(0)}$. Another way is to  
notice that the relations (\ref{b0b3}) for $M=\mathfrak{b}^{(0)}$ reduce the verification of 
$\mathfrak{b}^{(0)}\circ \mathfrak{b}^{(0)} (F)=0$ for a monomial $F\in 
{\cal{F}}\otimes {\cal{F}}\otimes {\cal{F}}$ to the case $F=1$, which is trivial.

\vskip .2cm
For the B\'ezout operator $\mathfrak{b}\equiv x\mathfrak{b}^{(0)}$ ($x$ here is the operator of multiplication by $x$), we have, for an arbitrary $F\in {\cal{F}}\otimes {\cal{F}}\otimes {\cal{F}}$,
\beq \begin{array}{lcl}\mathfrak{b}\circ\mathfrak{b}\, (F)
&=&x\mathfrak{b}_{12}^{(0)}(x\mathfrak{b}_{13}^{(0)}(F))+
x\mathfrak{b}_{13}^{(0)}(y\mathfrak{b}_{23}^{(0)}(F))- 
y\mathfrak{b}_{23}^{(0)}(x\mathfrak{b}_{12}^{(0)}(F))\\[1em]
&=&x\left(\mathfrak{b}_{13}^{(0)}(F)+y\mathfrak{b}_{12}^{(0)}\mathfrak{b}_{13}^{(0)}(F)\right)
+xy\mathfrak{b}_{13}^{(0)}\mathfrak{b}_{23}^{(0)}(F)- 
xy\mathfrak{b}_{23}^{(0)}\mathfrak{b}_{12}^{(0)}(F)\\[1em]
&=&x\mathfrak{b}_{13}^{(0)}(F)+xy\mathfrak{b}^{(0)}\circ\mathfrak{b}^{(0)}(F)=
\mathfrak{b}_{13}(F)\ .\end{array}\label{b0b4}\eeq
We used eq. (\ref{b0b}) for $\mathfrak{b}^{(0)}$ in the second equality.

\vskip .2cm
For the operator $r^{(s)}$, the identity
\beq \theta (i-k)\theta (i-j)+\theta (i-k)\theta (j-k)-\theta (i-j)\theta (j-k)=\theta (i-k)\eeq
for the step function $\theta$ is helpful.\hfill $\Box$

\vskip .2cm
In each of cases (\ref{bez4}-\ref{bez6}), the operator $r$ satisfies a quadratic equation
$r^2=u_1r+u_2I$, the relation $r+r_{21}=\alpha P+\beta I$ with some constants $\alpha$
and $\beta$ and the nhacYBe with some constant $c$. Several general comments about relations 
between the constants appearing in these equations are in order here.

\paragraph{1.} Assume that an $r$-matrix (a solution of the cYBe) satisfies $r\circ r =c 
r_{13}$. Then $r\circ' r =c r_{13}$. Taking the combinations $(r\circ r -c r_{13})-P_{23}
(r\circ' r -c r_{13})P_{23}$ and $(r\circ r -c r_{13})-P_{12}
(r\circ' r -c r_{13})P_{12}$, we find
\beq r_{13}(Sr)_{23}-(Sr)_{23}r_{12}=c(r_{13}-r_{12})\ ,\ (Sr)_{12}r_{13}-r_{23}(Sr)_{12}=
c(r_{13}-r_{23})\ ,\label{bez9}\eeq
where $(Sr)_{12}:=r_{12}+r_{21}$. If $(Sr)_{12}=\alpha P_{12}+\beta I$ with some constants
$\alpha$ and $\beta$, as in (\ref{bez4}-\ref{bez6}), then it follows from (\ref{bez9}) that 
$(\beta -c)(r_{13}-r_{12})=0$ thus
\beq c=\beta\ .\label{bez10}\eeq
This explains the value of the constant $c$ in lemma \ref{lemnhacYBe}.

\paragraph{2.} For an endomorphism $r$ of $V\otimes V$, assume that $r\circ r=\beta r_{13}$
and $(Sr)_{12}=\alpha P_{12}+\beta I$. Then
\beq\begin{array}{l}P_{23}(r\circ r-\beta r_{13})P_{23}=r_{13}r_{12}+r_{12}r_{32}-r_{32}r_{13}
-\beta r_{12}\\[.5em] \ \ \ \ \  \ \ \ \ \
=r_{13}r_{12}+r_{12}(\alpha P_{23}+\beta I-r_{23})-(\alpha P_{23}+\beta I-r_{23})r_{13}
-\beta r_{12}=r\circ' r-\beta r_{13}\ .\end{array}\label{bez11}\eeq
Thus, if $(Sr)_{12}=\alpha P_{12}+\beta I$ then $r\circ r=\beta r_{13}$ implies $r\circ' r=\beta 
r_{13}$.

\paragraph{3.} Assume that $r\circ r=cr_{13}$ for an endomorphism $r$ of $V\otimes V$.
Then for $\tilde{r}=r+aI+bP$, $a$ and $b$ are constants, we have
\beq\tilde{r}\circ\tilde{r}=(c+2a)\tilde{r}_{13}+bP_{13}(Sr)_{23}-a(a+c)I-bcP_{13}+b^2P_{23}P_{12}
\ .\label{bez12}\eeq
If, in addition, $(Sr)_{12}=\alpha P_{12}+\beta I$, then
\beq \tilde{r}\circ\tilde{r}=(c+2a)\tilde{r}_{13}-a(c+a)I+b(\beta -c)P_{13}+
b(\alpha +b)P_{23}P_{12}\ .\label{bez13}\eeq
This shows that the equation $r\circ r=c_1 r_{13}+c_2I+c_3P_{13}+c_4 P_{23}P_{12}$, $c_1$, $c_2$, 
$c_3$ and $c_4$ are constants, reduces to  $r\circ r=\tilde{c}_1r_{13}+\tilde{c}_3P_{13}$
by a shift $r\mapsto r+aI+bP$.

If $r\circ r=\beta r_{13}$ and $(Sr)_{12}=\alpha P_{12}+\beta I$ then
\beq \tilde{r}\circ\tilde{r}=(\beta +2a)\tilde{r}_{13}-a(\beta +a)I+
b(\alpha +b)P_{23}P_{12}\ .\label{bez14}\eeq
The combination $P_{23}P_{12}$ does not appear for $b=0$ or $b=-\alpha$. The choice $b=-\alpha$ 
corresponds, modulo a shift of $r$ by a multiple of $I$, to $r\mapsto r_{21}$, so we consider
only $b=0$. Then, with the choice $a=-\beta$ we find that the operator $\tilde{r}=r-\beta I$ satisfies the nhacYBe (and $(Sr)_{12}=\alpha P_{12}-\beta I$). For the choice $a=-\beta /2$ we 
find that the operator $\tilde{r}=r-\displaystyle{\frac{\beta}{2}}I$ satisfies
\beq \tilde{r}\circ\tilde{r}=\frac{\beta^2}{4}\ ,\ (S\tilde{r})_{12}=\alpha P_{12}\ .\label{bez15}\eeq
In particular, the operator
\beq \tilde{\mathfrak{b}}=\frac{x+y}{2(x-y)}I-\frac{x}{x-y}P\label{bez16}\eeq
satisfies (\ref{bez15}) with $\beta =1$ and $\alpha =-1$. Also, $\tilde{\mathfrak{b}}^2=
\displaystyle{\frac{1}{4}}I$.

\paragraph{4.} Assume that $r^2=ur+v$ and $r_{12}+r_{21}=\alpha P_{12}+\beta I$ for an 
endomorphism $r$ of $V\otimes V$. Squaring the relation $r_{12}-\beta I=\alpha P_{12}-r_{21}$ 
and using the same relation again, we obtain
\beq (u-\beta )(2r_{12}-\beta I-\alpha P_{12})=0\ .\label{bez17}\eeq
Thus, if $r$ is not a linear combination of $I$ and $P$ then
\beq u=\beta\ .\label{bez18}\eeq

\paragraph{5.} Assume that $r\circ r=cr_{13}$ and $rP=-r$ for an endomorphism $r$ of $V\otimes V$.
The nhacYBe has the following equivalent form:
\beq [ r_{13},r_{23} ] =(r_{12}-cI)r_{13}P_{23}\ .\label{bez23}\eeq
Indeed,
\beq r_{13}r_{23}-r_{23}r_{13}=(-r_{13}r_{23}+r_{23}r_{12})P_{23}=(r_{12}-cI)r_{13}P_{23}\ .
\label{bez24}\eeq
Here in the first equality we used $r_{23}P_{23}=-r_{23}$ and moved $P_{23}$ to the right;
in the second equality we used the nhacYBe $r\circ r=cr_{13}$.

\subsection{Linear quantization}

Consider an algebra with three generators $r_{12}$, $r_{13}$ and $r_{23}$ and relations
\beq\begin{array}{c}r_{13}r_{23}=r_{23}r_{12}-r_{12}r_{13}+\beta r_{13}\ ,\\[.3em]
r_{13}r_{12}=r_{12}r_{23}-r_{23}r_{13}+\beta r_{13}\ ,\\[.3em]
r_{12}^2=\beta r_{12}+v\ ,\ r_{13}^2=\beta r_{13}+v\ ,\ r_{23}^2=\beta r_{23}+v\ .
\end{array}\label{bez19}\eeq
Choose an order, say, $r_{13}>r_{23}>r_{12}$. Consider (\ref{bez19}) as ordering relations.
The overlaps in (\ref{bez19}) lead to exactly one more relation:
\beq r_{23}r_{12}r_{23}=r_{12}r_{23}r_{12}\ .\label{bez20}\eeq

Thus the algebra in question is 12-dimensional (it follows from (\ref{bez19}) and (\ref{bez20})
that a general element of the algebra is a product $AB$ of an 
element $A$ of the Hecke algebra generated by $r_{12}$ and $r_{23}$ and a polynomial $B$, of 
degree less than 2, in $r_{13}$). 

We conclude that the nhacYBe together with the quadratic equation for $r$ imply the 
YBe. Note that the other form of the YBe also follows:
\beq\begin{array}{l} r_{23}r_{13}r_{12}-r_{12}r_{13}r_{23}=(r_{12}r_{23}-r_{13}r_{12}+
\beta r_{13})r_{12}-r_{12}(r_{23}r_{12}-r_{12}r_{13}+\beta r_{13})r_{23}\\[.3em] 
\ \ \ \ \ \ \ \ \ \ \ \ \ \ \ \ \ \ \ \ \ \ \ \ \ \ \ \ 
=-r_{13}(\beta r_{12}+v)+\beta r_{13}r_{12}+(\beta r_{12}+v)r_{13}-\beta r_{12}r_{13}=0.
\end{array}\label{bez21}\eeq
Here in the first equality both nhacYBe for $r$ were used; the quadratic relation for $r$ was used
in the second equality.

\vskip .2cm
Therefore, the quantization of such $r$-matrix is "linear"\footnote{It was noted in \cite{EH} 
that the operator $\mathfrak{b}^{(0)}$ satisfies both forms of the YBe, squares to zero and
that its quantization has the simple form (\ref{bez22}).}: a combination
\beq R=I+\lambda r\ ,\label{bez22}\eeq
where $\lambda$ is an arbitrary constant, satisfies the YBe $R_{12}R_{13}R_{23}=
R_{23}R_{13}R_{12}$.

\subsection{Algebraic meaning}
 
We shall clarify the algebraic meaning of the non-homogeneous associative classical 
Yang--Baxter equation in the general context of associative algebras.

\vskip .2cm
Let ${\mathfrak{A}}$ be an algebra. Let $r\in 
{\mathfrak{A}}\otimes{\mathfrak{A}}$. The operation
\beq \delta^{(0)} :{\mathfrak{A}}\rightarrow {\mathfrak{A}}\otimes{\mathfrak{A}}\ ,\ 
\delta^{(0)} (u)=(u\otimes 1)\, r-r\, (1\otimes u)\label{um1}\eeq
(the algebra ${\mathfrak{A}}$ does not need to be unital, $(u\otimes 1)(a\otimes b)$
stands for $ua\otimes b$ and $(a\otimes b)(u\otimes 1)$ for $au\otimes b$) is coassociative if 
and only if \cite{A}
\beq (u\otimes 1\otimes 1)\, (r\circ' r)
=(r\circ' r)\, (1\otimes 1\otimes u)\ \ \ \forall\ u\in {\mathfrak{A}}\ .\label{um2}\eeq
In particular, $\delta^{(0)}$ is coassociative if $(r\circ' r)=0$.

\vskip .2cm
Assume now that the algebra ${\mathfrak{A}}$ is unital. Define the operations $\delta$
and $\tilde{\delta}$ $:{\mathfrak{A}}\rightarrow {\mathfrak{A}}\otimes{\mathfrak{A}}$,
\beq \delta (u):=(u\otimes 1)\, r-r\, (1\otimes u)-c\ (u\otimes 1)\ ,\label{um3}\eeq
\beq \tilde{\delta} (u):=(u\otimes 1)\, r-r\, (1\otimes u)+c\ (1\otimes u)\ ,\label{um4}\eeq
where $c$ is a constant. 

\begin{proposition}$\!\!\!${\bf .}\label{uum}
The coassociativity of each of the operations $\delta$ and $\tilde{\delta}$ is equivalent to
\beq (u\otimes 1\otimes 1)\, (r\circ' r-c\, r_{13})
=(r\circ' r-c\, r_{13})\, (1\otimes 1\otimes u)\ \ \ \forall\ u\in {\mathfrak{A}}\ .\label{um5}\eeq
\end{proposition}

\noindent {\bf Proof.} A straightforward calculation.\hfill $\Box$ 

\vskip .2cm
In particular, the operations $\delta$ and $\tilde{\delta}$ are coassociative if 
$r\circ' r=c\, r_{13}$. 

\vskip .2cm
The map (\ref{um1}) has the following property:
\beq \delta^{(0)}(uv)=(u\otimes 1)\, \delta^{(0)}(v)+\delta^{(0)}(u)\, (1\otimes v)\label{um6}
\ ;\eeq
that is, $\delta^{(0)}$ is a derivation with respect to the standard structure of
${\mathfrak{A}}\otimes{\mathfrak{A}}$ as a bi-module over ${\mathfrak{A}}$, $uU:=(u\otimes 1)U$ and 
$Uu:=U(1\otimes u)$ for $u\in {\mathfrak{A}}$ and $U\in {\mathfrak{A}}\otimes{\mathfrak{A}}$.

\vskip .2cm
{}For the operations $\delta$ and $\tilde{\delta}$, the analogue of the property (\ref{um6}) reads
\beq \delta (uv)=(u\otimes 1)\, \delta(v)+\delta(u)\, (1\otimes v)+c\, (u\otimes v)
\ ,\label{um7}\eeq
\beq \tilde{\delta} (uv)=(u\otimes 1)\, \tilde{\delta}(v)+\tilde{\delta}(u)\, (1\otimes v)-
c\, (u\otimes v)\ .\label{um8}\eeq

\subsection{Rota--Baxter operators}

Let ${\cal{A}}$ be an algebra. An operator $\mathfrak{r}: {\cal{A}}\rightarrow {\cal{A}}$ is 
called {\it Rota--Baxter operator of weight} $\alpha$ if 
\beq \mathfrak{r}(A)\mathfrak{r}(B)+\alpha\mathfrak{r}(AB)=\mathfrak{r}
\left(\mathfrak{r}(A)B+A\mathfrak{r}(B)\phantom{\frac{}{}}\right)\ \label{bez25}\eeq
for arbitrary $A,B\in {\cal{A}}$ ($\alpha$ is a constant). We refer to \cite{R} for further
information about the Rota--Baxter operators. 

\vskip .2cm
The Rota--Baxter operators of weight zero are closely related to the acYBe \cite{STS}. It turns 
out that the Rota--Baxter operators of non-zero weight are related to the nhacYBe. We shall 
discuss this relation and calculate the Rota--Baxter operators corresponding to the B\'ezout 
operators. 

\vskip .2cm
It is surprising that the B\'ezout operators, which rather have the sense of derivatives,
become, being interpreted as operators on matrix algebras, the Rota--Baxter operators
which are designed to axiomatize the properties of indefinite integrations and summations.

\paragraph{1.} For an endomorphism $r$ of $V\otimes V$, define two endomorphisms, $\mathfrak{r}$ 
and $\mathfrak{r}'$, of the matrix algebra Mat$(V)$:
\beq \mathfrak{r}(A)_1:=\mathrm{Tr}_2(r_{12}A_2)\ ,\ \mathfrak{r}'(A)_2:=\mathrm{Tr}_1(r_{12}A_1)
\ ,\ A\in \mathrm{Mat}(V)\ ,\label{bez26}\eeq
where $\mathrm{Tr}_i$ is the trace in the copy number $i$ of the space $V$.

\vskip .2cm
Assume that $r$ satisfies the nhacYBe (\ref{bez8}). Multiplying (\ref{bez8}) by $A_2B_3$,
$A,B\in \mathrm{Mat}(V)$, and taking traces in the spaces 2 and 3, we find
\beq \mathfrak{r}(A)\mathfrak{r}(B)+\mathfrak{r}\left(\mathfrak{r}'(A)B\phantom{\frac{}{}}\right)
-\mathfrak{r}\left( A\mathfrak{r}(B)\phantom{\frac{}{}}\right)=c\,\mathrm{Tr}(A)\mathfrak{r}(B)
\ .\label{bez27}\eeq
Assume, in addition, that $r_{12}+r_{21}=\alpha P_{12}+\beta I$. Then
\beq \mathfrak{r}(A)+\mathfrak{r}'(A)=\alpha A+\beta \mathrm{Tr}(A)\, 1\!\!1\ .\label{bez28}\eeq
If $c=\beta$ then expressing $\mathfrak{r}'(A)$ by (\ref{bez28}) and substituting into
(\ref{bez27}), we find that the term with $\mathrm{Tr}(A)$ drops out and $\mathfrak{r}$ is the
Rota--Baxter operator of weight $\alpha$ on the algebra of matrices. Similarly, $\mathfrak{r}'$ is 
the Rota--Baxter operator of weight $\alpha$ as well.

\paragraph{2.} We shall calculate the Rota--Baxter operators corresponding to the B\'ezout 
operators in the polynomial basis. 

\vskip .2cm
The action of the operator $\mathfrak{b}^{(0)}$
on monomials $x^ky^l$ reads
\beq \mathfrak{b}^{(0)}(x^ky^l)=\left\{\begin{array}{l}
-(x^{l-1}y^k+x^{l-2}y^{k+1}+\dots +x^ky^{l-1})\ \ ,\ k<l\ ,\\[.5em]
0\ \ ,\ k=l\ ,\\[.5em]
x^{k-1}y^l+x^{k-2}y^{l+1}+\dots +x^ly^{k-1}\ \ ,\ k>l\ .
\end{array}\right. \label{brr1}\eeq
The action of the operator $\mathfrak{b}$ on monomials $x^ky^l$ reads
\beq \mathfrak{b}(x^ky^l)=\left\{\begin{array}{l}
-(x^{l}y^k+x^{l-1}y^{k+1}+\dots +x^{k+1}y^{l-1})\ \ ,\ k<l\ ,\\[.5em]
0\ \ ,\ k=l\ ,\\[.5em]
x^{k}y^l+x^{k-1}y^{l+1}+\dots +x^{l+1}y^{k-1}\ \ ,\ k>l\ .
\end{array}\right. \label{brr2}\eeq
Shortly,
\beq \mathfrak{b}^{(0)}(x^ky^l)=\theta (k-l)\sum_{s=0}^{k-l-1}x^{l+s}y^{k-s-1}-
\theta (l-k)\sum_{s=0}^{l-k-1}x^{k+s}y^{l-s-1} \ ,\label{brr1p}\eeq
\beq \mathfrak{b}(x^ky^l)=\theta (k-l)\sum_{s=1}^{k-l}x^{l+s}y^{k-s}-
\theta (l-k)\sum_{s=1}^{l-k}x^{k+s}y^{l-s} \ .\label{brr2p}\eeq

We list several useful matrix forms of the operators $\mathfrak{b}^{(0)}$ and $\mathfrak{b}$ 
in the basis formed by monomials, $e_a\otimes e_b:=x^ay^b$; for the operator $\mathfrak{b}^{(0)}$:
\beq\begin{array}{lcl} \mathfrak{b}^{(0)}&=&\displaystyle{\sum_{i,j,a,b}}\ \theta (j-a)\,
\theta (j-b)\ \delta^{i+j}_{a+b+1}\ 
e^j_a\wedge e^i_b\\[2em]
&=&\displaystyle{\sum_{i,j,a,b}}\ \bigl( \theta (j-a)\,\theta (j-b)-
\theta (i-b)\,\theta (i-a)\bigr)\ \delta^{i+j}_{a+b+1}\ e^j_a\otimes e^i_b\\[1em]
&=&\displaystyle{\sum_{i,j:i<j}}\ \ \displaystyle{\sum_{a=1}^{j-i}}\ \ e^j_{i+a-1}\wedge e^i_{j-a}
\ \end{array}\label{brr3}\eeq
and for the operator $\mathfrak{b}$:
\beq\begin{array}{lcl} \mathfrak{b}&=&\displaystyle{\sum_{i,j,a,b}}\ \theta (j+1-a)\,
\theta (a-i)\ \delta^{i+j}_{a+b}\ 
(e^j_a\otimes e^i_b-e^i_a\otimes e^j_b)\\[2em]
&=&\displaystyle{\sum_{i,j,a,b}}\ \bigl( \theta (j+1-a)\,\theta (a-i)-
\theta (i+1-a)\,\theta (a-j)\bigr)\ \delta^{i+j}_{a+b}\ e^j_a\otimes e^i_b\\[1em]
&=&\displaystyle{\sum_{i,j:i<j}}\ \ \displaystyle{\sum_{a=1}^{j-i}}\ \  
(e^j_{i+a}\otimes e^i_{j-a}-e^i_{i+a}\otimes e^j_{j-a})\\[1em]
&=&\displaystyle{\sum_{i,j:i<j}}\ \ \Bigl(\displaystyle{\sum_{a=1}^{j-i-1}}\ 
e^j_{i+a}\wedge e^i_{j-a}\, +\, e^j_j\otimes e^i_i\, -\, e^i_j\otimes e^j_i\Bigr)\ ,
\end{array}\label{brr4}\eeq
where $x\wedge y=x\otimes y-y\otimes x$.

\vskip .2cm
The Rota--Baxter operator $\mathfrak{r}_{\mathfrak{b}^{(0)}}$ corresponding to ${\mathfrak{b}^{(0)}}$ reads
\beq \mathfrak{r}_{ \mathfrak{b}^{(0)}}(A)^i_j=\theta (j-i)\sum_{s\geq 0}
A^{i-s}_{j-s-1}-\theta (i+1-j)\sum_{s\geq 0}
A^{i+s+1}_{j+s}\ .\label{brr5}\eeq
In the right hand side of (\ref{brr5}), the summations are over those $s\geq 0$ for which
the corresponding matrix element in the sum makes sense; that is, the range of $s$ in the first 
sum is $s=0,1,\dots ,i-1$ and, in the second sum, $s=0,1,\dots n-i-1$; 

\vskip .2cm
The Rota--Baxter operator $\mathfrak{r}_{\mathfrak{b}}$ corresponding to
${\mathfrak{b}}$ reads (with the same convention about the summation ranges)
\beq \mathfrak{r}_{ \mathfrak{b}}(A)^i_j=\theta (j+1-i)\sum_{s\geq 0}
A^{i-s-1}_{j-s-1}-\theta (i-j)\sum_{s\geq 0}
A^{i+s}_{j+s}\ .\label{brr6}\eeq
Its weight is -1.

\vskip .2cm
For the operator $r^{(s)}$, given by eq. (\ref{bez3}), the corresponding Rota--Baxter 
operator 
$\mathfrak{r}^{(s)}$ is
\beq \mathfrak{r}^{(s)}
(A)^i_j=\left\{\begin{array}{cc}-\theta(j-i)A^i_j\ ,&i\neq j\ ,\\[.4em]
\displaystyle{\sum_{s:s<i}}A^s_s\ ,&i=j\ .\end{array}\right.\label{bez29}\eeq
Its weight is -1.

\vskip .2cm
We shall give also the Rota--Baxter operator for the B\'ezout $r$-matrix $\mathfrak{b}$ in the rime
basis, that is, for the $r$-matrix (\ref{rcb}); it has weight 1 (since $r_{12}+r_{21}=P-I$
for $r$ in (\ref{rcb})). The Rota--Baxter operator has the form
\beq \mathfrak{r} (A)^i_j=\left\{\begin{array}{cc}
\displaystyle{\frac{\phi_j}{\phi_j-\phi_i}}(A^i_j-A^j_j)\ ,&i\neq j\ ,\\[.8em]
\displaystyle{\sum_{s:s\neq i}\,\frac{\phi_i}{\phi_i-\phi_s}} (A^i_s-A^s_s)\ ,&i=j
\ .\end{array}\right.\label{bez30}\eeq

\subsection{$*$-multiplication}

\paragraph{1.} Let $\mathfrak{r}: {\cal{A}}\rightarrow {\cal{A}}$ be a Rota--Baxter operator of 
weight $\alpha$ (see eq.(\ref{bez25})) on an algebra ${\cal{A}}$. It is known that the operation
\beq A*B:=\mathfrak{r}(A)B+A\mathfrak{r}(B)-\alpha\, AB\ \ ,\ A,B\in\, {\cal{A}}\ ,\label{stm1}\eeq
defines an associative product on the space ${\cal{A}}$. This product is closely related to the
coproducts (\ref{um3}) and (\ref{um4}) by duality. We shall illustrate it in the context
of the matrix algebras. 

\vskip .2cm
Define an operation $\tilde{*}$ by
\beq \langle \tilde{\delta}(u), B\otimes A\rangle =\langle u,A\tilde{*}B\rangle\ ,\label{stm2}\eeq
where $\tilde{\delta}$ is given by (\ref{um4}). We have then
\beq\begin{array}{lcl} \langle \tilde{\delta}(u), B\otimes A\rangle &=&
\mathrm{Tr}_{12}\Bigl( u_1rB_1A_2-ru_2B_1A_2+c\, u_2B_1A_2\Bigr)\\[1em]
&=&\mathrm{Tr}_{1}\Bigl( u_1\, \mathrm{Tr}_2(rA_2)\, B_1\Bigr)-
\mathrm{Tr}_{1}\Bigl( u_1A_1\, \mathrm{Tr}_2(B_2\, r_{21})\Bigr) +
\mathrm{Tr}_1\Bigl( c\, u_1\, \mathrm{Tr}(B)\, A_1\Bigr)\\[1em]
&=&\mathrm{Tr}\Bigl( u\, \bigl[ \mathfrak{r}(A)\, B-A\, \mathfrak{r}'(B)+
c\, A\, \mathrm{Tr}(B)\bigr]\Bigr)\ ,\end{array}\label{stm3}\eeq
thus
\beq A\tilde{*}B=\mathfrak{r}(A)\, B-A\, \mathfrak{r}'(B)+c\, A\,\mathrm{Tr}(B)\ .\label{stm4}\eeq
In eq. (\ref{stm3}), $x_i$ stands for the copy of an element $x$ in the space number $i$ in 
${\cal{A}}\otimes {\cal{A}}$; the operators $\mathfrak{r}$ and $\mathfrak{r}'$ are given by 
(\ref{bez26}); to obtain the second and the third terms in the second line of (\ref{stm3}) we 
renumbered $1\leftrightarrow 2$ and then moved $r$ cyclically under the trace in the second term. 

\vskip .2cm
Assume, as before, that $r_{12}+r_{21}=\alpha P_{12}+\beta I$ and $c=\beta $. Then, expressing 
$\mathfrak{r}'(A)$ by (\ref{bez28}), we find that the term with $\mathrm{Tr}(B)$ drops out and 
it follows that 
\beq A\tilde{*}B=A*B\ .\label{stm5}\eeq

\paragraph{2.} We shall describe the $*$-multiplication in the simplest example of the Rota--Baxter operators \ref{brr5}) and (\ref{brr6}) corresponding to the B\'ezout operators for 
the the polynomials of degree less than 2 (that is, for $2\times 2$ matrices
$A=\left(\begin{array}{cc}a^1_1& a^1_2\\[.3em] a^2_1& a^2_2\end{array}\right)\equiv a^i_je^j_i$).

\vskip .2cm
For the operator $\mathfrak{b}^{(0)}=e^2_1\wedge e^1_1$, we have
\beq \mathfrak{r}_{\mathfrak{b}^{(0)}}(A)=\left(\begin{array}{cc}-a^2_1&a^1_1\\ 0&0
\end{array}\right)\label{stmn1}\eeq
and the $*$-multiplication reads
\beq A*^o\tilde{A}\equiv A\mathfrak{r}_{\mathfrak{b}^{(0)}}(\tilde{A})+\mathfrak{r}_{\mathfrak{b}^{(0)}}(A)\tilde{A}=
\left(\begin{array}{cc}-a^2_1\tilde{a}^1_1& -a^2_1\tilde{a}^1_2+
a^1_1(\tilde{a}^1_1+\tilde{a}^2_2)\\[.5em] 
-a^2_1\tilde{a}^2_1& a^2_1\tilde{a}^1_1\end{array}\right)\ .\label{stmn2}\eeq
This algebra is isomorphic to the algebra of $3\times 3$ matrices of the form
$$\left(\begin{array}{ccc}*&*&*\\ 0&0&*\\ 0&0&0\end{array}\right)\ ,$$
with the identification
\beq e^1_1\mapsto \left(\begin{array}{ccc}0&1&0\\ 0&0&1\\ 0&0&0\end{array}\right) ,\ 
e^1_2\mapsto \left(\begin{array}{ccc}-1&0&0\\ 0&0&0\\ 0&0&0\end{array}\right) ,\ 
e^2_1\mapsto \left(\begin{array}{ccc}0&0&1\\ 0&0&0\\ 0&0&0\end{array}\right) ,\
e^2_2\mapsto \left(\begin{array}{ccc}0&0&0\\ 0&0&1\\ 0&0&0\end{array}\right) .\label{stmn3}\eeq

For the operator $\mathfrak{b}=e^2_2\otimes e^1_1-e^1_2\otimes e^2_1$, we have
\beq \mathfrak{r}_{\mathfrak{b}}(A)=\left(\begin{array}{cc}0&0\\ -a^2_1&a^1_1
\end{array}\right)\label{stmn4}\eeq
and the $*$-multiplication reads
\beq A*\tilde{A}\equiv A\mathfrak{r}_{\mathfrak{b}}(\tilde{A})+\mathfrak{r}_{\mathfrak{b}}(A)\tilde{A}+A\tilde{A}=
\left(\begin{array}{cc}a^1_1\tilde{a}^1_1& a^1_1\tilde{a}^1_2+
a^1_2(\tilde{a}^1_1+\tilde{a}^2_2)\\[.5em] 
a^1_1\tilde{a}^2_1& a^1_1\tilde{a}^2_2+
a^2_2(\tilde{a}^1_1+\tilde{a}^2_2)\end{array}\right)
\ .\label{stmn5}\eeq
This algebra is isomorphic to the algebra of $3\times 3$ matrices of the form
$$\left(\begin{array}{ccc}*&*&*\\ 0&*&0\\ 0&0&0\end{array}\right)\ ,$$
with the identification
\beq e^1_1\mapsto \left(\begin{array}{ccc}1&0&0\\ 0&1&0\\ 0&0&0\end{array}\right) ,\ 
e^1_2\mapsto \left(\begin{array}{ccc}0&0&1\\ 0&0&0\\ 0&0&0\end{array}\right) ,\ 
e^2_1\mapsto \left(\begin{array}{ccc}0&1&0\\ 0&0&0\\ 0&0&0\end{array}\right) ,\
e^2_2\mapsto \left(\begin{array}{ccc}0&0&0\\ 0&1&0\\ 0&0&0\end{array}\right) .\label{stmn6}\eeq

\section{Rime Poisson brackets\vspace{.25cm}}
 
The Poisson brackets having the form
\beq \{ x^i,x^j\} =f_{ij}(x^i,x^j)\quad ,\quad i,j=1,2,\dots ,n\ ,
\label{rpb0}\eeq
with some functions $f_{ij}$ of two variables, we shall call {\it rime}. In this section we study quadratic rime Poisson brackets, 
\beq \{ x^i,x^j\} =a_{ij}(x^i)^2-a_{ji}(x^j)^2+2\nu_{ij}x^ix^j\quad ,\quad i,j=1,2,\dots ,n\ .
\label{rpb1}\eeq
We show that there is a three-dimensional pencil of such Poisson brackets and then find
the invariance group and the normal form of each individual member of the pencil.

\subsection{Rime pencil}

In this subsection we establish that the quadratic rime Poisson brackets form a three-dimensional 
Poisson pencil.

\vskip .2cm
The left hand side of (\ref{rpb1}) contains a matrix $a_{ij}$ with zeros on the diagonal, 
$a_{ii}=0$, and an anti-symmetric matrix $\nu_{ij}$, $\nu_{ij}=-\nu_{ji}$. The Jacobi identity 
constraints these matrices to satisfy
\beq a_{ij}a_{jk}+a_{ik}(\nu_{ij}+\nu_{jk})=0\ ,\ i\neq j\neq k\neq i\ .\label{rpb2}\eeq
We shall describe a general solution of eq. (\ref{rpb2}) in the strict situation, that is,
when all $a_{ij}$ and $\nu_{ij}$ are different from zero for $i\neq j$. 

The left hand side of $\nu_{ij}+\nu_{jk}=-a_{ij}a_{jk}/a_{ik}$ is anti-symmetric with respect to $(i,k)$, that is $\Upsilon_{ij}\Upsilon_{jk}\Upsilon_{ki}=1$ for $\Upsilon_{ij}=-a_{ij}/a_{ji}$, 
which readily implies the existence of a vector $\phi_i$ such that $\Upsilon_{ij}=
\phi_i^2/\phi_j^2$. Therefore,
\beq a_{ik}=\phi_ic_{ik}\phi_k^{-1}\ , \label{rpb3}\eeq
where the matrix $c_{ij}$ is anti-symmetric, $c_{ij}=-c_{ji}$. Next, 
$2\nu_{ki}=-(\nu_{ij}+\nu_{jk})+(\nu_{jk}+\nu_{ki})+(\nu_{ki}+\nu_{ij})$; using (\ref{rpb2}) to 
express each bracket in the right hand side, we find 
\beq\nu_{ki}=\frac{1}{2}\left(\frac{c_{ij}c_{ki}}{c_{jk}}+\frac{c_{jk}c_{ki}}{c_{ij}}-
\frac{c_{ij}c_{jk}}{c_{ki}}\right)\ .\label{rpb4}\eeq
The right hand side of eq. (\ref{rpb4}) does not depend on $j$ which imposes further restrictions 
on the matrix $c_{ij}$ when $n>3$. Writing the sum $\nu_{ij}+\nu_{jk}+\nu_{kl}+\nu_{li}$ in two 
ways, as $(\nu_{ij}+\nu_{jk})+(\nu_{kl}+\nu_{li})$ and as $(\nu_{jk}+\nu_{kl})+(\nu_{li}+
\nu_{ij})$, and using (\ref{rpb2}) to express each bracket in terms of the matrix $c$, we obtain
\beq\frac{c_{ij}c_{jk}-c_{il}c_{lk}}{c_{ik}}=\frac{c_{jk}c_{kl}-c_{ji}c_{il}}{c_{jl}}\ 
.\label{rpb5}\eeq 
Replacing $j$ by $m$ in (\ref{rpb4}) gives the condition on the matrix $c$:
\beq \frac{c_{ij}c_{ki}}{c_{jk}}+\frac{c_{jk}c_{ki}}{c_{ij}}-
\frac{c_{ij}c_{jk}}{c_{ki}}=\frac{c_{im}c_{ki}}{c_{mk}}+\frac{c_{mk}c_{ki}}{c_{im}}-
\frac{c_{im}c_{mk}}{c_{ki}}\ .\label{rpb6}\eeq
Using eq. (\ref{rpb5}) to rewrite the combination $\displaystyle{\frac{c_{ij}c_{jk}}{c_{ki}}}-
\displaystyle{\frac{c_{im}c_{mk}}{c_{ki}}}$, we find
\beq (c_{jk}c_{km}-c_{ji}c_{im})\Psi_{ijkm}=0\ \ ,\ \ {\mathrm{where}}\ \ \ \ 
\Psi_{ijkm}=\left(\frac{1}{c_{jk}c_{im}}+\frac{1}{c_{ij}c_{km}}+
\frac{1}{c_{ki}c_{jm}}\right) \ .\label{rpb7}\eeq
The quantity $\Psi_{ijkm}$ is totally anti-symmetric with respect to its indices. Therefore, if 
$\Psi_{ijkm}\neq 0$ then the combinations $(c_{jk}c_{km}-c_{ji}c_{im})$ vanish for all 
permutations of indices. This is however impossible: the system of three linear equations 
\beq\begin{array}{l} c_{ij}c_{jk}-c_{im}c_{mk}=0\ ,\\
c_{ik}c_{km}-c_{ij}c_{jm}=0\ ,\\
c_{im}c_{mj}-c_{ik}c_{kj}=0\ \end{array}\eeq
for unknowns $\{ c_{jk},c_{km},c_{mj}\}$ has, by definition, a non-zero solution but the 
determinant of the system is different from zero. Thus the Pfaffian $\Psi_{ijkm}$
vanishes for each quadruple $\{ i,j,k,m\}$; in other words, the coefficients
of the matrix $1/c_{ij}$ satisfy the Pl\"ucker relations; therefore the form $1/c_{ij}$ is
decomposable, $c_{ij}^{-1}=s_it_j-s_jt_i$,
for some vectors $\vec{s}$ and $\vec{t}$. For each $i$, at least one of $s_i$ or $t_i$ is 
different from zero. Making, if necessary, a change of basis in the two dimensional plane
spanned by $\vec{s}$ and $\vec{t}$, we can therefore always assume that all components of, say,
the vector $\vec{s}$ are different from zero, $s_i\neq 0$ $\forall$ $i$. We represent the
bivector $1/c_{ij}$ in the form ($u_i^{-1}=s_i$ and $\psi_i=-t_i/s_i$)
\beq \frac{1}{c_{ij}}=u_i^{-1}u_j^{-1}(\psi_i-\psi_j)\ .\label{rpb9}\eeq
Substituting (\ref{rpb9}) into (\ref{rpb4}) we obtain
\beq \nu_{ki}+\frac{1}{2}\frac{u_k^2+u_i^2}{\psi_k-\psi_i}=-\frac{1}{2}\left(
\frac{u_i^2-u_j^2}{\psi_i-\psi_j}-\frac{u_k^2-u_j^2}{\psi_k-\psi_j}\right)\ .\label{rpb10}\eeq
Replacing $j$ by $m$ in the right hand side and equating the resulting expressions, we find that 
the independency of the right hand side on $j$ implies:
\beq\begin{array}{rcl} E_{ijkm}&:=&{\displaystyle{\frac{u_i^2}{(\psi_i-\psi_j)(\psi_i-\psi_k)(\psi_i-\psi_m)}}}+
{\displaystyle{\frac{u_j^2}{(\psi_j-\psi_i)(\psi_j-\psi_k)(\psi_j-\psi_m)}}}\\[1em]
&+& {\displaystyle{\frac{u_k^2}{(\psi_k-\psi_i)(\psi_k-\psi_j)(\psi_k-\psi_m)}}}+
{\displaystyle{\frac{u_m^2}{(\psi_m-\psi_i)(\psi_m-\psi_j)(\psi_m-\psi_k)}}}=0\ . \end{array} 
\label{rpb11}\eeq
for every quadruple $\{ i,j,k,m\}$.

The quantity $E_{ijkm}$ is totally symmetric. Selecting three values of the index, say, 1,2
and 3, we can form the quadruple $\{ i,1,2,3\}$ for each $i$. Solving $ E_{i123}=0$, we obtain
the following expression for $u_i^2$:
\beq\begin{array}{rcl} 
u_i^2&=&A_1{\displaystyle{\frac{(\psi_i-M_2)(\psi_i-M_3)}{(M_1-M_2)(M_1-M_3)}}}
+A_2{\displaystyle{\frac{(\psi_i-M_1)(\psi_i-M_3)}{(M_2-M_1)(M_2-M_3)}}}\\[1em]
&+&A_3{\displaystyle{\frac{(\psi_i-M_1)(\psi_i-M_2)}{(M_3-M_1)(M_3-M_2)}}}\ 
\end{array}\label{rpb13}\eeq
for some constants $A_1,A_2,A_3,M_1,M_2$ and $M_3$. The right hand side is the value,
at the point $\psi_i$, of a quadratic polynomial which equals to $A_a$ at the points $M_a$,
$a=1,2,3$. Since $A_1,A_2,A_3,M_1,M_2$ and $M_3$ are arbitrary, we can simply write
\beq u_i^2=a\psi_i^2+b\psi_i+c\ .\label{rpb14}\eeq
With the expressions (\ref{rpb14}) for $u_i^2$, the equalities (\ref{rpb11}) are identically
satisfied which shows that (\ref{rpb14}) is the general solution.

Upon rescaling $x^i\mapsto \phi_i u_ix^i$ with $\phi_i$ from (\ref{rpb3}), the Poisson brackets 
(\ref{rpb1}) simplify. The following statement is established (for $n=2$ or 3, (\ref{rpb14}) does 
not impose a restriction on the anti-symmetric matrix $c_{ij}$ with all off-diagonal entries 
different from zero).

\begin{proposition}$\!\!\!${\bf .}\label{rimepb} Up to a rescaling of variables, the general 
strict quadratic rime Poisson brackets have the form
\beq \{ x^i,x^j\} =\frac{\varrho (\psi_j)(x^i)^2+\varrho (\psi_i)(x^j)^2}{\psi_i-\psi_j}
+\left( (\psi_i-\psi_j)\, a-\frac{\varrho (\psi_i)+\varrho 
(\psi_j)}{\psi_i-\psi_j}\right) x^ix^j\ ,\label{rpb15}\eeq
where $\vec{\psi}$ is an arbitrary vector with pairwise distinct components and $\ \varrho (t)=
at^2+bt+c\ $ is an arbitrary quadratic polynomial\,\footnote{To have nonvanishing coefficients in 
the formula (\ref{rpb15}) one has to impose certain inequalities for the components of the vector 
$\vec{\psi}$ and the coefficients of the polynomial $\varrho$; however, the formula
(\ref{rpb15}) defines Poisson brackets without these inequalities.}.\end{proposition}

Thus the strict quadratic rime Poisson brackets form the three-dimensional pencil
(parameterized by the polynomial $\varrho$).

\vskip .2cm
The Poisson brackets (\ref{rpb15}) can be rewritten in the following forms:
\beqa &\{ x^i,x^j\} =\displaystyle{\frac{1}{\psi_i-\psi_j}}\, 
\left( \varrho (\psi_j)x^i-\varrho (\psi_i)x^j\phantom{\frac{}{}}\!\!\right) (x^i-x^j)+
a\, (\psi_i-\psi_j)\, x^ix^j\ ,&\label{rpb16}\\
&\{ x^i,x^j\} =\displaystyle{\frac{au_{ij}^2+bu_{ij}v_{ij}+cv_{ij}^2}{\psi_i-\psi_j}}
\equiv\displaystyle{\frac{v_{ij}^2}{\psi_i-\psi_j}}\ \varrho\left(\displaystyle{\frac{u_{ij}}{v_{ij}}}\right)\, ,&\label{rpb17}\eeqa
where $u_{ij}=\psi_jx^i-\psi_ix^j$ and $v_{ij}=x^i-x^j$.

\vskip .2cm
\noindent {\bf Remark 1.}
For $\varrho (t)=b t$ (respectively, $\varrho (t)=c$) these Poisson brackets appear in the 
classical limit of the commutation relations (\ref{qp}) in the non-unitary (respectively, unitary) 
case (with the parameterization  $\beta_{ij}=-\displaystyle{\frac{\beta\psi_j}{\psi_i-\psi_j}}$
in the non-unitary case).

\vskip .2cm
\noindent {\bf Remark 2.} The strict rime linear Poisson brackets 
\beq \{ x^i,x^j\} =a_{ij}x^i-a_{ji}x^j\ ,\ \ \ a_{ij}\neq 0\ \ \ {\mathrm{for}}\ {\mathrm{all}}\ 
\ i,j=1,2,\dots,n:\ i\neq j\eeq
(or strict rime Lie algebras) are less interesting. The Jacobi identity is
\beq\label{jsla}a_{ik}a_{kj}=a_{ij}a_{jk}\ \ \ {\mathrm{for}}\ {\mathrm{all}}
\ i\neq j\neq k\neq i\ .\eeq 
Rescale variables $x^2,x^3,\dots,x^n$ to have $a_{1i}=1$, $i=2,\dots,n$. Then the condition
(\ref{jsla}) with one of $i,j,k$ equal 1 implies $a_{ij}=a_{ji}$ and $a_{ij}=a_{i1}/a_{j1}$, 
$i,j=2,\dots,n$; it follows that $a_{i1}^2=a_{j1}^2$, $i,j=2,\dots,n$. For $n>3$, the condition
(\ref{jsla}) with $i,j,k>1$ forces $a_{i1}=a_{j1}$, $i,j=2,\dots,n$. Denote by $\nu$ this common 
value, $a_{i1}=\nu$, $i,j=2,\dots,n$. After a rescaling $x^1\mapsto \nu x^1$ we find a unique strict rime Lie algebra, $ [ x^i,x^j ] =x^i-x^j$ for all $i$ and $k$, which is almost trivial: 
$ [ x^i,x^k-x^l ] =-(x^k-x^l)$ for all $i,k$ and $l$ and $ [ x^i-x^j,x^k-x^l ] =0$ for all $i,j,k,l$.

\vskip .2cm
{}For $n=3$, there is one more possibility: $a_{31}=-a_{21}$. After a rescaling 
$x^1\mapsto a_{21} x^1$, the solution reads
\beq [x^1,x^2]=x^1-x^2\ ,\ [x^1,x^3]=x^1+x^3\ ,\ [x^2,x^3]=-x^2+x^3\ .\eeq
This Lie algebra is isomorphic to $\mathfrak{sl}(2)$; the isomorphism is given, for example, by
$h\mapsto x^1-x^3$, $e\mapsto x^1+x^3$ and $f\mapsto x^2-(x^1+x^3)/4$ (here $h,e$ and $f$ are the 
standard generators of $\mathfrak{sl}(2)$, $[h,e]=2e$, $[h,f]=-2f$ and $[e,f]=h$).

\subsection{Invariance}\label{sinv}

In this subsection we analyze the invariance group of each individual member of the Poisson pencil
from the proposition \ref{rimepb}. We find that the Poisson brackets (\ref{rpb15}), with arbitrary
(non-vanishing) $\varrho$, admit a non-trivial 1-parametric invariance group.

\vskip .2cm
The transformation law of Poisson brackets $\{ x^i,x^j\} =f^{ij}(x)$ under an infinitesimal change 
of variables, $\tilde{x}^i=x^i+\epsilon\,\varphi^i(x)$, $\epsilon^2=0$, is 
$\{ \tilde{x}^i,\tilde{x}^j \} =f^{ij}(\tilde{x})+\epsilon\,\delta_x f^{ij}$,
where $\delta_x f^{ij}=\{ \varphi^i,x^j\} +\{ x^i,\varphi^j\} +\varphi^k\partial_kf^{ij}$.
For a linear infinitesimal transformation, $\varphi^i(x)=A^i_jx^j$, we have 
\beq\label{ri1}\delta_x f^{ij}=A^i_k\{ x^k,x^j\} +A^j_k\{ x^i, x^k \} -
x^lA^k_l\partial_k\{ x^i,x^j\}\ .\eeq
Specializing to the Poisson brackets (\ref{rpb15}), we find
\beq \delta_x f^{ij}=U_{ji}-U_{ij}\ \label{ri1p}\eeq
with
\beq\label{ri2}\begin{array}{l}
U_{ij}:=\displaystyle{\sum_s}\left( 2A^i_s\displaystyle{\frac{\varrho_j}{\psi_{ij}}}+
A^j_s(\psi_{ij}a-\displaystyle{\frac{\varrho_i+\varrho_j}{\psi_{ij}}})\right) x^ix^s\\[1em]
\ \ \ \ \ \ +\displaystyle{\sum_{s:s\neq i}}A^j_s\left(\displaystyle{\frac{\varrho_i}{\psi_{si}}}(x^s)^2+
\displaystyle{\frac{\varrho_s}{\psi_{si}}}(x^i)^2+
(\psi_{si}a-\displaystyle{\frac{\varrho_s+\varrho_i}{\psi_{si}}})x^ix^s\right)
\ ,\end{array}\eeq
where $\psi_{ij}=\psi_i-\psi_j$ and $\varrho_s=\varrho (\psi_s)$.

\vskip .2cm
The Poisson brackets (\ref{rpb15}) remain rime under the infinitesimal linear transformation
with the matrix $A$ if the coefficients in $(x^s)^2$, $x^sx^i$ and
$x^sx^j$, $s\neq i,j$, in (\ref{ri1p}) vanish which gives the following system:
\begin{eqnarray}
\label{ris3}(x^s)^2\ ,\ s\neq i,j &\Rightarrow&A^j_s\frac{\varrho_i}{\psi_{si}}
-A^i_s\frac{\varrho_j}{\psi_{sj}}=0\ ,\\
\label{ris2}x^ix^s\ ,\ s\neq i,j &\Rightarrow&2A^i_s\frac{\varrho_j}{\psi_{ij}}+A^j_s\left(\psi_{sj}a-
\frac{\varrho_i+\varrho_j}{\psi_{ij}}-\frac{\varrho_s+\varrho_i}{\psi_{si}}\right) =0\ .
\end{eqnarray}
Eq. (\ref{ris3}) implies that $A^l_k=\nu_k \varrho_l/\psi_{lk}$, $l\neq k$, with arbitrary 
constants $\nu_k$. For a quadratic polynomial $\varrho$, this solves eq. (\ref{ris2}) as well. The 
coefficient in $x^jx^s$ vanishes due to the anti-symmetry.

\vskip .2cm
The Poisson brackets (\ref{rpb15}) are invariant under the infinitesimal linear transformation
with the matrix $A$ if, in addition to (\ref{ris3}) and (\ref{ris2}), the coefficients in 
$(x^i)^2$, $(x^j)^2$ and  $x^ix^j$ in (\ref{ri1p}) vanish which gives:
\begin{eqnarray}
\label{ris1}x^ix^j&\Rightarrow&\varrho_jA^i_j+\varrho_iA^j_i=0\ ,\\
\label{ris4}(x^i)^2&\Rightarrow&A^i_i\frac{\varrho_j}{\psi_{ij}}+A^j_i\left(\psi_{ij}a-
\frac{\varrho_i+\varrho_j}{\psi_{ij}}\right) +\sum_{s:s\neq i}A^j_s\frac{\varrho_s}{\psi_{si}}=0\ .\end{eqnarray}
Eq. (\ref{ris1}) implies that $\nu_k$ are equal, $\nu_k=\nu$. The matrix $A$ is defined up
to a multiplicative factor, so we can set $\nu$ to 1. Since the Poisson brackets (\ref{rpb15}) 
are quadratic, a global rescaling leaves them invariant, so we can add to $A$ a matrix,
proportional to the identity matrix and make $A$ traceless. The traceless condition, together
with eq. (\ref{ris4}) determines the diagonal entries, $A^i_i=a(n-1)\psi_i+
\displaystyle{\frac{n-1}{2}b+\varrho_i\sum_{s:s\neq i}\frac{1}{\psi_{si}}}$. The coefficient in 
$(x^j)^2$ vanishes due to the anti-symmetry. We summarize the obtained results.

\begin{proposition}$\!\!\!${\bf .} (i) The infinitesimal linear transformation with the matrix 
$A$ leaves the Poisson brackets (\ref{rpb15}) rime if and only if
\beq A^l_k=\frac{\nu_k \varrho_l}{\psi_{lk}}\ ,\ \ l\neq k\ ,\label{ris5}\eeq
with arbitrary constants $\nu_k$.

\vskip .2cm
(ii) Up to a global rescaling of coordinates, the invariance group of the Poisson brackets 
(\ref{rpb15}) is 1-parametric, with a generator ${\cal{A}}$,
\beq {\cal{A}}(\varrho )^i_j=\frac{\varrho_i}{\psi_{ij}}\ ,\ \ i\neq j\ ,\ \ \ {\mathrm{and}}\ \ \ \ 
{\cal{A}}(\varrho )^i_i=\displaystyle{\frac{n-1}{2}\varrho_i'+\varrho_i\xi_i\ ,\ \ \xi_i:=
\sum_{s:s\neq i}\frac{1}{\psi_{si}}}\ ,\label{ris6}\eeq
where $\varrho_i'$ is the value of the derivative of the polynomial $\varrho$ at the point $\psi_i$.
\end{proposition}

Since the Poisson brackets transformed with the matrix (\ref{ris5}) are still rime, it follows 
from the proposition \ref{rimepb} that they can be written, after an appropriate rescalings of coordinates, in the form (\ref{rpb15}). In other words, the variation $\delta_x$ can be compensated by a variation of $\psi$'s and $\varrho$ and a diagonal transformation of the
coordinates. We have 
\beq -\delta_x f^{ij}=\delta^{(1)}+\delta^{(2)}\ ,\label{ris7}\eeq
where
\beq \delta^{(1)}=\frac{\varrho_i\varrho_j(x^i-x^j)^2}{\psi_{ij}^2}(\nu_i-\nu_j)
+a\left( \nu_j\varrho_j(x^i)^2-\nu_i\varrho_i(x^j)^2\phantom{\frac{}{}}\!\! \right)\label{ris8}\eeq
and 
\beq\delta^{(2)}=(\tilde{A}^i_i-\tilde{A}^j_j)\frac{\varrho_j(x^i)^2-\varrho_i(x^j)^2}{\psi_{ij}}\ ,\ \ 
\tilde{A}^i_i:=A^i_i-\varrho_i'\nu_i-\sum_{s:s\neq i}\frac{\nu_s\varrho_s}{\psi_{si}}\ .\label{ris9}\eeq
Choose $A^i_i$ to set $\tilde{A}^i_i$ to 0; this is a diagonal transformation of the
coordinates. Then $\delta^{(2)}$ vanishes and the variation of $f^{ij}$ is reduced to
$\delta^{(1)}$.

On the other hand, under a variation of $\psi's$, $\psi_i\mapsto \psi_i+\delta\psi_i$, the 
Poisson brackets (\ref{rpb15}) transform in the following way:
\beq \delta_\psi f^{ij}=\frac{(x^i-x^j)^2}{\psi_{ij}^2}(\varrho_i\delta\psi_j-\varrho_j\delta\psi_i)
+a\left( (x^j)^2\delta\psi_i-(x^i)^2\delta\psi_j\phantom{\frac{}{}}\!\!\right)\label{ris10}\eeq
and we conclude that with the choice
\beq\delta\psi_i=\epsilon \varrho_i\nu_i\label{ris11}\eeq
the variation $\delta^{(1)}$ is compensated by the variation $\delta_\psi$. The coefficients 
of the polynomial $\varrho$ stay the same. In the next subsection we will study relations between
the variation of $\psi$'s and the polynomial $\varrho$.

\vskip .2cm
\noindent {\bf Remark.} With $\xi_i$ as in (\ref{ris6}), define three operators,
\beqa &(B^-)^i_j=\displaystyle{\frac{1}{\psi_{ij}}}\ ,\ i\neq j\ ,\ \ {\mathrm{and}}\ \ 
(B^-)^i_i=-\xi_i\ ,&\label{ops1}\\
&(B^0)^i_j=\displaystyle{\frac{\psi_i}{\psi_{ij}}}\ ,\ i\neq j\ ,\ \ {\mathrm{and}}\ \ (B^0)^i_i=-
(\frac{n-1}{2}+\psi_i\xi_i)\ ,&\label{ops2}\\
&(B^+)^i_j=\displaystyle{\frac{\psi_i^2}{\psi_{ij}}}\ ,\ i\neq j\ ,\ \ {\mathrm{and}}\ \ (B^+)^i_i=-
((n-1)\psi_i+\psi_i^2\xi_i)\ .&\label{ops3}\eeqa
The operators $B^+,B^0$ and $B^-$ generate an action of the Lie algebra $sl(2)$,
\beq [ B^0,B^-] =-B^-\ ,\ [ B^0,B^+] =B^+\ ,\ [ B^+,B^-] =-2B^0\ \label{ops4}\eeq
(to obtain the usual commutation relations for the generators of $sl(2)$, change the sign
of $B^+$).
 
\vskip .2cm
This is the usual projective action of $sl(2)$ on polynomials $f(t)$ of degree less than $n$, 
\beq B^-\, :\, f(t)\,\mapsto\, f'(t)\ ,\ \ B^0\, :\, f(t)\,\mapsto\, tf'(t)-\frac{n-1}{2}\, 
f(t)\ ,\ \ B^+\, :\, f(t)\,\mapsto\, t^2f'(t)-(n-1)t\, f(t)\ ,\eeq
written in the basis of the non-normalized Lagrange polynomials, $l_i(t)=
\displaystyle{\prod_{s:s\neq i}}(t-\psi_s)$, at points $\{ \psi_i\}$, $i=1,2,\dots ,n$.
Indeed, in the basis $\{ l_i(t)\}$, a polynomial $f(t)$, deg$(f)\leq n-1$, takes the form
$f=\sum f^il_i$, where $f^i=l_i(\psi_i)^{-1}f(\psi_i)$. We have 
\beq l_i'(t)=\displaystyle{\sum_{a:a\neq i}\prod_{b:b\neq a,i}}(t-\psi_b)\ ,\ \ \ {\mathrm{so}}\ \ \  l_i'(\psi_k)=\displaystyle{\prod_{b:b\neq k,i}\psi_{kb}=l_k(\psi_k)\frac{1}{\psi_{ki}}}\ ,
\ k\neq i\ .\eeq 
Also,
\beq l_i'(t)=l_i(t)\displaystyle{\sum_{s:s\neq i}\frac{1}{t-\psi_s}}\ ,\ \ \ {\mathrm{so}}\ \ \  
l_i'(\psi_i)=-l_i(\psi_i)\xi_i\ .\eeq 
Therefore, $l_i'(t)=\displaystyle{-\xi_il_i(t)+\sum_{k:k\neq i}\frac{1}{\psi_{ki}}l_k(t)}$, which 
is exactly (\ref{ops1}). For functions on the set of points $\{ \psi_i\}$, the operator of
multiplication by $t$ acts as a diagonal matrix Diag$(\psi_1,\psi_2,\dots ,\psi_n)$ and
(\ref{ops2})-(\ref{ops3}) follow.

\vskip .2cm
Define an involution $\varpi$ on the space of matrices\footnote{The involution $\varpi$ 
is the difference of two complementary projectors. The involution $\varpi$ satisfies
\beq\varpi (Y_1)\varpi (Y_2)+\varpi (Y_1Y_2)=\left\{\begin{array}{lll}
\varpi (\varpi (Y_1)\varpi (Y_2))+Y_1Y_2\ ,\\
\varpi (\varpi (Y_1)Y_2)+Y_1\varpi (Y_2)\ ,\\
\varpi (Y_1\varpi (Y_2))+\varpi (Y_1)Y_2\ 
\end{array}\right.\label{trid}\eeq
for arbitrary $Y_1,Y_2\in {\mathrm{Mat_n}}$. All other linear dependencies between
$Y_1Y_2$, $\varpi (Y_1)Y_2$, $Y_1\varpi (Y_2)$, $ \varpi (Y_1)\varpi (Y_2)$,
$\varpi (Y_1Y_2)$, $\varpi (\varpi (Y_1)Y_2)$, $\varpi (Y_1\varpi (Y_2))$ and
$\varpi (\varpi (Y_1)\varpi (Y_2))$
are consequences of the three identities (\ref{trid}).},
\beq \varpi (Y)^i_j= Y^i_j\ ,\ i\neq j\ \ \ {\mathrm{and}}\ \ \ \varpi (Y)^i_i=-Y^i_i\ \ ,\ \ 
Y\in {\mathrm{Mat_n}}\ .\label{ops5}\eeq

Let 
\beq B(\varrho )=aB^++bB^0+cB^-\ \ ,\ \ B(\varrho )\, :f\,\mapsto\, \varrho (t)f'(t)-\frac{n-1}{2}\varrho'(t)f(t)\ 
.\label{ops6}\eeq 
In the basis $\{ l_i(t)\}$ for $B$, the generator (\ref{ris6}) of the invariance group is
\beq {\cal{A}}(\varrho )=\varpi (B(\varrho ))\ .\label{ops7}\eeq
Note that the operators $\varpi (B^-)$, $\varpi (B^0)$ and $\varpi (B^+)$ do not form a Lie
algebra.

\subsection{Normal form}

In this subsection we derive a normal form of each individual member of the Poisson pencil
from the proposition \ref{rimepb}. It depends only on the discriminant of the polynomial $\varrho$.
When the discriminant of $\varrho$ is different from zero, the Poisson brackets (\ref{rpb15}) are 
equivalent to the Poisson brackets defined by the $r$-matrix (\ref{clcr}). When the polynomial 
$\varrho$ is different from zero but its discriminant is zero, the Poisson brackets (\ref{rpb15}) 
are equivalent to the Poisson brackets defined by the $r$-matrix (\ref{bee}).

\vskip .2cm
Under a variation of the polynomial $\varrho$, $\varrho (t)\mapsto (a+\delta a)t^2+(b+\delta b)t+(c+\delta c)$, we have for the Poisson brackets (\ref{rpb17}):
\beq \delta_\varrho f^{ij}=\displaystyle{\frac{u_{ij}^2\delta a+u_{ij}v_{ij}\delta b+
v_{ij}^2\delta c}{\psi_i-\psi_j}}\ .\label{ich}\eeq

The variation of $\varrho$ can be compensated by a variation (\ref{ris10}) of $\psi$'s if the 
coefficients in $(x^j)^2$, $x^ix^j$ and $(x^i)^2$ in the combination 
$(\delta_\psi +\delta_\varrho )f^{ij}$ vanish, which gives the following system:
\begin{eqnarray} (x^j)^2&\Rightarrow&\displaystyle{\frac{\varrho_i\delta\psi_j-
\varrho_j\delta\psi_i}{\psi_{ij}^2}}+a\delta
\psi_i+\displaystyle{\frac{\psi_i^2\delta a+\psi_i\delta b+\delta c}{\psi_{ij}}}=0\ ,\label{ich2}\\
x^ix^j&\Rightarrow&-\displaystyle{\frac{2(\varrho_i\delta\psi_j-
\varrho_j\delta\psi_i)}{\psi_{ij}^2}}-
\displaystyle{\frac{2\psi_i\psi_j\delta a+(\psi_i+\psi_j)\delta b+2\delta c}{\psi_{ij}}}=0
\ .\label{ich3}\end{eqnarray}
A combination 2$\times$(\ref{ich2})+(\ref{ich3}) gives
\beq 2a\delta\psi_i+2\psi_i\delta a+\delta b=0\ .\label{ich4}\eeq
Substituting the expression (\ref{ich4}) for $\delta\psi$'s into (\ref{ich2}) gives
\beq \delta D(\varrho )=0\ ,\ \ {\mathrm{where}}\ \ D(\varrho )=b^2-4ac\ .\label{ich5}\eeq
The coefficient in $(x^i)^2$ in $(\delta_\psi +\delta_\varrho )f^{ij}$ vanishes due to the anti-symmetry.

\vskip .2cm
Therefore, a necessary condition for a variation of $\varrho$ to be compensated by a variation
of $\psi$'s is that the discriminant $D(\varrho )$ does not vary. We shall now see that the
discriminant is the unique invariant.

Explicitly, under a shift $\psi_j\mapsto\psi_j+\zeta$, we have $u_{ij}\mapsto u_{ij}+\zeta v_{ij}$ 
and $v_{ij}\mapsto v_{ij}$ (in the notation (\ref{rpb17})), which produces the following effect 
on the coefficients of the polynomial $\varrho$:
\beq a\mapsto a\ ,\ b\mapsto b+2\zeta a\ ,\ c\mapsto c+\zeta b+\zeta^2 a\ .\label{ich6}\eeq
A dilatation $\psi_j\mapsto\lambda\psi_j$ produces the following effect on the 
coefficients of $\varrho$:
\beq a\mapsto \lambda a\ ,\ b\mapsto b\ ,\ c\mapsto \lambda^{-1}c\ .\label{ich7}\eeq
The inversion $\psi_j\mapsto\psi_j^{-1}$ accompanied by a change of variables 
$\tilde{x}^i=\psi_i^{-1}x^i$ produces the following effect on the coefficients of $\varrho$:
\beq a\mapsto -c\ ,\ b\mapsto -b\ ,\ c\mapsto -a\ .\label{ich8}\eeq
The set of operators (\ref{ich6}) and (\ref{ich7}) generates the action of the affine group
on the space of the polynomials $\varrho$. The affine group, together with the inversion 
(\ref{ich8}) generates an action\footnote{Let $e_+$ be the generator of the 1-parametric group 
(\ref{ich6}) and $h$ the generator of the 1-parametric group (\ref{ich7}). Denote by ${\cal{I}}$ 
the inversion (\ref{ich8}). The remaining generator $e_-$ is
${\cal{I}}e_+{\cal{I}}$.} of $so(3)$ (the spin 1 representation of $sl(2)$) on the space of the 
polynomials $\varrho$ and the classification reduces to that of orbits. The orbits (in the complex 
situation) of non-zero polynomials are of two types: "massive", $D(\varrho )\neq 0$, or 
"light-like", $D(\varrho )=0$. Particular representatives of both types appear in the Poisson 
brackets, corresponding to the rime $r$-matrices (see the remark 1 after the proposition 
\ref{rimepb}) and thus to the $r$-matrices studied in subsections \ref{nskc} and \ref{sksc}. We 
obtain the following statement.

\begin{proposition}$\!\!\!${\bf .} Let $\varrho (t)$ be a non-zero quadratic polynomial.

\vskip .1cm
If the discriminant of $\varrho$ is different from zero, $D(\varrho)\neq 0$, then there exists a 
change of the parameters $\psi^i$ in the Poisson brackets (\ref{rpb15}) which sets $\varrho (t)$ 
to $bt$, $\varrho (t)\mapsto bt$; these are the Poisson brackets corresponding to the $r$-matrix 
$r_{CG}$ (subsection \ref{nskc}). 
 
\vskip .1cm
If the discriminant of $\varrho$ is zero, $D(\varrho )=0$, then there exists a change of the 
parameters $\psi^i$ in the Poisson brackets (\ref{rpb15}) which sets $\varrho (t)$ to $c$, 
$\varrho (t)\mapsto c$; these are the Poisson brackets corresponding to the $r$-matrix $b_{CG}$ 
(subsection \ref{sksc}).\end{proposition}

The generator ${\cal{A}}(\varrho )$ of the invariance group can be easily described in both cases,
$D(\varrho )\neq 0$ and $D(\varrho )=0$, in the parameter-free basis (that is, for the 
$r$-matrices $r_{CG}$ and $b_{CG}$; in the rime basis the generators are given by (\ref{inr2}) and 
(\ref{inr4}), respectively). For $D(\varrho )\neq 0$ (respectively, $D(\varrho )=0$), it coincides 
with the matrix of the operator $B_0$ (respectively, $B_-$), as in the remark in subsection 
\ref{sinv}, in the basis $\{ t^i\}$ of powers of the variable $t$. This implies somewhat
unexpectedly that for an arbitrary polynomial $\varrho (t)$ the matrices ${\cal{A}}(\varrho )$ and 
$\varpi ({\cal{A}}(\varrho ))$ are related by a similarity transformation. Note that in the basis
$\{ t^i\}$ of powers, the operators $aB^++bB^0+cB^-$ and $\varpi (aB^++bB^0+cB^-)$ are also
related by a similarity transformation for arbitrary $a,b$ and $c$ but here it is obvious:
$\varpi (aB^++bB^0+cB^-)=aB^+-bB^0+cB^-$, so the operator $\varpi (aB^++bB^0+cB^-)$ belongs
to $sl(2)$ and moreover lies on the same (complex) orbit as $aB^++bB^0+cB^-$ with respect to the
adjoint action.

\section{Orderable quadratic rime associative algebras\vspace{.1cm}}

Consider an associative algebra ${\cal{A}}$ defined by quadratic relations giving a 
lexicographical order. This means that $x^jx^k$ for $j<k$ is a linear combination of terms 
$x^ax^b$ with $a\geq b$ and either $a>j$ or $a=j$ and $b>k$.

\vskip .2cm
We shall say that such algebra ${\cal{A}}$ is {\it rime} if $\{ a,b\}\subset\{ j,k\}$. In other 
words, the relations in the algebra are
\beq x^jx^k=f_{jk}x^kx^j+g_{jk}x^kx^k\ ,\ j<k\ .\label{qra1}\eeq

We shall classify the {\it strict} rime algebras ${\cal{A}}$ (that is, the algebras for which 
all coefficients $f_{ij}$ and $g_{ij}$ are different from zero for $i<j$). 

\vskip .2cm
The only possible overlaps for the set of relations (\ref{qra1}) are of the form 
$(x^jx^k)x^l=x^j(x^kx^l)$, $j<k<l$. The ordered form of the expression $(x^jx^k)x^l$ is
\beq \begin{array}{lcl}(x^jx^k)x^l&=&f_{jk}f_{jl}f_{kl}\ x^lx^kx^j+
f_{jk}f_{jl}g_{kl}\ x^lx^lx^j+f_{kl}^2g_{jk}\ x^lx^kx^k\\[.5em]
&+&(f_{kl}g_{jk}g_{kl}+f_{kl}^2(f_{jk}g_{jl}+g_{jk}g_{kl}))\, x^lx^lx^k\\[.5em]
&+&(f_{jk}g_{jl}g_{kl}+g_{jk}g_{kl}^2+f_{kl}g_{kl}(f_{jk}g_{jl}+g_{jk}g_{kl}))\, x^lx^lx^l
\ .\end{array}\label{qra2}\eeq
The ordered form of the expression $x^j(x^kx^l)$ is
\beq \begin{array}{lcl}x^j(x^kx^l)&=&f_{jk}f_{jl}f_{kl}\ x^lx^kx^j+
f_{jl}^2g_{kl}\ x^lx^lx^j+f_{kl}f_{jl}g_{jk}\ x^lx^kx^k\\[.5em]
&+&f_{kl}g_{jl}\ x^lx^lx^k+
(g_{kl}g_{jl}+f_{jl}g_{kl}g_{jl})\, x^lx^lx^l\ .\end{array}\label{qra3}\eeq
Equating coefficients, we find
\beqa &x^lx^lx^j\, :\ f_{jk}f_{jl}g_{kl}=f_{jl}^2g_{kl}\ ,&\label{qra4}\\
&x^lx^kx^k\, :\ f_{kl}^2g_{jk}=f_{kl}f_{jl}g_{jk}\ ,&\label{qra5}\\
&x^lx^lx^k\, :\ f_{kl}g_{jk}g_{kl}+f_{kl}^2(f_{jk}g_{jl}+g_{jk}g_{kl})=f_{kl}g_{jl}\ ,&
\label{qra6}\\
&x^lx^lx^l\, :\ f_{jk}g_{jl}g_{kl}+g_{jk}g_{kl}^2+f_{kl}g_{kl}(f_{jk}g_{jl}+g_{jk}g_{kl})=g_{kl}g_{jl}+
f_{jl}g_{kl}g_{jl}\ .&\label{qra7}\eeqa
In the strict situation, eqs. (\ref{qra4}) and (\ref{qra5}) simplify, respectively, to
\beqa &f_{jk}=f_{jl}\ ,\ \ {\mathrm{for}}\ j<k\ \ {\mathrm{and}}\ j<l\ ,&\label{qra8}\\
&f_{kl}=f_{jl}\ ,\ \ {\mathrm{for}}\ j<l\ \ {\mathrm{and}}\ k<l\ .\ &\label{qra9}\eeqa
Eqs. (\ref{qra8}) and (\ref{qra9}) imply that $f_{jk}$'s are all equal,
\beq f_{jk}=:f\ .\label{qra10}\eeq
The substitution of (\ref{qra10}) into (\ref{qra6}) gives (in the strict situation)
\beq (f+1)\Bigl( g_{jk}g_{kl}+g_{jl}(f-1)\Bigr) =0\ \ {\mathrm{for}}\ j<k<l\ .\label{qra11}\eeq
Eq. (\ref{qra7}) follows from (\ref{qra10}) and (\ref{qra11}).

\vskip .2cm
We have thus two cases: 

\vskip .1cm\noindent
(i) $f=-1$ and no extra conditions on $g_{jk}$'s; 

\vskip .1cm\noindent
(ii) $f\neq -1$ and 
\beq g_{jk}g_{kl}=(1-f)\, g_{jl}\ \ {\mathrm{for}}\ j<k<l\ ;\label{qra12}\eeq
$1-f\neq 0$ since $g_{jk}\neq 0$ and $g_{kl}\neq 0$.

\vskip .2cm
In the case (ii), make an appropriate rescaling of generators, $x^i\mapsto d_ix^i$ to achieve
\beq g_{i,i+1}=1-f\ \ {\mathrm{for}}\ \ {\mathrm{all}}\ i=1,\dots ,n-1\ .\label{qra13}\eeq
It then follows from eq. (\ref{qra12}) that
\beq g_{ij}=1-f\ \ {\mathrm{for}}\ \ {\mathrm{all}}\ i<j\ .\label{qra14}\eeq

We summarize the obtained results.

\begin{proposition}$\!\!\!${\bf .} Up to a rescaling of variables, the general orderable 
quadratic strict rime algebra has relations 

\vskip .15cm
\noindent
(i) either of the form
\beq x^jx^k=-x^kx^j+g_{jk}x^kx^k\ ,\ j<k\ ,\label{qra15}\eeq
with no conditions on the coefficients $g_{jk}$; 

\vskip .15cm
\noindent
(ii) or of the form
\beq x^jx^k=fx^kx^j+(1-f)x^kx^k\ ,\ j<k\ ,\ \label{qra16}\eeq
with arbitrary $f$ (it is strict when $f\neq 0,1$).
\end{proposition}

By construction, the algebras of types (i) and (ii) possess a basis formed by ordered 
monomials and thus have the Poincar\'e series of the algebra of commuting variables.

\vskip .2cm
The algebra with defining relations (\ref{qra16}) is the quantum space for the $R$-matrix 
(\ref{rstcl}). The relations (\ref{qra16}) can be written in the form 
\beq (x^j-x^k)x^k=fx^k(x^j-x^k)\ ,\ j<k\ ;\label{qra17}\eeq
this is a quantization of the Poisson brackets 
\beq \{ x^j,x^k\} =x^k(x^j-x^k)\ ,\ j<k\ .\label{qra18}\eeq

\vskip .2cm
It would be interesting to know if the algebra with the defining 
relations (\ref{qra15}) admits an $R$-matrix description.

\section*{Acknowledgements}
It is our pleasure to thank L\'aszl\'o Feh\'er, Alexei Isaev and Milen Yakimov for enlightening 
discussions. The work was partially supported by the ANR project GIMP No.ANR-05-BLAN-0029-01. The 
second author (T. Popov) was also partially  supported by the Program ``Bourses d'\'echanges 
scientifiques pour les pays de l'Est europ\'een'' and by the Bulgarian National Council for 
Scientific Research project PH-1406.

\section*{Appendix A.$\ $ Equations}\addcontentsline{toc}{section}{Appendix A.$\ $ Equations}

Here we give the list of the equations $Y\!\! B(\hat{R})^{ijk}_{abc}=0$ for the rime matrix
\beq \hat{R}_{kl}^{ij}=\alpha_{ij}\delta^i_l\delta^j_k +\beta_{ij}\delta^i_k\delta^j_l 
+\gamma_{ij}\delta^i_k\delta^i_l+\gamma'_{ij}\delta^j_k\delta^j_l\ ,\eeq
with a convention $\alpha_i=\alpha_{ii}$ and $\beta_{ii}=\gamma_{ii}=\gamma_{ii}'=0$.

The rime Ansatz implies that $Y\!\! B(\hat{R})^{ijk}_{abc}$ can be different from zero only if
the set of lower indices is contained in the set of upper indices. Therefore, the equations split
into two lists: the first one with two different indices among $\{ i,j,k\}$ and the second one
with three different indices.

The full set of equations $Y\!\! B(\hat{R})^{ijk}_{abc}=0$ is invariant under the 
involution $\iota$,
\beq\iota\, :\ \alpha_i\leftrightarrow\alpha_i\ ,\  \alpha_{ij}\leftrightarrow \alpha_{ji}\ ,\ 
\beta_{ij}\leftrightarrow \beta_{ji}\ ,\ \gamma_{ij}\leftrightarrow\gamma'_{ji}\ ,\label{invo}\eeq
for if $\hat{R}$ is a solution of the YBe then $\hat{R}_{21}=P\hat{R}P$ is a solution of the
YBe as well. We shall write only the necessary part of the equations, the rest can be obtained
by the involution $\iota$.

\vskip .2cm
The equations $Y\!\! B(\hat{R})^{ijk}_{abc}=0$ with two different indices are:
\beq\label{yb1}\alpha_{ij}\gamma_{ij}(\gamma_{ji}+\gamma'_{ij})=0\ ,\eeq
\vspace{-.6cm}
\beq\label{yb2}\alpha_{ij}(\beta_{ij}\beta_{ji}+\gamma_{ij}\gamma'_{ij})=0=
\alpha_{ij}(\beta_{ij}\beta_{ji}-\gamma_{ij}\gamma_{ji})\ ,\eeq
\beq\label{yb3}\alpha_{ij}\gamma_{ij}(\alpha_{ij}+\beta_{ji}-\alpha_j)=0=
\alpha_{ij}\gamma_{ij}(\alpha_{ji}+\beta_{ij}-\alpha_j)\ ,\eeq
\beq\label{yb4}\beta_{ij}(\alpha_i^2-\alpha_{ij}\alpha_{ji}-\alpha_i\beta_{ij})
+(\alpha_i-\beta_{ij})\gamma_{ij}\gamma'_{ij} =0\ ,\eeq
\beq(\alpha_i-\alpha_j)\gamma^2_{ij}+\alpha_{ij}\gamma_{ij}(\gamma_{ij}+\gamma'_{ji})=0
\ ,\label{eI}\eeq
\beq\label{eI1}\alpha_{ij}\beta_{ij}\gamma'_{ji}+(\alpha_i\beta_{ij}+
\gamma'_{ij}\gamma_{ij})\gamma_{ij}=0\ ,\eeq
\beq\label{eI2}(\alpha_{ij}-\alpha_{ji}-\beta_{ij}+\beta_{ji})\gamma_{ij} \gamma'_{ji}=0
=(\alpha_{ij}-\alpha_{ji}-\beta_{ij}+\beta_{ji})\beta_{ij} \beta_{ji}\ ,\eeq
\beq\label{eI3}\alpha_{ij}\gamma'_{ji}(\alpha_j-\alpha_{ij})+
\beta_{ji}\gamma_{ij}(\alpha_i-\beta_{ji})+
\gamma_{ij}(\beta_{ij}\beta_{ji}+\gamma_{ji}\gamma'_{ji})=0\ ,\eeq
\beq\label{eI4}(\alpha^2_i-\alpha_i(\alpha_{ji}+\beta_{ji}) +\beta_{ij}\beta_{ji} 
-\gamma_{ij}\gamma_{ji})\gamma_{ij}  
=(\alpha^2_i - \alpha_i(\alpha_{ij}+\beta_{ij}) +\beta_{ij}\beta_{ji} 
-\gamma'_{ij}\gamma'_{ji})\gamma'_{ji}\ .\eeq

\vskip .2cm
The equations with three different indices $\{i,j,k \}$ are:
\beq (\alpha_{ij}-\alpha_{ki}-\beta_{ij}+\beta_{ki})\gamma_{ij}\gamma'_{ki}=0\ ,\label{ee1}\eeq
\vspace{-.7cm}
\beq\alpha_{ij}(\beta_{ij}\beta_{jk}+\beta_{ik}\beta_{ji}-\beta_{ik}\beta_{jk})=0\ ,\label{ee2}\eeq
\beq\alpha_{ij}(\gamma_{ij}\gamma_{jk}+\gamma_{ik}(\beta_{jk}-\beta_{ji}))=
\alpha_{ij}(\gamma_{ij}\gamma'_{kj}+\gamma_{ik}(\beta_{kj}-\beta_{ij}))=0\ ,\label{ee3}\eeq
\beq (\alpha_{ij}\alpha_{ji}-\alpha_{jk}\alpha_{kj})\beta_{ik}
+\beta_{ij}\beta_{jk}(\beta_{ij}-\beta_{jk})=0\ ,\eeq
\beq (\alpha_i+\beta_{ik}-\beta_{ji})\beta_{ji}\gamma_{ik}+\gamma_{ik}\gamma_{ji}\gamma'_{ji}
+\alpha_{ik}(\gamma_{jk}\gamma'_{ji}+\beta_{jk}\gamma'_{ki})=0\ ,\eeq
\beq (\alpha_{i}+\alpha_{ij}-\alpha_{kj}-\beta_{kj})\gamma_{ij}\gamma_{ik}
-\gamma_{ik}^2\gamma_{kj}+ \gamma_{ij}(\alpha_{ik}\gamma'_{ki}-\gamma_{ij}\gamma'_{kj})=0\ ,\eeq
\beq (\alpha_{i}-\beta_{kj})\beta_{ij}\gamma_{ik}+(\beta_{ik}\beta_{kj}+\gamma_{ij}\gamma'_{ij})
\gamma_{ik}+ \alpha_{ik}\beta_{ij}\gamma'_{ki}-(\beta_{ij}-\beta_{ik})\gamma_{ij}\gamma'_{kj}=0
\ ,\eeq
\beq\alpha_{ij}(\gamma_{ij}\gamma_{jk}+\gamma_{ik}(\alpha_{jk}-\alpha_{ji}))\! =\!
\alpha_{ji}(\gamma_{ij}\gamma_{jk}+\gamma_{ik}(\alpha_{jk}-\alpha_{ji}))\! =\! 
\alpha_{ij}(\gamma_{ij}\gamma'_{kj}+\gamma_{ik}(\alpha_{kj}-\alpha_{ij}))\! =\! 0.\eeq
 
\section*{Appendix B.$\ $ Blocks}
\addcontentsline{toc}{section}{Appendix B.$\ $ Blocks\vspace{.25cm}}

We analyze here the structure of 4$\times$4 blocks of an invertible and skew-invertible
rime $R$-matrix corresponding to two-dimensional coordinate planes. 

We denote the matrix elements as in (\ref{fri}).

\vskip .2cm
The skew-invertibility of a rime $R$-matrix imposes restrictions on its entries: in the line 
$\hat{R}^{i*}_{j*}$ only two entries can be non-zero, $\hat{R}^{ij}_{ji}$ and $\hat{R}^{ij}_{jj}$; 
in the line $\hat{R}^{*j}_{*i}$ only two entries can be non-zero, $\hat{R}^{ij}_{ji}$ and 
$\hat{R}^{ij}_{ii}$. Therefore,
\beq \alpha_{ij}=0\ \ \Rightarrow\ \ \gamma_{ij}\gamma_{ij}'\neq 0\ \ \ \ {\mathrm{and}}\ \ \ \ 
\gamma_{ij}\gamma_{ij}'=0 \ \ \Rightarrow\ \ \alpha_{ij}\neq 0\ . \label{skinv}\eeq

Dealing with a single block, this becomes especially clear: to skew invert a 4$\times$4
block is the same as to invert the matrix
\beq\label{skb}
\left(\begin{array}{cccc}\alpha_1&0&\gamma_{12}&\beta_{12}\\0&0&\alpha_{12}&\gamma_{12}'\\
\gamma_{21}'&\alpha_{21}&0&0\\ \beta_{21}&\gamma_{21}&0&\alpha_2\end{array}\right)\ ,\eeq
whose determinant is
\beq (\alpha_{12}\beta_{12}-\gamma_{12}\gamma_{12}')(\alpha_{21}\beta_{21}-\gamma_{21}\gamma_{21}')
-\alpha_1\alpha_2\alpha_{12}\alpha_{21}\ .\eeq

\subsection*{B.1\hspace{.3cm} Solutions}
\addcontentsline{toc}{subsection}{B.1\hspace{.2cm} Solutions}

Here we classify solutions which are neither ice nor strict rime. For an ice $R$-matrix,
$\alpha_{12}\neq 0$ and $\alpha_{21}\neq 0$. For a rime $R$-matrix, $\alpha_{ij}$ might vanish
and we consider the cases according to the number of $\alpha_{ij}$'s which can be zero.

\paragraph{1.} Both $\alpha_{12}$ and $\alpha_{21}$ do not vanish, $\alpha_{12}\alpha_{21}\neq 0$. 

\vskip .2cm
If $\gamma_{12}\gamma_{21}\neq 0$ then by (\ref{yb1}), $\gamma_{12}'\gamma_{21}'\neq 0$.
This is strict rime.

\vskip .2cm
If both $\gamma_{12}=0$ and $\gamma_{21}=0$ then eq. (\ref{yb3}) implies $(\alpha_{ji}+\beta_{ij}
-\alpha_j)\gamma_{ji}'=0$; eq. (\ref{eI}) implies $(\alpha_i-\alpha_j+\alpha_{ji})\gamma_{ji}'=0$
and eq. (\ref{eI1}) implies $\beta_{ij}\gamma_{ji}'=0$. Combining these, we find $\gamma_{ij}'=0$,
this is ice.

\vskip .2cm
It is left to analyze the situation when only one of $\gamma$'s is different from zero, say  
$\gamma_{12}\neq 0$ and $\gamma_{21}=0$. We have the following chain of implications:
\beq\begin{array}{ccc}\label{ca1}(\ref{yb1})&\Rightarrow&\ \ \ \ \gamma_{12}'=0\ ,\end{array}\eeq
\beq\begin{array}{ccc}\label{ca3}(\ref{yb3})&\Rightarrow&\left\{\begin{array}{l}\beta_{12}=\alpha_2
-\alpha_{21}\ ,\\ 
\beta_{21}=\alpha_2-\alpha_{12}\ ,\end{array}\right.\end{array}\eeq
\beq\begin{array}{ccc}\label{ca2}(\ref{yb2})&\Rightarrow&\ \ \ \ (\alpha_2-\alpha_{12})(\alpha_2-\alpha_{21})=0\ ,\end{array}\eeq
\beq\begin{array}{ccc}\label{ca4}(\ref{yb4})&\Rightarrow&\left\{\begin{array}{l}(\alpha_1-\alpha_2)
(\alpha_2-\alpha_{21})(\alpha_1+\alpha_{21})=0\ ,\\
(\alpha_1-\alpha_2)(\alpha_2-\alpha_{12})(\alpha_1+\alpha_{12})=0\ , \end{array}\right.\end{array}\eeq
\beq\begin{array}{ccc}\label{ca5}(\ref{eI})\& (\ref{eI4})&\Rightarrow&\left\{\begin{array}{l}
(\alpha_1-\alpha_2+\alpha_{12})\gamma_{12}+\alpha_{12}\gamma_{21}'=0\ ,\\
(\alpha_1-\alpha_2+\alpha_{21})\gamma_{21}'+\alpha_{21}\gamma_{12}=0\ .
\end{array}\right.\end{array}\eeq
Eqs. (\ref{eI1}), (\ref{eI2}) and (\ref{eI3}) are satisfied. By the second line in (\ref{ca5}), 
$\gamma_{21}'\neq 0$.

\vskip .2cm
Now the system of inequalities and equations is invariant under $\hat{R}\leftrightarrow\hat{R}_{21}$, so up to this transformation we can solve eq. (\ref{ca2})
by setting $\alpha_{21}=\alpha_2$. Then, by (\ref{ca5}), $\gamma_{21}'=-
\gamma_{12}\alpha_2/\alpha_1$, $\beta$'s are expressed in terms of $\alpha$'s by (\ref{ca3}) 
and the remaining system for $\alpha$'s reduces to a single equation
$(\alpha_1-\alpha_2)(\alpha_1+\alpha_{12})=0$. We obtain two solutions:

\vskip .2cm
{\bf 1a.} $\alpha_2=\alpha_1$; $\alpha_1$, $\alpha_{12}$ and $\gamma_{12}$ are arbitrary non-zero 
numbers; we rescale the $R$-matrix to set $\alpha_1\alpha_{12}=1$ and denote $q=\alpha_1$,
$\gamma =\gamma_{12}$:
\beq\label{rbl1}\hat{R}_{(q;\gamma)}=
\left(\begin{array}{cccc}\alpha_1&0&0&0\\ \gamma_{12}&0&\alpha_{12}&0\\
-\gamma_{12}&\alpha_1&\alpha_1-\alpha_{12}&0\\ 0&0&0&\alpha_1\end{array}\right)
=\left(\begin{array}{cccc}q&0&0&0\\ \gamma&0&q^{-1}&0\\
-\gamma&q&q-q^{-1}&0\\ 0&0&0&q\end{array}\right)\ .\eeq
The $R$-matrix (\ref{rbl1}) is semi-simple if (and only if) $q+q^{-1}\neq 0$ and it is 
then an $R$-matrix of $GL(2)$-type, ${\mathrm{Spec}}(\hat{R})=\{ q,q,q,-q^{-1}\}$. This solution 
is a specialization of (\ref{subst})-(\ref{subst2}).

\vskip .2cm
{\bf 1b.} $\alpha_{12}=-\alpha_1$; $\alpha_1$, $\alpha_2$ and $\gamma_{12}$ are arbitrary non-zero 
numbers; we rescale the $R$-matrix to set $\alpha_1\alpha_2=-1$ and denote $q=\alpha_1$,
$\gamma =\gamma_{12}/q$:
\beq\label{rbl2}\hat{R}_{(q;\gamma)}=
\left(\begin{array}{cccc}\alpha_1&0&0&0\\ 
\gamma_{12}&0&-\alpha_1&0\\ -\gamma_{12}\alpha_2/\alpha_1&\alpha_2&\alpha_1+\alpha_2&0\\ 
0&0&0&\alpha_2\end{array}\right) =\left(\begin{array}{cccc}q&0&0&0\\ 
q\gamma&0&-q&0\\ q^{-1}\gamma&-q^{-1}&q-q^{-1}&0\\ 0&0&0&-q^{-1}\end{array}\right)\ .\eeq
The $R$-matrix (\ref{rbl2}) is semi-simple if (and only if) $q+q^{-1}\neq 0$ and it is then an 
$R$-matrix of $GL(1|1)$-type, ${\mathrm{Spec}}(\hat{R})=\{ q,q,-q^{-1},-q^{-1}\}$.

\paragraph{2.} Assume that $\alpha_{12}=0$. 

\vskip .2cm
By the invertibility, $\beta_{12}\beta_{21}\neq 0$; by 
the skew-invertibility, $\gamma_{12}\gamma_{12}'\neq 0$; now
eqs. (\ref{yb1}) and (\ref{yb2}) imply $\beta_{12}\beta_{21}=\gamma_{12}\gamma_{21}$, 
$\gamma_{12}'=-\gamma_{21}$ and $\gamma_{21}'=-\gamma_{12}$. Eq. (\ref{yb3}) implies
$\alpha_2=\alpha_1$, $\beta_{12}=\alpha_1-\alpha_{21}$ and $\beta_{21}=\alpha_1$. 

\vskip .2cm
The rest is satisfied and we obtain a solution, in which $\alpha_1$, $\beta_{12}$  and 
$\gamma_{12}$ are arbitrary non-zero numbers; we rescale the $R$-matrix to set 
$\alpha_1\beta_{12}=-1$ and denote $q=\alpha_1$, $\gamma =\gamma_{12}$:
\beq\label{rbl3}\hat{R}_{(q;\gamma)}=\left(\begin{array}{rccr}
\alpha_1&0&0&0\\ \gamma_{12}&\beta_{12}&0&-\alpha_1\beta_{12}/\gamma_{12}\\[.2em]
-\gamma_{12}&\alpha_1-\beta_{12}&\alpha_1&\alpha_1\beta_{12}/\gamma_{12}\\ 
0&0&0&\alpha_1\end{array}\right) =
\left(\begin{array}{rccr}q&0&0&0\\ \gamma&-q^{-1}&0&1/\gamma\\[.2em]
-\gamma&q+q^{-1}&q&-1/\gamma\\ 0&0&0&q\end{array}\right)\ .\eeq
The $R$-matrix (\ref{rbl3}) is semi-simple if (and only if) $q+q^{-1}\neq 0$ and it is 
then an $R$-matrix of $GL(2)$-type, ${\mathrm{Spec}}(\hat{R})=\{ q,q,q,-q^{-1}\}$.
This solution is a specialization of (\ref{subst})-(\ref{subst2}).

\paragraph{3.} Finally, assume that $\alpha_{12}=\alpha_{21}=0$. 

\vskip .2cm
By the invertibility, $\beta_{12}\beta_{21}\neq 0$; by the
skew-invertibility, $\gamma_{12}\gamma_{12}'\gamma_{21}\gamma_{21}'\neq 0$; now eq. (\ref{eI})
implies $\alpha_2=\alpha_1$, eq. (\ref{eI2}) implies $\beta_{21}=\beta_{12}$; eq. (\ref{eI1}) 
implies $\gamma_{12}\gamma_{12}'=\gamma_{21}\gamma_{21}'=-\alpha_1\beta_{12}$; eq. (\ref{eI4})
implies that $\gamma_{12}\gamma_{21}$ can take three values: $\alpha_1^2$, $\beta_{12}^2$ or
$(-\alpha_1\beta_{12})$. 

\vskip .2cm
The rest is satisfied and we obtain a solution, in which $\alpha_1$, $\beta_{12}$  and $\gamma_{12}$ are arbitrary non-zero numbers;
we rescale the $R$-matrix to set $\alpha_1\beta_{12}=-1$ and denote $q=\alpha_1$,
$\gamma =\gamma_{12}$:
\beq\label{rbl4}\hat{R}_{(q,\omega ;\gamma)}=\left(\begin{array}{cccc}\alpha_1&0&0&0\\ 
\gamma_{12}&\beta_{12}&0&-\alpha_1\beta_{12}/\gamma_{12}\\
-\alpha_1\beta_{12}\gamma_{12}/\omega&0&\beta_{12}&\omega /\gamma_{12}\\ 
0&0&0&\alpha_1\end{array}\right) =\left(\begin{array}{cccc}q&0&0&0\\ \gamma&-q^{-1}&0&1/\gamma\\ 
\gamma/\omega&0&-q^{-1}&\omega /\gamma\\ 0&0&0&q\end{array}\right)\ ,\eeq
where $\omega =q^2,1,q^{-2}$. The $R$-matrix (\ref{rbl4}) is semi-simple if (and only if) 
$q+q^{-1}\neq 0$ and it is then an $R$-matrix of $GL(1|1)$-type, ${\mathrm{Spec}}(\hat{R})=\{ 
q,q,-q^{-1},-q^{-1}\}$.

\vskip .25cm
It follows from the analysis above that if $\gamma_{ij}\neq 0$ in an invertible and 
skew-invertible rime $R$-matrix then $\gamma_{ji}'\neq 0$.
 
\vskip .2cm
In each of the cases (\ref{rbl1})-(\ref{rbl4}), the parameter $\gamma\neq 0$ can be set to an 
arbitrary (non-zero) value by a diagonal change of basis. The $R$-matrices 
(\ref{rbl1})-(\ref{rbl4}) are skew-invertible.
 
\subsection*{B.2\hspace{.3cm} GL(2) and GL(1$|$1) $R$-matrices}
\addcontentsline{toc}{subsection}{B.2\hspace{.2cm} GL(2) and GL(1$|$1) $R$-matrices}

\paragraph{\bf 1.} In dimension 2, except the standard $R$-matrices of $GL$-type,
\beq \hat{R}^{GL(2)}_{(q,p)}=\left(\begin{array}{cccc}q&0&0&0\\0&0&p&0\\0&p^{-1}&q-q^{-1}&0\\ 
0&0&0&q\end{array}\right)\ \ ,\ \ \hat{R}^{GL(1|1)}_{(q,p)}=\left(\begin{array}{cccc}q&0&0&0\\0&0&p&0\\0&p^{-1}&q-q^{-1}&0\\ 
0&0&0&-q^{-1}\end{array}\right)\ ,\label{li2a}\eeq
there are two non-standard one-parametric families of non-unitary $R$-matrices of the type 
$GL(1|1)$: the eight-vertex one,
\beq \hat{R}^{I}_{(q)}=\frac{1}{2}\left(\begin{array}{cccc}q-q^{-1}+2&0&0&q-q^{-1}\\ 0&q-q^{-1}&q+q^{-1}&0\\
0&q+q^{-1}&q-q^{-1}&0\\ 
q-q^{-1}&0&0&q-q^{-1}-2\end{array}\right)\ ,\label{li2b}\eeq
and the matrix $\hat{R}^{(II)}$ for which the matrix $R=P\hat{R}$ can be given an upper-triangular form,
\beq \hat{R}^{II}_{(q,\varepsilon 
)}=\left(\begin{array}{cccc}q&0&0&q+q^{-1}\\0&0&\varepsilon q^{-1}&0\\0&\varepsilon q&q-q^{-1}&0\\ 
0&0&0&-q^{-1}\end{array}\right)\ ,\label{li2c}\eeq
where $\varepsilon =\pm 1$.

\vskip .2cm
The $R$-matrices (\ref{li2a}), (\ref{li2b}) and (\ref{li2c}) are semi-simple if (and only if) 
$q+q^{-1}\neq 0$.

\vskip .2cm
Up to the transformations $\hat{R}\leftrightarrow\hat{R}_{21}$ and 
$\hat{R}\leftrightarrow\hat{R}^t$ (the transposition), basis changes and rescalings 
$\hat{R}\mapsto c\,\hat{R}$ (where $c$ is a constant), the complete list of semi-simple invertible and skew-invertible $R$-matrices includes (see \cite{Hi} for a description of all solutions of the 
Yang--Baxter equation in two dimensions and \cite{EOW} for the classification of $GL(2)$-type
$R$-matrices), in addition to (\ref{li2a})-(\ref{li2c}), the one-parametric family of Jordanian 
solutions $\hat{R}^{(J)}_{(h_1:h_2)}$, 
\beq \hat{R}^{(J)}_{(h_1:h_2)}=\left(\begin{array}{cccc}1&h_1&-h_1&h_1h_2\\0&0&1&-h_2\\0&1&0&h_2\\ 
0&0&0&1\end{array}\right)\label{rjo}\eeq
(the Jordanian $R$-matrix is of $GL(2)$-type; it is unitary; the essential parameter is the
projective vector $(h_1:h_2)$), as well as the permutation-like 
solution $\hat{R}^{(P)}_{(a,b,c)}$ and one more solution $\hat{R}^{(')}_{(a)}$,
\beq \hat{R}^{(P)}_{(a,b,c)}=\left(\begin{array}{cccc}1&0&0&0\\0&0&a&0\\0&b&0&0\\ 
0&0&0&c\end{array}\right)\ \ ,\ \ 
\hat{R}^{(')}_{(a)}=\left(\begin{array}{cccc}0&0&0&a\\0&1&0&0\\0&0&1&0\\ 
a&0&0&0\end{array}\right)\ .\eeq
The $R$-matrix $\hat{R}^{(P)}_{(a,b,c)}$ is Hecke when $ab=1$ and $c=\pm 1$ and it is then 
standard (and unitary). The $R$-matrix $\hat{R}^{(')}_{(a)}$ is Hecke when $a^2=1$; it is then 
unitary and related to the standard $R$-matrix by a change of basis with 
the matrix $\left(\begin{array}{cc}1&1\\-1&1\end{array}\right)$ .
 
\vskip .2cm
Without the demand of semi-simplicity, the full list of invertible and skew-invertible 
$R$-matrices contains two more solutions,
\beq \hat{R}^{('')}_{(h_1:h_2:\sqrt{h_3})}=
\left(\begin{array}{cccc}1&h_1&h_2&h_3\\0&0&1&h_1\\0&1&0&h_2\\ 
0&0&0&1\end{array}\right)\ \ ,\ \ \hat{R}^{(''')}=\left(\begin{array}{cccc}1&0&0&1\\0&0&-1&0\\0&-1&0&0\\ 
0&0&0&1\end{array}\right)\ .\eeq
The essential parameter for the $R$-matrix $\hat{R}^{('')}_{(h_1:h_2:\sqrt{h_3})}$ is the 
projective vector $(h_1:h_2:\sqrt{h_3})$. The $R$-matrix $\hat{R}^{('')}_{(h_1:h_2:\sqrt{h_3})}$ 
is semi-simple if and only if $h_2=-h_1$ and $h_3=-h_1^2$; it then belongs to the family 
(\ref{rjo}) of Jordanian $R$-matrices.

\paragraph{\bf 2.} For the $R$-matrices from the list above, the transformations $\hat{R}\leftrightarrow\hat{R}_{21}$,
$\hat{R}\leftrightarrow\hat{R}^t$ and $\hat{R}\leftrightarrow\hat{R}^{-1}$ partly overlap or 
reduce to parameter or basis changes. We shall write formulas for the Hecke $R$-matrices only. 

\vskip .1cm
For the standard $R$-matrix $\hat{R}_{(q,p)}:=\hat{R}^{GL(2)}_{(q,p)}$,
\beq \hat{R}_{(q,p)}^t=\hat{R}_{(q,p^{-1})}\ ,\
(\hat{R}_{(q,p)})_{21}=(\pi\otimes\pi )\hat{R}_{(q,p)}(\pi^{-1}\otimes\pi^{-1})
\ ,\ \hat{R}_{(q,p)}^{-1}=(\hat{R}_{(q^{-1},p^{-1})})_{21}\ ,\eeq
where $\pi =\left(\begin{array}{cc}0&1\\1&0\end{array}\right)$.

\vskip .1cm
For the standard $R$-matrix $\hat{R}_{(q,p)}:=\hat{R}^{GL(1|1)}_{(q,p)}$,
\beq \hat{R}_{(q,p)}^t=\hat{R}_{(q,p^{-1})}\ ,\
(\hat{R}_{(q,p)})_{21}=(\pi\otimes\pi )\hat{R}_{(-q^{-1},p)}(\pi^{-1}\otimes\pi^{-1})
\ ,\ \hat{R}_{(q,p)}^{-1}=(\hat{R}_{(q^{-1},p^{-1})})_{21}\ .\eeq

\vskip .1cm
For the non-standard $GL(1|1)$ $R$-matrix $\hat{R}_{(q)}:=\hat{R}^{I}_{(q)}$,
\beq \hat{R}_{(q)}^t=\hat{R}_{(q)}\ ,\ (\hat{R}_{(q)})_{21}=\hat{R}_{(q)}
\ ,\ \hat{R}_{(q)}^{-1}=({D}\otimes{D})\hat{R}_{(q^{-1})}
({D}\otimes{D})^{-1}\ ,\eeq
where ${D} =\left(\begin{array}{cc}1&0\\0&\sqrt{-1}\end{array}\right)$.

\vskip .1cm
For the non-standard $GL(1|1)$ $R$-matrix $\hat{R}_{(q,\varepsilon )}:=
\hat{R}^{II}_{(q,\varepsilon )}$, 
\beq 
\hat{R}_{(q,\varepsilon )}^t=(\tilde{\pi}\otimes\tilde{\pi} )(\hat{R}_{(-q^{-1},-\varepsilon )})_{21}(\tilde{\pi}^{-1}\otimes\tilde{\pi}^{-1})
\ ,\ \hat{R}_{(q,\varepsilon )}^{-1}=(\hat{R}_{(q^{-1},\varepsilon )})_{21}\ ,\eeq
where $\tilde{\pi}=\left(\begin{array}{cc}0&1\\\sqrt{-1}& 0\end{array}\right)$.

\vskip .1cm
For the Jordanian $R$-matrix $\hat{R}_{(h_1:h_2)}:=\hat{R}^{(J)}_{(h_1:h_2)}$,
\beq \hat{R}_{(h_1:h_2)}^t=(\pi\otimes\pi )\hat{R}_{(h_2:h_1)}(\pi^{-1}\otimes\pi^{-1})\ ,\ 
(\hat{R}_{(h_1:h_2)})_{21}=\hat{R}_{(-h_1:-h_2)}\ ,\ 
\hat{R}_{(h_1:h_2)}^{-1}=\hat{R}_{(h_1:h_2)}\ .\eeq

\vskip .2cm
\subsection*{B.3\hspace{.3cm} Riming}\addcontentsline{toc}{subsection}{B.3\hspace{.2cm} Riming}

We shall now identify the rime $R$-matrices (\ref{rbl1})-(\ref{rbl4}).

\paragraph{\bf 1. GL(2)}$\ $

\nopagebreak
\vskip .4cm
\noindent The $R$-matrices (\ref{rbl1}) and (\ref{rbl3}) are related by a change of basis
(the number in brackets refers to the corresponding equation),
\beq \hat{R}_{(q;\gamma)}^{(\ref{rbl1})}\ T\otimes T=T\otimes T\ \hat{R}_{(q;\gamma)}^{(\ref{rbl3})}\ \ ,\ \ 
T=\left(\begin{array}{cc}q&-1/\gamma \\ \gamma &0\end{array}\right)\ .\eeq
In turn, the $R$-matrix (\ref{rbl1}) is related to the standard $R$-matrix 
$\hat{R}^{GL(2)}_{(q,q^{-1})}$ by a change of basis,
\beq\hat{R}_{(q;\gamma)}^{(\ref{rbl1})}\ T\otimes T=T\otimes T\ \hat{R}^{GL(2)}_{(q,q^{-1})}\ \ ,\ \ 
T=\left(\begin{array}{cc}q-q^{-1}&0\\\gamma &\gamma \end{array}\right)\ .\eeq

In the unitary situation (that is, $q-q^{-1}=0$), the $R$-matrix 
$\hat{R}_{(q;\gamma)}^{(\ref{rbl1})}$ belongs to the family of Jordanian $R$-matrices.
 
\vskip .2cm
Note that for the $R$-matrices (\ref{rbl1}) and (\ref{rbl3}), the left even quantum spaces are
classical.

\paragraph{\bf 2. GL(1$|$1)}$\ $

\nopagebreak
\vskip .4cm
\noindent The $R$-matrix (\ref{rbl2}) is related to the $R$-matrix (\ref{rbl4}) with the choice
$\omega =\beta_{12}^2$,
\beq \hat{R}_{(q;\gamma)}^{(\ref{rbl2})}\ T\otimes T=T\otimes T\  \hat{R}_{(-q^{-1},q^2;1)}^{(\ref{rbl4})}\ \ ,\ \ T=\left(\begin{array}{cr}1&q\\0&  \gamma q
\end{array}\right)\ .\eeq

We have
\beq\begin{array}{ccc}\label{uu1}&\hat{R}_{(q,1;\gamma)}^{(\ref{rbl4})}\ T\otimes T=T\otimes T\ 
\hat{R}^{I}_{(q)}\ \ ,\ \ T=\left(\begin{array}{cr}1&\tau\\\gamma&-\gamma \tau\end{array}\right)\ ,
\ \ {\mathrm{where}}\ \ \tau^2=\displaystyle{\frac{q-1}{q+1}}\ ,\end{array}\eeq
\beq\begin{array}{ccc}\label{uu2}&\hat{R}_{(q,q^2;\gamma)}^{(\ref{rbl4})}\ T\otimes T=T\otimes T\ 
\hat{R}^{II}_{(q,1)}\ \ ,\ \ T=\left(\begin{array}{cc}1&1\\\gamma q^{-1} & -\gamma q^{-1} \end{array}\right)\ ,\end{array}\eeq
\beq\begin{array}{ccc}\label{uu3}&\hat{R}_{(q,q^{-2};\gamma)}^{(\ref{rbl4})}\ T\otimes T=
T\otimes T\ (\hat{R}^{II}_{(q,1)})_{21}\ \ ,\ \ T=\left(\begin{array}{cc}1&1\\ \gamma q& - \gamma q\end{array}\right)\ .\end{array}\eeq

In the unitary situation (that is, for $q=\pm 1$) only eq. (\ref{uu1}) changes; but now different 
choices for $\omega$ coincide.

\paragraph{3.} Since the standard $R$-matrices are rime as well, we conclude that in dimension 2, 
all non-unitary Hecke $R$-matrices fit into the rime Ansatz. When $h_1=0$, the Jordanian 
$R$-matrix $\hat{R}^{(J)}_{(0:h_2)}$ is rime as well. However, when $h_1\neq 0$, the Jordanian 
$R$-matrix $\hat{R}^{(J)}_{(h_1:h_2)}$ cannot be rimed. Indeed, assume that $h_1\neq 0$ and let 
$A=(T\otimes T)\hat{R}^{(J)}_{(h_1:h_2)}(T\otimes T)^{-1}$ with some invertible matrix $T$. Then 
\beq\begin{array}{l}
({\mathrm{Det}}(T))^2\, A^{11}_{12}=h_1\, (T^1_1)^2\, 
({\mathrm{Det}}(T)-h_2\, T^1_1T^2_1)\ ,\\[1em]
({\mathrm{Det}}(T))^2\, A^{11}_{21}=-h_1\, (T^1_1)^2\, 
({\mathrm{Det}}(T)+h_2\, T^1_1T^2_1)\ ,\\[1em]
({\mathrm{Det}}(T))^2\, A^{22}_{12}=h_1\, (T^2_1)^2\, 
({\mathrm{Det}}(T)-h_2\, T^1_1T^2_1)\ ,\\[1em]
({\mathrm{Det}}(T))^2\, A^{11}_{21}=-h_1\, (T^2_1)^2\, ({\mathrm{Det}}(T)+h_2\, T^1_1T^2_1)\ .
\end{array}\label{nre}\eeq
For an invertible $T$, the non-rime entries (\ref{nre}) of $A$ cannot vanish 
simultaneously.

\paragraph{4.}
We remark also that all non-standard $R$-matrices of $GL(1|1)$-type are uniformly described by 
the formula (\ref{rbl4}). The right quantum spaces for the $R$-matrix 
$\hat{R}_{(q,\omega;\gamma )}^{(\ref{rbl4})}$, with $\gamma =1$, read
\begin{eqnarray}
\label{nsq1}(\hat{R}-q 1\!\! 1\otimes 1\!\!1)^{ij}_{kl}\ x^kx^l=0\ \ \ &:&
\left\{\begin{array}{l}(q+q^{-1})xy=x^2+y^2\ ,\\[.5em] 
(q+q^{-1})xy=\omega^{-1}x^2+\omega y^2\ ;\end{array}\right.\\[1em]
\label{nsq2}(\hat{R}+q^{-1} 1\!\! 1\otimes 1\!\!1)^{ij}_{kl}\ x^kx^l=0&:&
\left\{\begin{array}{l}x^2=0\ ,\\[.5em]
y^2=0\ .\end{array}\right. \end{eqnarray}
Using the diamond lemma, it is straightforward to verify that the Poincar\'e series of the quantum 
space (\ref{nsq1}) is of $GL(1|1)$-type if and only if $\omega =q^{-2},1$ or $q^2$.

\section*{Appendix C.$\ $ Rimeless triple}
\addcontentsline{toc}{section}{Appendix C.$\ $ Rimeless triple}

We sketch here a proof that the triple (\ref{ftr}) cannot be rimed. Relations 
$x^iy^j=\hat{R}^{ij}_{kl}y^kx^l$, where $\hat{R}$ is the $R$-matrix for the triple (\ref{ftr}) 
with arbitrary multiparameters, read 
\beq x^iy^i=y^ix^i\ \ ,\ \ i=1,2,3,4 \label{orr1}\eeq
and
\beq\label{orr2}\begin{array}{ll}\begin{array}{l}x^1y^2=\displaystyle{\frac{p}{q}}y^2x^1\ ,\\[.9em]
x^1y^3=\displaystyle{\frac{r}{q^2}}y^3x^1\ ,\\[.9em]
x^1y^4=\displaystyle{\frac{pr}{q}}y^4x^1-\displaystyle{\frac{rs}{q}}y^3x^2\ ,\end{array}&
\!\begin{array}{l}x^2y^1=\displaystyle{\frac{1}{pq}}\ y^1x^2+(1-q^{-2})y^2x^1\ ,\\[.9em]
x^2y^3=\displaystyle{\frac{s}{pq}}\ y^3x^2\ ,\\[.9em]
x^2y^4=\displaystyle{\frac{s}{q^2}}\ y^4x^2\ ,\end{array}\end{array}
\eeq
\beq\label{orr2p}\begin{array}{ll}
\begin{array}{l}x^3y^1=\displaystyle{\frac{1}{r}}\ y^1x^3+(1-q^{-2})y^3x^1\ ,\\[.9em]
x^3y^2=\displaystyle{\frac{p}{qs}}\ y^2x^3+(1-q^{-2})y^3x^2\ ,\\[.9em]
x^3y^4=\displaystyle{\frac{pr}{qs}}\ y^4x^3\ ,\end{array}&
\!\begin{array}{l}x^4y^1=\displaystyle{\frac{1}{pqr}}\ y^1x^4+(1-q^{-2})y^4x^1+\displaystyle{\frac{1}{q}}y^2x^3\ ,\\[.9em]
x^4y^2=\displaystyle{\frac{1}{s}}\ y^2x^4+(1-q^{-2})y^4x^2\ ,\\[.9em]
x^4y^3=\displaystyle{\frac{s}{pqr}}\ y^3x^4+(1-q^{-2})y^4x^3\ .\end{array}\end{array}\eeq
The parameter $q$ enters the characteristic equation for $\hat{R}$, $\hat{R}^2=
(1-q^{-2})\hat{R}+q^{-2}1\!\!1\otimes 1\!\!1$; $p,r$ and $s$ are the multiparameters.
The only needed restriction is $q^2\neq 1$.

\vskip .2cm
Denote by $\langle l^{(1)},l^{(2)}\rangle$ a two-dimensional plane spanned by $l^{(1)}$ and 
$l^{(2)}$. We say that two linear forms $l^{(1)}$ and $l^{(2)}$ (in four variables) form a rime 
pair if, for the ordering relations (\ref{orr1}) and (\ref{orr2})-(\ref{orr2p}), each product
$l^{(\alpha )}(x)l^{(\beta )}(y)$, $\alpha =1,2$, $\beta =1,2$, is a linear combination of 
$l^{(1)}(y)l^{(1)}(x)$, $l^{(1)}(y)l^{(2)}(x)$, $l^{(2)}(y)l^{(1)}(x)$ and $l^{(2)}(y)l^{(2)}(x)$. 
If, in addition, $l^{(\alpha )}(x)l^{(\alpha )}(y)$ is proportional to 
$l^{(\alpha )}(y)l^{(\alpha )}(x)$ for $\alpha =1$ and $2$, we say that $l^{(1)}$ and $l^{(2)}$ 
form a rime basis in the plane $\langle l^{(1)},l^{(2)}\rangle$. We call a plane {\it rime} if it 
admits a rime basis.

\vskip .4cm\noindent {\bf Fork Lemma.} {\em Assume that $l^{(1)}(x)=x^1+a_2x^2+a_3x^3$ and 
$l^{(4)}(x)=b_2x^2+b_3x^3+x^4$ form a rime pair for some $a_2,a_3,b_2$ and $b_3$. Then either
$a_3b_2\neq 0$ and $a_2=b_3=0$ or $a_2b_3\neq 0$ and $a_3=b_2=0$.

If $a_3b_2\neq 0$ then
\beq r=s=1\ ,\ l^{(1)}(x)=x^1+w x^3\ \ {\mathrm{and}}\ \  l^{(4)}(x)=x^4+\displaystyle{\frac{1}{q-q^{-1}}}
\displaystyle{\frac{1}{w}}x^2\ ,\ w\neq 0\ \ {\mathrm{is\ arbitrary}}\ .\label{btr1}\eeq
\indent If $a_2b_3\neq 0$ then
\beq p=\displaystyle{\frac{1}{q}}\ ,\ r=s\ ,\ l^{(1)}(x)=x^1+w x^2\ \ {\mathrm{and}}\ \  l^{(4)}(x)=x^4+\displaystyle{\frac{s}{q-q^{-1}}}
\displaystyle{\frac{1}{w}}x^3\ ,\  w\neq 0\ \ {\mathrm{is\ arbitrary}}\ .\label{btr2}\eeq

Moreover, if $r=s=1$ and $p\neq q^{-1}$ then the rime plane $\langle l^{(1)},l^{(4)}\rangle$ 
admits a unique, up to rescalings, rime basis $\{ l^{(1)},l^{(4)}\}$; if $p=q^{-1}$ and 
$r=s\neq 1$ then the rime plane $\langle l^{(1)},l^{(4)}\rangle$ admits a unique, up to 
rescalings, rime basis $\{ l^{(1)},l^{(4)}\}$; if $p=q^{-1}$ and $r=s=1$ then any two independent
linear combinations of $l^{(1)}$ and $l^{(4)}$ form a rime basis in the plane 
$\langle l^{(1)},l^{(4)}\rangle$.}

\vskip .2cm\noindent{\bf Proof.} A straightforward calculation.\hfill$\Box$

\vskip .2cm
Assume that a rime basis $\{ \tilde{x}^i\}$ for the triple (\ref{ftr}) exists, 
$\tilde{x}^i=A^i_jx^j$, the matrix $A^i_j$ is invertible. Rename the rime  variables $\tilde{x}^i$ 
in such a way that the minor $\left|\begin{array}{cc}A^1_1&A^1_4\\ A^4_1&A^4_4\end{array}\right|$ 
is non-zero and $A^1_1A^4_4\neq 0$; normalize the variables $\tilde{x}^1$ and $\tilde{x}^4$ to 
have $A^1_1=A^4_4=1$. The plane $\langle\tilde{x}^1,\tilde{x}^4\rangle$ is, by definition, rime, 
with a rime basis $\{ \tilde{x}^1,\tilde{x}^4\}$.

Suppose that $r=s=1$ and $p\neq q^{-1}$ or $p=q^{-1}$ and $r=s\neq 1$. Then, by Fork Lemma, the
rime basis in the plane $\langle\tilde{x}^1,\tilde{x}^4\rangle$ is, up to proportionality,
unique, so we know the variables $\tilde{x}^1$ and $\tilde{x}^4$. The variables
$\tilde{x}^1$ and $\tilde{x}^2$ form a rime plane. Therefore, if the variable $\tilde{x}^2$ 
contains $x^4$ with a non-zero coefficient then, by Fork Lemma, $\tilde{x}^2$ must
be proportional to $\tilde{x}^4$, contradicting to the linear independence of the variables 
$\tilde{x}^2$ and $\tilde{x}^4$. Similarly, the variable $\tilde{x}^2$ cannot contain
$x^1$ with a non-zero coefficient (the plane $\langle\tilde{x}^2,\tilde{x}^4\rangle$ is rime).
Thus, $\tilde{x}^2$ is a linear combination of $x^2$ and $x^3$. Same for $\tilde{x}^3$: it is a 
linear combination of $x^2$ and $x^3$. One of the variables, $\tilde{x}^2$ or $\tilde{x}^3$, say, $\tilde{x}^2$, contains $x^2$ with a non-zero coefficient. Writing rime equations for the plane
$\langle\tilde{x}^1,\tilde{x}^2\rangle$ in the case $r=s=1$ and $p\neq q^{-1}$ (for the plane 
$\langle\tilde{x}^2,\tilde{x}^4\rangle$ in the case $p=q^{-1}$ and $r=s\neq 1$) 
quickly leads to a contradiction.

\vskip .2cm
Therefore, if the relations (\ref{orr1}) and (\ref{orr2})-(\ref{orr2}) can be rimed then $p=q^{-1}$ and $r=s=1$.
It follows from Fork Lemma that $$\tilde{x}^4=(q-q^{-1})c_2c_3x^1+c_2x^2+c_3x^3+x^4$$
for some $c_2$ and $c_3$. The planes $\langle\tilde{x}^a,\tilde{x}^4\rangle$, $a=1,2,3$, are rime.
Subtracting from the variables $\tilde{x}^a$ the variable $\tilde{x}^4$ with appropriate
coefficients, we find three linearly independent combinations 
\beq l(x)=d_1x^1+d_2x^2+d_3x^3\ ,\label{r3c}\eeq 
each forming a rime pair with $\tilde{x}^4$. We must have: $l(x)l(y)$ is a linear combination of
$l(y)l(x)$, $l(y)\tilde{x}^4$, $\tilde{y}^4 l(x)$ and $\tilde{y}^4\tilde{x}^4$. It follows, 
after a straightforward calculation, that $d_2d_3=0$. Moreover, $d_2=d_3=0$ is excluded by Fork
Lemma. In the case $d_2\neq 0$ and $d_3=0$ (respectively, $d_3\neq 0$ and $d_2=0$), the rime
condition implies that $d_1=(q-q^{-1})c_2d_3$ (respectively, $d_1=(q-q^{-1})c_3d_2$). Thus,
only two linearly independent combinations (\ref{r3c}) can form a rime pair with $\tilde{x}^4$, 
the final contradiction.

\end{document}